\documentclass{conm-p-l}
\usepackage{amssymb,amsmath,amsthm,amscd} 
\usepackage{hyperref}

\hypersetup{
  colorlinks   = true, 
  urlcolor     = blue, 
  linkcolor    = blue, 
  citecolor   = red 
}

\usepackage[sorted, msc-links]{amsrefs}

\usepackage{mathtools}

\usepackage{csquotes}

\usepackage[capitalise]{cleveref}

\crefname{thm}{Thm.}{}
\crefname{prop}{Prop.}{}
\crefname{lem}{Lem.}{}
\crefname{cor}{Cor.}{}

\usepackage[linesnumbered,lined, ruled]{algorithm2e}

\MakeOuterQuote{"}

\newtheorem{thm}{Theorem}
\newtheorem{prop}{Proposition}
\newtheorem{lem}{Lemma}
\newtheorem{cor}{Corollary}
\newtheorem{rem}{Remark}
\newtheorem{exa}{Example}

\newtheorem{problem}{Problem}
\newtheorem{quest}{Question}
\newtheorem{corollary}{Corollary}

\theoremstyle{definition}
\newtheorem{defi}[thm]{Definition}

\def\P{\mathbb P}
\def\H{\mathbb H}
\def\O{\mathcal O}
\def\<{\langle}
\def\>{\rangle}

\def\a{\alpha}

\def\l{\lambda}

\def\CC{\mathcal{C}}
\def\AA{\mathcal A}
\def\t{\tau}

\DeclareMathOperator\Aut{Aut }
\DeclareMathOperator\Mat{Mat }
\DeclareMathOperator\Img{Img }
\DeclareMathOperator\Jac{Jac }  
\DeclareMathOperator\Gal{Gal }
\DeclareMathOperator\Hom{Hom }
\DeclareMathOperator\End{End }

\DeclareMathOperator\tr{tr }

\DeclareMathOperator\ord{ord }

\DeclareMathOperator\Pic{Pic}
\DeclareMathOperator\Div{Div}
\DeclareMathOperator\PDiv{PDiv}

\DeclareMathOperator\SL{SL}

\DeclareMathOperator\Spec{Spec}
\DeclareMathOperator\conorm{conorm}

\newcommand\N{\mathbb N}
\newcommand\Z{\mathbb Z}
\newcommand\Q{\mathbb Q}
\newcommand\R{\mathbb R}
\newcommand\C{\mathbb C}
\newcommand\A{\mathbb A}
\newcommand\F{\mathbb F}

\newcommand\iso{{\, \cong\, }}

\newcommand\M{\mathcal M}

\newcommand\p{\mathfrak p}

\newcommand\X{\mathcal X}
\newcommand\Y{\mathcal Y}

\newcommand\om{\omega}

\newcommand\G{\Gamma}

\newcommand\E{\mathcal{E}}

\newcommand\J{\mathcal{J}}
\newcommand\Fp{{\F}_p}

\newcommand\Fq{\mathbb{F}_q}

\newcommand\f{\mathfrak f}

\newcommand\B{\mathcal B}
\newcommand\K{\mathcal K}
\def\char{\mbox{char }}

\newcommand\ra{\rightarrow}

\newcommand\car{\rm char}

\newcommand\PP{\frak{P}}

\newcommand\g{\gamma}

\DeclarePairedDelimiter\floor{\lfloor}{\rfloor}

\newcommand\ov{\overline}
\newcommand\cS{\mathcal{S}}

\newcommand\Mor{\mathrm{Mor}}

\newcommand\Ss{\mathcal{S}}

\def\L{\mathcal L}
\newcommand\D{\mathcal{D}}

\newcommand\WP{{\mathbb{WP}}^3_{(2, 4, 6, 10)}}
\newcommand\x{\mathbf x}

\newcommand\T{\theta}

\newcommand{\ch}[2]
{\begin{bmatrix}
 #1 \\
 #2\\
\end{bmatrix}}

\newcommand{\chr}[4]
{\begin{bmatrix}
 #1 & #2\\
 #3 & #4
\end{bmatrix}}

\newcommand\m{\mathfrak m}
\newcommand\bn{\mathfrak b}
\newcommand\embd{\hookrightarrow}

\begin{document}
%
%

\title{Curves, Jacobians, and Cryptography}

\author{Gerhard Frey}

\address{Institut für Experimentelle Mathematik, \\
Universit\"at Duisburg-Essen, 45326 Essen, Germany.}
\email{frey@iem.uni-due.de}

\author{Tony Shaska}
\address{Department of Mathematics and Statistics, \\
Oakland University, Rochester, MI 48309, USA.}
\email{shaska@oakland.edu}

\begin{abstract}
The main purpose of this paper is to give an overview over the   theory  of abelian varieties,
with main focus on Jacobian varieties of curves   reaching from well-known results till to latest developments  and their usage in cryptography. 
In the first part we provide the necessary mathematical background  on abelian varieties, their torsion points, Honda-Tate theory, Galois representations, with emphasis on Jacobian varieties and hyperelliptic Jacobians. 
In the second part we focus on applications of abelian varieties on cryptography and treating separately, elliptic curve cryptography,   genus 2 and 3 cryptography, including Diffie-Hellman Key Exchange, index calculus in Picard groups, isogenies of Jacobians via correspondences and applications to discrete logarithms. Several open problems and new directions are suggested.  
\end{abstract}


\maketitle
 
\section*{Preface}

There has been a continued interest on Abelian varieties in mathematics during the last century.  Such interest is  renewed   in the last few years, mostly due to   applications of abelian varieties in cryptography.  In these notes we  give a brief introduction to the mathematical background on abelian varieties and their applications on cryptography with the twofold aim of introducing abelian varieties to the experts in cryptography and introducing methods of cryptography to the mathematicians working in algebraic geometry and related areas. 
%
\subsection*{A word about cryptography}
Information security will continue to be one of the greatest challenges of the modern world with implications in technology, politics, economy, and every aspect of everyday life. Developments and drawbacks of the last decade in the area will continue to put emphasis on searching for safer and more efficient crypto-systems. The idea and lure of the quantum computer makes things more exciting, but at the same time frightening. 

There are two main methods to achieve secure transmission of information: \textit{secret-key cryptography} (\textit{symmetric-key}) and \textit{public-key cryptography} (\textit{asymmetric-key}).  The main disadvantage of symmetric-key cryptography is that a shared key must be exchanged beforehand in a secure way.  In addition, managing keys in a large public network becomes a very complex matter.  Public-key cryptography is used as a complement to secret-key cryptography for signatures, authentication and key-exchange. There are two main methods used in public-key cryptography, namely RSA and the discrete logarithm problem (DLP) in cyclic groups of prime order which are embedded  in rational points of Abelian varieties, in particular of Jacobian varieties  of curves.
 The last method is usually referred to as \textit{curve-based cryptography}. 

In addition, there is always the concern about the post-quantum world. What will be the crypto-systems which can resist the quantum algorithms? Should we develop such systems now? There is enthusiasm in the last decade that some aspects of curve-based cryptography can be adapted successfully to the post-quantum world.  Supersingular Isogeny Diffie-Hellman (SIDH), for example,  is based on isogenies of supersingular elliptic curves  and is one of the promising schemes for post-quantum cryptography. Isogenies of hyperelliptic Jacobians of dimension 2 or 3 have also been studied extensively in the last decade and a lot of progress has been made.  In this paper we   give an overview of recent developments in these topics.   

\subsection*{Audience}  Computer security and cryptography courses for mathematics and computer science majors are being introduced in all major universities. Curve-based cryptography has become a big part of such courses and a popular area even among professional mathematicians who want to get involved in cryptography.   The main difficulty that these newcomers is the advanced mathematical background needed to be introduced to curve-based cryptography. 

Our target audience is advanced graduate students and researchers from mathematics or computer science departments who work with curve-based cryptography.  Many researchers from other areas of mathematics who want to learn about abelian varieties and their use in cryptography will find these notes useful.

\subsection*{Notations and bibliography}

The symbols $\N$ and $\Z$ will denote the natural numbers and the ring of integers while $\Q$, $\R$, and $\C$   the fields of rationals, reals, and complex numbers. $\F_p = \Z/p\Z$ will denote the field of $p$ elements, for a prime number $p$ and $\F_q$ a field of $q$ elements, where $q=p^n$ is a power of $p$. 

In general a field will be denoted by $k$.  We shall always assume that $k$ is a perfect field, i.e., every algebraic extension of $k$ is separable.  Denote  its characteristic by $\char k$, and its algebraic closure by $\bar k$. The Galois group $\Gal (\bar k / k)$ will be denoted by $G_k$. A number field will be denoted by $K$ and its ring of integers by $\O_K$. 

By a  "curve" we mean an irreducible, smooth, projective algebraic curve.  A genus $g \geq 2$ curve defined over $k$ will be denoted by $\CC/k$ or sometimes   by  $\CC_g$ and its Jacobian by $\Jac _\CC$. The automorphism group of $\CC$  is denoted by $\Aut  \CC$ and it means the full group of automorphisms of $\CC$ over the algebraic closure $\bar k$. 

We will use $\AA$, $\B$ to denoted Abelian varieties defined over a field $k$ and $k (\AA)$, $k(\B)$ their function fields.  The set of $n$-torsion points of $\AA$ will be denoted by $\AA[n]$.  For a subscheme $G$ of $\AA$, the quotient variety is denoted by $\AA/G$. 

\subsection*{Background and preliminaries}

We assume the reader is familiar with the basic tools from algebra and algebraic geometry. Familiarity with algebraic curves   is expected and the ability to read some of the classical works on the subject \cite{Mum}. 
\subsection*{Organization of these notes}
In \cref{part-1} we give the mathematical background on    Abelian varieties, their torsion points, endomorphisms and isogenies. We focus mostly on Abelian varieties defined over    fields of positive characteristic. The main references here are  \cite{Mum}, and \cite{Frey}.  

We give a brief introduction of abelian varieties as complex tori with period matrices in \cref{sect-1}.  Of course, of special interest for us are Jacobian varieties, hence we define in detail algebraic curves, constant field extensions, group schemes, and principally polarized varieties.

In \cref{sect-2}   we focus on endomorphism rings of abelian varieties and isogenies,  the characteristic polynomial of the Frobenius, $l$-adic Tate module, and Tate's result on determining necessary and sufficient conditions for two Abelian varieties to be isogenous.  

Algebraic curves and Jacobian varieties are treated in detail in  \cref{jacobians}, including   Picard groups on curves, the group of divisors, canonical divisors, Riemann-Roch theorem,  and the definition of   Jacobians of curves.   

In \cref{sect-4} we focus on applications of the Riemann-Roch theorem, including the Hurwitz genus formula, gonality of curves and Hurwitz spaces, Cantor's algorithm on Jacobians of hyperelliptic curves, and automorphism groups of curves and their Jacobians.  As illustration and for later applications curves of small genus are discussed in more detail.

In \cref{sect-5} we give a brief description of the theory of modular curves over $\C$, modular polynomials, and the arithmetic theory of modular curves. 

In \cref{part-2} we focus on applications of abelian varieties on cryptography. Our main reference is \cite{book} and the material provided in \cref{part-1}. 

In \cref{Diffie-Hellman}  we go over the preliminaries of the Diffie-Hellman Key Exchange and the mathematical challenges including $Q$-bit security, Key Exchange with $G$-sets, and the abstract setting of Key Exchange. 

In \cref{index-calculus} we describe the methods of index calculus in Picard groups and their use in cryptography. Such methods have been quite successful due to work of Diem, Gaudry, et al.  As consequence one sees that only elliptic and hyperelliptic curves  of genus $\leq 3$ provide candidates  for secure crypto systems based on discrete logarithms. Hence we shall discuss these curves in detail. In \cref{sect-correspondences} we focus on isogenies of Jacobians via correspondences.  We discuss the Weil descent, modular correspondences, and correspondences via monodromy groups. It is an open and difficult problem to find interesting correspondences of low degree between Jacobian varieties induced by correspondences between curves. 

In  \cref{gen-3} we study hyperelliptic Jacobians of dimension 3. We give a short introduction of non-hyperelliptic and hyperelliptic genus 3 curves and their plane equations. Then we define Picard groups of genus 3 curves and their use in cryptography and results of Diem and Hess.  In the following part we describe the index-calculus attacks applied to  genus 3 and results of Diem, Gaudry, Thom\'{e}, Th\'{e}riault. We also discuss isogenies via $S_4$-covers and work of Frey and Kani \cite{FK1}, \cite{FK2}  and Smith \cite{Smith}. 

In  \cref{gen-2} we focus on dimension 2 Jacobians and their use in cryptography.   Methods based on \cite{lombardo} of how to compute the endomorphism ring of a dimension 2 Jacobian are described and in particular isogenies of Abelian surfaces via Donagi-Livn\'{e} approach and some recent results of Smith \cite{Smith}. Further we give details of point counting algorithms on genus 2 Jacobians and explicit formulas for $[n]D$, when $D$ is a reduced divisor. Work of Gaudry, Harley, Schost and others is briefly described. 
In  \cref{ell-curves} we focus on the elliptic curves and elliptic curve cryptography. We  give an explicit description of the methods used in supersingular isogeny-based cryptography.   We describe the necessary background including Velu's formula,  ordinary and supersingular elliptic curves and the more recent results \cite{de-feo-1}, \cite{de-feo-2}, \cite{DJP} among others. 

\part{Abelian varieties}\label{part-1}

In the first  part of these notes we give the basic theory of abelian varieties, their endomorphisms, torsion points, characteristic polynomial of the Frobenius, Tate models, and then focus on Jacobian varieties and hyperelliptic Jacobians. While there are many good references on the topic, we mostly use \cite{Mum},   \cite{Tate}. 
\section{Definitions and basic properties}\label{sect-1}
We shall use projective respectively affine \emph{schemes} defined over $k$. Let $n\in \N$ and $I_h$ (respectively $I$) be a homogeneous   ideal in $k[Y_0,\cdots, Y_n]$ different from  $\<Y_0,\dots Y_n \>$  (respectively an arbitrary ideal in $k[X_1,\dots,X_n]$).    

Let $R_h:=k[Y_0,\dots,Y_n]/I_h$ (respectively $R:=k[X_1,\dots,X_n]/I$) be the quotients.    By assumption, $R_h$ is a graded ring, and so localizations $R_{h,\mathfrak{A}}$ with respect to homogeneous ideals $\mathfrak{A}$ are graded, too. Let $R_{h,\mathfrak{A}}$ $_0$ be the ring of elements of grade $0$.
 
The projective scheme $\Ss_h$ (respectively the affine scheme $\Ss$) defined by $I_h$ ($I$) consists of
\begin{enumerate}
\item the  topological space $V_h:=\mathrm{Proj}(R_h)$ ($V:=\mathrm{Spec(R)}$) consisting of homogeneous prime ideals in $R_h$ with pre-image in $k[Y_0,\dots,Y_n]$ different from $\<Y_0,\dots Y_n \>$ (prime ideals in $R$) endowed with the Zariski topology and

\item the sheaf of rings of holomorphic functions given on Zariski-open sets $U\subset V_h$ ($U\subset V$) as elements of grade $0$ in localization  of $R_{h,0}$ ($R$) with respect to the elements that become invertible when restricted to $U$. 
\end{enumerate}
 
\noindent \textbf{Examples:} 

\begin{enumerate}
\item The projective space $\P^n$ over $k$ of dimension $n$ is given by the ideal $\<0\>\subset k[Y_0,\dots,Y_n]$.  The ring of holomorphic functions on $\P^n$ (take $U=\P^n$) is $k$.

Next take $U=\emptyset$ to get the ring of \emph{meromorphic} functions on $\P^n$:
It consists of the quotients 
\[ 
f/g\mbox{ with  }f, g \mbox{  homogeneous of degree d with }g\neq 0.
\]

\item The affine space $\A^n$ of dimension $n$ over $k$ is the topological space 
 \[\Spec(k[X_1,\dots X_n]).\]
The ring of holomorphic functions on $A^n$ is $k[X_1,\dots,X_n]$, where polynomials are interpreted as polynomial functions.   The ring of meromorphic functions on $\A^n$ (take $U=\emptyset$) is the field of rational functions $k(X_1,\dots ,X_n)$.

\item The easiest but important example for an affine scheme: Take $n=1$, $I= \<X_1 \>$, $V=\Spec (k)=\{(0)\}$ and  $O_{(0)}=k^*$. 
\end{enumerate}

\emph{Morphisms} of affine or projective schemes are continuous maps between the underlying topological spaces induced (locally) by (in the projective case, quotients of the same degree) of polynomial maps of the sheaves.

\emph{Rational maps} $f$ between affine or projective schemes $\Ss$ and $\mathcal{T}$ are equivalence classes of morphisms defined on open subschemes $U_i$ of $\Ss$ with image in $\mathcal{T}$ and compatible with restrictions to $U_i\cap U_j$.   If $f$ is invertible (as rational maps from $\mathcal{T}$ to $\Ss$), then $f$ is \emph{birational}, and $\Ss$ and $\mathcal{T}$ are birationally equivalent.

The $k$-rational points $\Ss(k)$ of a scheme $\Ss$ is the set of morphisms from $\Spec (k)$ to $\Ss$. The reader should verify that  for projective schemes defined by the ideal $I_h$ the set $\mathcal{S}(k)$ is, in a natural way, identified with  points $(y_0:y_1\dots  :y_n) $ with $k$-rational homogeneous coordinates in the projective space of dimension $n$ which are common zeros of the polynomials in $I_h$, and an analogous statement holds for affine schemes. \\

\noindent \textbf{Constant field extensions:}  \,  Let $k\stackrel{\iota}{\hookrightarrow} L$ be an embedding of $k$ into a field $L$ (or $k\subset L$ if the embedding is clear) of $k$. Let $\Ss$ be a projective (affine) scheme defined over $k$ with ring $R$. $\iota$ induces a  morphism $\f_\iota$ from $R$ in $R\otimes_{k}L=:R_\iota$  given by the interpretation via $\iota$ of polynomials  with coefficients in $k$ as polynomials with coefficients in $L$. The   ideal $I_\Ss$ extends to a  ideal  in  $R_\iota$ and hence we get in a natural way a projective  scheme  $\Ss_\iota$ with a morphism
\[ 
\Ss_\iota\ra \Ss
\]
as $\Spec (k)$ schemes. $\Ss_\iota $ is again a  projective (affine) scheme now defined over $L$, which is denoted as  \emph{scalar extension} by $\iota$.  If there is no confusion possible (for instance if $k\subset L\subset \bar{k}$ and $\iota$ is the inclusion) we denote $\Ss_\iota$ by $\Ss_L$. 

A scheme $\Ss$ is irreducible if the ideal $I_h$ (respectively $I$) is a prime ideal. $\Ss$ is absolutely irreducible if $\Ss_{\bar{k}}$ is irreducible. This is the case if and only if $k$ is algebraically closed in $R$. Classically, irreducible schemes are called \emph{irreducible varieties}.   \\

\noindent \textbf{Affine covers} \,   There are many possibilities to embed $\A^n$ into $\P^n$, and there is no "canonical" way to do this. But after having chosen coordinates there is a standard way to construct a covering of $\P^n$ by $n+1 $ copies of $\A^n$. Every homogeneous polynomial $P(Y_0,...,Y_n)$ can be transformed into $n+1$ polynomials $p_j(X)$ ($j=0,...,n$) in $n$ variables by the transformation  
\[t_j: Y_i\mapsto X_i:=Y_i/Y_j.\] 
Define $U_j$ as open subscheme of $\P^n$ which is the complement of the projective scheme attached to the ideal $\<Y_j\>$. Then $t_j$ $_{|U_j}$ is holomorphic and bijective and its image is isomorphic to $\A^n$.

By the inverse transform $\iota_j$ we embed $\A^n$ into $\P^n$ and so $U_j$ is isomorphic to $\A^n$ as affine variety.  Taking the collection $(\iota_0,\dots \iota_n)$ we get a finite open covering of $\P^n$  by $n+1$ affine subspaces.

Having an affine cover $U_j$ of $\P^n$ one can intersect it with projective varieties $V$ and get
\[V=\bigcup _j V_{j,a}\;  \; \text{with }   \;   V_{j,a}:=V\cap U_j\]
as union of affine varieties. \newline

\noindent \emph{Converse process:} \,  Given a polynomial $p(X_1,...,X_n)$ of degree $d$ we get a homogeneous polynomial $p^h(Y_0,...,Y_n)$ of degree $d$ by the transformation 
\[ X_i\mapsto Y_i/Y_0\mbox{ for }i=1, \dots , n \]
and then clearing denominators.  Assume that $V_a$ is an affine variety with ideal $I_a\subset k[X_1, \dots ,X_n]$. By applying the homogenization explained above to all polynomials in $I_a$ we get a homogeneous ideal $I^h_a\subset k[Y_0,...,Y_n]$ and a projective variety $V$ with ideal $I^h_a$ containing $V_a$ in a natural  way.   $V$ is called a projective closure of $V_a$.    By abuse of language one  calls $V\cap U_0=V\setminus V_a$ "infinite points" of $V_a$. \\

\noindent \textbf{Function Fields:} \,  Let $\Ss\subset \A^n$ be an affine irreducible variety with ring $R$. In particular, $R$ is an integral domain. The function field $k(\Ss) $ is the quotient field of $R$. It consists of the meromorphic functions  of $\A^n$ restricted to $\Ss$. $\mathcal{T}$ is birational equivalent to $\Ss$ if and only if $k(\Ss)=k(\mathcal{T})$.

If $U \neq \emptyset$ is affine and open in a projective variety $\Ss$ then the field of meromorphic functions $k(\Ss) = k(U)$. In particular, it is independent of   $U$.

\begin{defi}
Let $\Ss$ be an irreducible variety.  The dimension of $\Ss$ is the transcendental degree of $k(\Ss)$ over $k$.
\end{defi}

\noindent \textbf{Group schemes:} \,   A projective (affine) group scheme $G$ defined over $k$ is a projective (affine) scheme  over $k$
endowed with

\par i) addition, i.e., a morphism 
\[ m   : \; G \times G \to G \]
\par ii) inverse, i.e., a morphism 
\[ i  : \;  G \to G \]
\par iii) the identity, i. e., a $k$-rational point $0 \in G(k)$,  \\

\noindent such that it satisfies group laws. The group law is uniquely determined by the choice of the identity element.  A morphism of group schemes that is compatible with the addition law is a homomorphism.

Let $L$ be a field extension of $k$.  $G (L)$ denotes the set of $L$-rational points of $G$ and it is also a group. A homomorphism between groups schemes induces a homomorphism between the group of rational points.   If $G$ is an absolutely irreducible projective variety, then the group law $m$ is commutative. 

\begin{defi} 
An Abelian variety defined over $k$ is an absolutely irreducible projective variety defined over $k$ which is  a group scheme.
\end{defi}

We will denote an Abelian variety  defined over a field $k$ by $\AA_k$ or simply $\AA$ when there is no confusion. From now on the addition $m(P, Q)$ in an abelian variety will be denoted by $P\oplus Q$ or simply $P+Q$ and the inversion $i(P)$ by $\ominus P$ or simply by $-P$.  \\

\noindent \textbf{Fact:} A morphism from the  Abelian varieties $\AA_1$ to the Abelian variety $\AA_2$ is a homomorphism if and only if it maps the identity element of $\AA_1$ to the identity element of $\AA_2$.  

An abelian variety over a field $k$ is called \textbf{simple} if it has no   proper non-zero Abelian subvariety over $k$,  it is called \textbf{absolutely simple} (or \textbf{geometrically simple}) if it   is simple over the algebraic closure of $k$. 

\subsection{Complex tori and abelian varieties}
Though we are interested in Abelian varieties over  arbitrary fields $k$ or in particular, over finite fields, it is helpful to look at the origin of the whole theory, namely the theory of Abelian varieties over the complex numbers. Abelian varieties are connected, projective algebraic  group schemes.   Their analytic counterparts are the connected compact Lie groups.  
 
Let $d$ be a positive integer and $\C^d$ the complex Lie group (i.e., with vector addition as group composition). The group $\C^d$ is not compact, but we can find quotients which are compact. Choose a lattice $\Lambda \subset \C^d$ which is a $\Z$-submodule of rank $2d$.   The quotient $\C^d/\Lambda$ is a complex, connected Lie group which is called a \textit{complex $d$-dimensional torus}. Every connected, compact Lie group of dimension $d$ is isomorphic to a torus $\C^d/\Lambda$. 

A Hermitian form $H$ on $\C^d \times \C^d$ is a form that can be decomposed as 
\[ H(x, y) = E(ix, y) + i \, E(x, y), \]
where $E$ is a skew symmetric real form on $\C^d$ satisfying $E(ix, iy)=E(x, y)$. $E$ is called the imaginary part  $\Img (H)$ of $H$. 
The torus $\C^d/\Lambda$ can be embedded into a projective space if and only if there exists a positive Hermitian form $H$ on $\C^d$ with $E=\Img (H)$ such that restricted to $\Lambda \times \Lambda$ has values in $\Z$.  Let  $\mathbb H_g$ be the Siegel upper half plane 
\[\mathbb H_d =\{\tau \in \Mat_d  (\C) \ | \ \tau^T=\tau, \,  \Img (\tau)>0 \}.\] 
Then, we have the following.

\begin{lem}
Let $\C^d/\Lambda$ be a complex torus attached to an abelian variety $\AA$. Then $\Lambda$ is isomorphic to $\Z^d \oplus \Omega \cdot \Z^d$, where $\Omega \in \mathbb{H}_d$.
\end{lem}

The matrix $\Omega$ is called the \textbf{period matrix} of $\AA$.   The lattice $\hat \Lambda$ given by 
\[   \hat \Lambda := \{  x\in \C^d \, | \, E(x, y)\in \Z, \; \; \text{for all } \; y\in \Lambda    \} \]
is called the \textbf{dual lattice} of $\Lambda$.   If $\hat \Lambda = \Lambda$ then $E$ is called a \textit{principal polarization} on $\AA$ and the pair $(\AA, E)$ is  called a \textbf{principally polarized} abelian variety; we may also say that $\AA$ admits a principal polarization. 

For a principally polarized abelian variety $(\AA, E)$ there exists a basis $\{ \mu_1, \dots , \mu_{2d} \}$ of $\Lambda$ such that 
\[
J:= \left[ E(\mu_i, \mu_j)   \right]_{1 \leq i, j \leq 2d} \, = \, 
\begin{bmatrix}
  0 & I_d \\
  -I_d & 0
\end{bmatrix}.
\]
The symplectic group  
\[ 
Sp(2d, \Z) = \{ M\in GL(2d,\Z) \ | \ MJM^T=J\}\]     acts on  $\mathbb{H}_d$, via 
\[
\begin{split}
 Sp   (2d, \Z) \times \H_d  & \rightarrow \H_d   \\
\begin{bmatrix} a & b \\ c & d \end{bmatrix} \times \tau & \rightarrow (a \tau + b)(c \tau + d)^{-1}
\end{split}
\]
where $a, b, c, d, \tau$ are $d \times d$ matrices. The moduli space of $d$-dimensional abelian varieties is 
\[ \mathbf A_g := \mathbb H_d / Sp(2d, \Z).\]
Jacobian  varieties of  projective irreducible nonsingular curves admit  canonical principal polarizations.  
These Abelian varieties are in the center of our interest and will be discussed in detail  in \cref{jacobians}. 
%
\subsubsection{Elliptic curves over $\C$}
Take $d=1$ and a lattice $\Lambda_\tau:= \Z+\Z\tau$ with $\tau$ in the upper half plane $\H_1=\H$. The  torus $\C/\Lambda_\tau$ is a compact Riemann surface and so an algebraic projective curve $\E_\tau$ over $\C$.

The function field $\C(\E_\tau)$ is generated by the   Weierstrass function $\wp _\tau$ and its derivative $\wp'_\tau$, which are meromorphic functions on $\C$ with periods $1, \tau$ and poles of order $2$ respectively $3$ in $\Lambda_\tau.$  $\wp_\tau$  satisfies a differential equation $$W_\tau: \wp_\tau'^2=\wp_\tau^3-g_2(\tau)\wp_\tau^2-g_3(\tau).$$ This is an affine equation for $\E_\tau$, by introducing homogeneous coordinates $(X:Y:Z)$ by  $\wp_\tau=X/Z,\wp'_\tau=Y/Z, $ we get the projective plane curve $\E_\tau$ with equation 
\[Y^2Z=X^3-g_2(\tau)XZ^2-g_3(\tau)Z^3\]
with coefficients $g_2(\tau),g_3(\tau)$ depending on $\Lambda_\tau$ in a very specific way: $g_2$ and $g_3$ are Eisenstein series in $\tau$.  It follows that $\Delta_\tau=4g_2(\tau)^3-27  g_3(\tau)^2\neq 0$ and so $\E$ is without singularities.  We get a parametrization
\[\phi:\C\ra \E_\tau(\C)\]
by 
\[z\mapsto (\wp_\tau(z):\wp'_\tau(z):1) \; \;     \mbox{ if and only if }  z\notin \Lambda_\tau \]
and $\phi(\Lambda_\tau)=(0:1:0)$, the point at infinity. This parametrization yields on $\E_\tau$ an addition and makes $\E_\tau$ to an Abelian variety of dimension $1$ over $\C$.

\begin{defi}
An Abelian variety of dimension $1$ is called an \textbf{elliptic curve}.
\end{defi}

We have seen that we can attach to every elliptic curve $\E$ an element $\tau\in \mathbb{H}$ such that $\E$ is isomorphic to $\E_\tau$. Let $\E_{\tau'}$  be another  elliptic curve. Then $\E_\tau$ is isomorphic to $\E_{\tau'}$ if and only if $\tau$ is equivalent to $\tau'$ under the action of $Sp(2d, \Z)=Sl(2,\Z)$ on $\mathbb{H}$.

Since $\H/Sl(2,\Z)$ is as Riemann surface isomorphic to $\mathbb{A}^1$ we get a one-to-one correspondence between  isomorphism classes of elliptic curves over $\C$ and  points on the affine line.  This correspondence is given by a modular function (i.e. a holomorphic function on $\mathbb{H}$ invariant under $Sl(2,\Z)$): the $j$-function.

\begin{defi}
The absolute invariant of $\E_\tau$ is given  by $j(\tau):=12^3\frac{4g_2^3(\tau)}{\Delta_\tau}$.
\end{defi}

\begin{thm}
$\E_\tau$ is isomorphic to $\E_{\tau'}$ if and only it $j(\tau)=j(\tau')$.   Hence the $j$-function is an analytic map from $\A^1$ to $\mathbb{A}^1$.
\end{thm}

We remark that we shall define elliptic curves $\E$ in a purely algebraic setting over arbitrary fields $k$ (cf. \cref{section-4-2-3}) and that we shall define an absolute invariant $j$ for such curves, which coincides with $j(\tau)$ if $k=\C$, and which also  has the property: If $\E$ is isomorphic to $\E'$ then $j_\E=j_\E'$, and the converse holds if $k$ is algebraically closed.

\section{Endomorphisms and isogenies}\label{sect-2}

Let $\AA$, $\B$ be abelian varieties over a  field $k$.  We denote the $\Z$-module of homomorphisms  $\AA \mapsto  \B$  by $\Hom( \AA, \B)$  and the ring of endomorphisms $\AA \mapsto \AA$ by $\End \AA$. In the context of Linear Algebra it can be more convenient  to   work with the $\Q$-vector spaces $\Hom^0 (\AA, \B):= \Hom(\AA, \B) \otimes_\Z \Q$, and $\End^0 \AA:= \End \AA\otimes_\Z \Q$. Determining $\End \AA$ or $\End^0 \AA$ is an interesting problem on its own; see  \cite{Oort}.

For any abelian variety $\AA$ defined over a  number field $K$, computing $\End_K (\AA)$ is a harder problem than computation of $\End_{\bar K} (\AA)$; see \cite[lemma~5.1]{lombardo} for details.

\subsection{Isogenies}
A homomorphism $f : \AA \to  \B$ is called an \textbf{isogeny} if $\Img f = \B$  and $\ker f$ is a finite group scheme. If an isogeny $\AA \to \B$ exists we say that $\AA$ and $\B$ are isogenous.  We remark that this relation is symmetric, see  \cref{dual}.

The degree of an isogeny $f : \AA \to \B$ is the degree of the function field extension
\[ \deg f  := [k(\AA) : f^\star k(\B)].\]
It is equal to the order of the group scheme $\ker (f)$, which is, by definition, the scheme theoretical inverse image $f^{-1}(\{0_\AA\})$.

The group of $\bar{k}$-rational points has order $\#(\ker f)(\bar{k}) = [k(A) : f^\star k(B)]^{sep}$, where $[k(A) : f^\star k(B)]^{sep}$ is the degree of the maximally separable extension in $k(\AA)/ f^\star k(\B)$.   $f$ is a \textbf{separable isogeny} if and only if 
\[ \# \ker f(\bar{k}) = \deg f.\]
Equivalently: The group scheme $\ker f$ is \'{e}tale. 
The following result   should be compared with the well known result of quotient groups of abelian groups.

\begin{lem}\label{noether} 
For any Abelian variety $\AA/k$ there is a one to one correspondence between the finite subgroup schemes $\K \leq \AA$ and   isogenies $f : \AA \to \B$, where $\B$ is determined up to isomorphism.   Moreover, $\K = \ker f$ and $\B = \AA/\K$. 

$f$ is separable if and only if $\K$ is   \'{e}tale, and then $\deg f =\#\K(\bar{k})$.
\end{lem}

\noindent Isogenous Abelian varieties have  commensurable endomorphism rings.

\begin{lem} 
If  $\AA$ and $\B$ are isogenous then $\End^0 (\AA) \iso \End^0 (\B)$. 
\end{lem}

\begin{lem}
If $\AA$ is a  absolutely simple Abelian variety then every endomorphism  not equal $0$ is an isogeny. 
\end{lem}

We can assume that $k=\bar{k}$. Let $f$ be a nonzero  isogeny   of $\AA$. Its kernel $\ker f$ is a subgroup scheme of $\AA$ (since it is closed in the Zariski topology because of continuity and under $\oplus$ because of homomorphism). It contains $0_\AA$ and so its connected component, which is, by definition, an Abelian variety. 

Since $\AA$ is simple and $f\neq 0$ this component is equal to $\{0_\AA\}$. But it has finite index in $\ker f$ (Noether property) and so $\ker f$ is a finite group scheme.

\subsubsection{Computing isogenies between Abelian varieties}\label{sect-3}

Fix a field $k$ and let $\AA$ be an Abelian variety over $k$.  Let $H$ denote a finite subgroup scheme of $\AA$.  From the computational point of view we have the following problems:

\begin{itemize}
\item Given $\AA$ and $H$, determine   $\B:=\AA/H$ and the isogeny $\AA \to \B$.

\item Given two Abelian varieties $\AA$ and $\B$, determine if they are isogenous and compute a rational expression for an isogeny $\AA \to \B$. 
\end{itemize}

There is a flurry of research activity in the last decade to solve these problems explicitly for low dimensional Abelian varieties; see \cite{L-R}, \cite{L-R-2} among many others.   For a survey and some famous conjectures on isogenies see \cite{Frey}.  

\begin{rem} For elliptic curves (Abelian varieties of dimension $1$) and for Jacobians of curves of genus $2$  we shall come back to these questions in more detail.
\end{rem}

\subsection{Torsion points and Tate modules}

The most classical example of an isogeny is the  scalar multiplication by $n$ map  $[n] : \, \AA \to \AA$.
The kernel of $[n]$ is a group scheme of order $n^{2\dim \AA}$ (see \cite{Mum}). We denote  by $\AA [n]$  the group $\ker  [n] (\bar{k}) $.   The elements in $\AA[n]$ 
are called $n$-\textbf{torsion points} of $\AA$. 

\begin{lem}\label{dual}
Let $f : \AA \to \B$ be a degree $n$ isogeny.  Then there exists an isogeny $\hat f : \B \to \AA$ such that
\[ f \circ \hat f = \hat f \circ f = [n]. \]
\end{lem}

\begin{cor} 
Let $\AA$ be an absolutely  simple Abelian variety.  Then $\End(\AA)^0$ is a skew field. 
\end{cor}

\proof Every endomorphism $\neq 0$ of $\AA$ is an isogeny, hence invertible in $\End(\AA)^0$.
\qed

\begin{thm}\label{thm-1} 
Let $\AA/k$ be an Abelian variety, $p = \char k$, and $\dim \AA= g$. 

\begin{itemize}
\item[i)] If $p \nmid \, n$, then $[n]$ is separable, $\# \AA[n]= n^{2g}$ and $\AA[n]\iso (\Z/n\Z)^{2g}$.

\item[ii)] If $p \mid n$, then $[n]$ is inseparable.  Moreover, there is an integer $0 \leq i \leq g$ such that 
\[ \AA [p^m] \iso (\Z/p^m\Z)^i, \; \text{for all } \; m \geq 1. \]
\end{itemize}
\end{thm}

If $i=g$ then $\AA$ is called \textbf{ordinary}.  If $\AA[p^s](\bar k)= \Z/p^{ts}\Z$ then the abelian variety has \textbf{$p$-rank} $t$. If $\dim \AA=1$ (elliptic curve) then it is called \textbf{supersingular} if it has $p$-rank 0.\footnote{For an alternative definition see \cref{lift}.}
An abelian variety $\AA$ is called \textbf{supersingular} if it is isogenous to a product of supersingular elliptic curves. 

\begin{rem}
If $\dim \AA \leq 2$ and $\AA$ has $p$-rank 0 then $\AA$ is supersingular. 
This is not true for $\dim \AA\geq 3$.
\end{rem}

Let $l$ be a prime 
that is (here and in the following)
 different from $p=\char k$ and $k \in \N$.  Then,
\[ [l] \AA \left[l^{k+1}\right] = \AA [l^k].\]
Hence, the collection of groups 
\[ \dots,  \AA[l^{k+1}] , \dots , \AA[l^k], \dots \]
forms a projective system. The $l$-adic Tate module of $\AA$ is 
\[ T_l (\AA) := \varprojlim   \,  \AA[l^k].\]
\begin{lem}\label{tate-module}
The Tate module $T_l (\AA)$ is a $\Z_l$-module isomorphic to $\Z_l^{2 \dim \AA}$. 
\end{lem}

\subsection{$l$-adic  representations and characteristic polynomials}

\subsubsection{Galois representations}
Torsion points on abelian varieties are used to construct very important representations of the Galois group of $k$.  Let $n$ be relatively prime to $p$ and $g= \dim \AA$. Then  $G_k$  acts on $\AA[n]$  which gives rise to a representation 
\[ \rho_{\AA, n} :  G_k \to \Aut \left( \AA[n] \right) \]
and after a choice of basis in $\AA[n]$   yields a  representation
\[ \rho_{\AA, n} : G_k \to  GL_{2g} (\Z/n\Z) \]
This action extends in  a natural way to   $T_l (\AA)\otimes \Q_\ell$ and therefore to a  $\ell$-adic representation $\tilde\rho_{\AA, l}$ which is called the \textbf{$l$-adic Galois representation attached to $\AA$}. 

\subsubsection{Representations of endomorphisms}
Let $\phi$ be an endomorphism of the $g$-dimensional Abelian variety  $\AA$. By restriction $\phi$ induces a $\Z$-linear map $\phi_n$  on $\AA[n]$. Since the collection ($\phi_{\ell^k}$) is compatible with the system defining $T_\ell(\AA)$ it yields a $\Z_\ell$-linear map $\tilde{\phi}_\ell$ on $T_\ell(\AA)$ .

Applying this construction to all elements in $\End(\AA)$ we get an injection (since $\AA[\l^\infty]:=\cup_{k\in \N}\AA[\ell^k])$ is Zariski-dense in $\AA$) from $\End(\AA)$ into $Gl(2g,\Z_\ell)$. By tensorizing with $\Q_\ell$ we get the $\ell$-adic representation 
\[  \tilde{\eta}_\ell: \End(\AA)\otimes \Q_\ell \rightarrow Gl_{2g}(\Q_\ell).\]

\begin{thm}
$\tilde{\eta}_\ell$ is injective.
\end{thm}

For a proof see   \cite[Theorem 3, p.176]{Mum}.  This result has important consequences for the structure of $\End^0(\AA)$, more precisely     $\End^0(\AA)$ is a  $\Q$-algebra of dimension $\leq 4\dim (\AA)^2$. 

Adding more information (see Corollary 2 in \cite{Mum}) one gets that $\End^0 (\AA)$ is a semi-simple algebra, and by duality (key word Rosati-involution) one can apply  a complete classification due to Albert of \emph{possible} algebra structures on $\End^0(\AA)$, which can be found on \cite[pg. 202]{Mum}. 

The question is: Which algebras occur as endomorphism algebras? The situation is well understood if $k$ has characteristic $0$ (due to Albert)  but wide open in characteristic $p>0$. For $g=1$ (elliptic curves) everything is explicitly known due to M. Deuring. We describe the results in  \cref{lift}.  For curves of genus $2$ we give an overview on results in  \cref{gen-2}

\noindent \emph{Characteristic Polynomial:}  For $\phi\in \End^0(\AA)$ let $\tilde{\phi}_\ell $ be its $\ell$-adic representation.  Denote its characteristic polynomial by  $\chi_{\ell,\phi}(T)\in \Z_\ell[T]$.

\begin{thm}[Weil]\label{weil-thm}
  $\chi_{\ell,\phi}(T)$ is a monic polynomial  $\chi_\phi(T) \in \Z[T]$ which is independent of $\ell$.
We have 
\[ \chi_\phi(\phi)   \equiv 0  \mbox{ on } \AA,\]
and so it is justified to call $\chi_\phi(T)$ the \textbf{characteristic polynomial} of $\phi$.
\end{thm}

The degree of $\chi_\phi(T)$ is $2\dim (\AA)$, the second-highest coefficient is the negative of the trace of $\phi$, and the constant coefficient is equal to the determinant of $\phi$. 

\subsubsection{Frobenius representations}
Let $\AA$ be a $g$-dimensional Abelian variety defined over $\F_q$, where $q=p^d$ for a prime $p$ and $\bar \F_q$ the algebraic closure of $\F_q$.  Let  $\pi \in \Gal (\bar \F_q / \F_q)$ be the Frobenius automorphism of $\Fq$,  given by
\[ \pi \, : \, x \to  x^q.\]
%
%
%

Since $\Gal (\bar \F_q / \F_q)$ is topologically generated by $\pi$  and because of continuity  $\rho_{\AA,n}$ is determined by $\rho_{\AA, n} (\pi)$. We define
\begin{equation}\label{Frob-1}
 \chi_{\AA, q } (T):= \chi (T) \left( \tilde\rho_{\AA, l} (\pi) \right) \in  \Z_\ell[T]
\end{equation}
as the characteristic polynomial of the image of $\pi$ under $\tilde\rho_{\AA, l}$.

\begin{lem}[Weil]  
  $\chi_{\AA, q } (T)$   is a monic polynomial of degree $2g$   in $\Z[T]$,  independent of $\ell$,  and for all $n\in \N$ we get
\[ \chi_{\AA, q } (T) \equiv  \chi \left( \rho_{\AA, n} (\pi)  \right)  \mod n .\]
\end{lem}

\begin{lem}[Tate] 
Let $k=\Fq$. The $\ell$-adic representation $\tilde\rho_{\AA, l}$ is semi-simple and so is determined by  $\chi (T) \left( \tilde\rho_{\AA, l} (\pi) \right)$.\footnote{An analogous result for $k=K$ a number field is the main result of Faltings on his way to prove Mordell's conjecture.}
\end{lem}

\paragraph{\emph{Geometric Interpretation:}}   We continue to assume that $\AA$ is an Abelian variety defined over $\Fq$. Hence $\pi$ acts on the algebraic points of $\AA$ by exponentiation on coordinates with $q$. This action induces an action on the function field $\Fq(\AA)$ given again by exponentiation by $q$.

This action is polynomial, and so it induces a morphism on $\AA$. Without loss of generality we can assume that this morphism  fixes $0_\AA$ and so is an endomorphism $\phi_q$ called the \textbf{Frobenius endomorphism}.

By definition it follows that the characteristic polynomial of the $\ell$-adic representation of $\phi_q$ is equal to the characteristic polynomial $\chi_{\AA,q}(T)$ of the $\ell-adic $ Galois representation of $\pi$.

So for given $\AA$, the Frobenius automorphism plays a double role as Galois element and as endomorphism, and this is  of great importance for the arithmetic of Abelian varieties over finite fields.

The explicit knowledge of $\phi_q$ yields immediately that it is purely inseparable and
\[ \deg \phi_q  =  [k(\AA) : \pi^\star k(\AA)] =  q^g.\]
%
%
\begin{defi}
$\chi_{\AA,q}(T)$ is the characteristic polynomial of the Frobenius endomorphism $\phi_q$ of $\AA$.
\end{defi}

Its importance for the arithmetic of Abelian varieties over finite fields becomes evident by the following theorem.
\begin{thm}[Tate]
Let $\AA$ and $\B$ be Abelian varieties over a finite field $\F_q$ and $\chi_{\AA}$ and $\chi_{\B}$ the characteristic polynomials of their Frobenius endomorphism and $l \neq p$ a prime. The following are equivalent. 

\begin{itemize}
\item[i)] $\AA$ and $\B$ are isogenous.

\item[ii)] $\chi_{\AA, q} (T) \equiv  \chi_{\B, q} (T)$

\item[iii)]  The zeta-functions for $\AA$ and $\B$ are the same.  Moreover, $\# \AA (\F_{q^n}) = \# \B (\F_{q^n)}$ for any positive integer $n$. 

\item[iv)] $T_l (\AA) \otimes \Q \iso T_l (\B) \otimes \Q$
\end{itemize}
\end{thm}

$\chi_{\AA, q} (T)$ is the most important tool  for \textbf{counting points} on $\AA(\Fq)$:
Since $\phi_q$ is purely inseparable the endomorphism $\phi_q -id_\AA$ is separable, and hence $\deg\ker(\phi_q-id_\AA) $ is equal to the number of elements in its kernel. Since $\pi$ fixes exactly the elements of $\Fq$ the endomorphism $\phi_q$ fixes exactly  $\AA(\Fq)$ and so    $\ker(\phi_q-id_\AA)(\bar{\Fq})=\AA(\Fq)$.
By linear algebra it follows that: 
\begin{thm}\label{count} The number of points over $\Fq$ is given by
\[ 
\# (\AA (\Fq)  )=\chi_{\AA,q}(1).
\]
\end{thm}

The importance of this observation for \emph{algorithms} for the computation of $\#(\AA(\Fq)$ is due to one of the deepest results (\textbf{Hasse} for $g=1$ and \textbf{Weil} for general $g$) in the arithmetic of 
Abelian varieties over finite fields, which is the analogue of the Riemann Hypothesis in number theory.

\begin{thm}\label{Hasse-Weil} 
Let $\AA$ be an Abelian variety of dimension $g$ over $\Fq$.  The zeroes $\lambda_1,\dots,\lambda_2g$  of the characteristic polynomial of the Frobenius endomorphism $\chi_{\AA, q} (T)$ have the following properties:

\begin{itemize}
\item Each $\lambda_i$ is an algebraic integer.
\item After a suitable numeration one gets for $1\leq i\leq g$ 
\[ \lambda_i\cdot \lambda_{i+g}=q.\]

\item The complex absolute value $|\lambda_i|$ is equal to $\sqrt{q}$.
\end{itemize}
\end{thm}
For the proof we refer to \cite{Mum}. It is evident that this theorem yields bounds for the size of the coefficients of $\chi_{\AA, q} (T)$ deepening only on $g$ and $q$ and so estimates for the size of $\# (\AA (\Fq))$. We state the
\begin{cor} 
\[ |\# (\AA (\Fq)) - q^g | =  \mathcal{O} (q^{g-1/2}).\]
\end{cor}

\begin{rem}
If $\AA$ is the Jacobian of a curve $\CC$ of genus $g$ one can use this result to prove the Riemann Hypothesis for curves over finite fields:
\[ 
| \CC (\Fq) | - q - 1 |  \leq 2g \sqrt{q}.
\]
\end{rem}

For a proof and refinements see \cite{Sti}.   In the next few sections we will focus on some special cases of Abelian varieties, namely Jacobian varieties and more specifically  Jacobians   of hyperelliptic curves. 

\section{Projective Curves and Jacobian Varieties}\label{jacobians}

\subsection{Curves}
First let us establish some notation and basic facts about algebraic curves.   In this paper the notion \emph{ curve}   is an absolutely irreducible projective variety of dimension $1$ without singularities.

At some rare points of the following discussion it is convenient to have that $\CC(k)\neq \emptyset$, and without loss of generality we then can assume that there is a point $P_\infty$ "at infinity", i.e.  in $\CC(k)\setminus U_0$. If we have to study curves with different properties (like being affine or having singularities) we shall state this explicitly.

Let $\CC$ be a curve defined over $k$.  Hence there is $n\in \N$ and a homogeneous \emph{prime} ideal   $ I_\CC \subset k[X_0,\dots,X_n]$ such that, with $R=k[X_0,\dots,X_n]/I_\C$, we have

\begin{enumerate}
\item  $\CC$ is the scheme consisting of the topological space $\mathrm{Proj}( R)$ and the sheaf of holomorphic functions given on open subsets $U$ of $ \mathrm{Proj}(R)$ by the localization with respect to the functions in $R$ not vanishing on  $U$.  

\item The dimension of $\CC$ is one, i.e. for every non-empty affine open subset $U\subset \mathrm{Proj}( R)$ the  ring  of holomorphic functions $R_U$ on $U$ is a ring with Krull dimension $1$.

\item $\CC$ is regular, i.e. the localization of $R$ with respect to every maximal  ideal $M$ in $R$ is a discrete valuation ring $R_M$ of rank $1$.  The equivalence class of the valuations attached to $R_M$ is  the  \emph{place} $\p$ of $\CC$, in this class the valuation with value group $\Z$ is denoted by $w_M$. Alternatively we use the notation $R_\p$ and $w_\p$.  A place $\p$ of $\CC$ is also called \emph{prime divisor} of $\CC$.

\item  (Absolute irreducibility) $\,\, I_\CC\cdot\bar{k}[X_0,\dots,X_n]$ is a prime ideal in $\bar{k}[X_0,\dots,X_n]$.  This is equivalent with: $k$ is algebraically closed in $\mathrm{Quot}(R)$.
\end{enumerate}
As important consequence we note that for all open $\emptyset\neq U\neq \CC$ the ring $R_U$ is a \emph{Dedekind domain}.

\subsubsection{{Prime Divisors and Points}} \,  The set of all places $\p$ of the curve $\CC$ is denoted by $\Sigma_\CC(k)$.  The \emph{completeness} of projective varieties yields:
\begin{prop}
There is a one-to-one correspondence between $\Sigma_\CC(k)$ and the equivalence classes of valuations of $k(\CC)$,
which are trivial on $k$.
\end{prop}

Let $ \p\in \Sigma_\CC(k)$ be a prime divisor with corresponding maximal ideal $M_\p$ and valuation ring $R_\p$.  We have a homomorphism
\[r_\p: R_\p\ra R_p/M_\p =:L \]
where $L$ is a finite algebraic extension of $k$. 

\begin{defi}
The \textbf{degree} of the prime divisor $\p$ is $\deg(\p):=[L:k]$.
\end{defi}

If $\deg(\p)=1$ then $L=k$ and $r_\p$ induces a morphism from $\Spec(k)$ into $\CC$ and so corresponds  to a point $P\in \CC(k)$, uniquely determined by $\p$. 
More explicitly,  the point $P$ has the homogeneous coordinates $(y_0:y_1:\dots,:y_n)$ with $y_i=r_\p(Y_i)$. 

\begin{lem}
The set $\Sigma^1_\CC(k)$ of prime divisors of $\CC$ of degree $1$ is in bijective correspondence with  the set of $k$-rational points $\CC(k)$ of the curve $\CC$.
\end{lem}

Now look at $\CC_{\bar{k}}$, the curve obtained from $\CC$ by constant field extension to the algebraic closure of $k$.   Obviously, every prime divisor of $\CC_{\bar{k}}$ has degree $1$.

\begin{cor}
The set of prime divisors of $\CC_{\bar{k}}$ corresponds one-to-one to the points in $\CC_{\bar{k}}(\bar{k})$.
\end{cor}

Let's go back to $k$. Since $\bar{k}/k$ is separable,  every equivalence class $\p$ of valuations of $k(\CC)$ trivial on $k$  has $\deg(\p)=d$ extensions to $\bar{k}$ and these extensions are conjugate under the operation of $G_k$ (Hilbert theory of valuations). Denote these extension by  $(\tilde{\p}_1, \dots,\tilde{\p}_d)$ and the corresponding points in $ \CC_{\bar{k}}(\bar{k})$ by $(P_1,\dots,P_d)$. Then $\{P_1,\dots,P_d\}$ is an orbit under the action of $G_k$ and we have: 

\begin{cor}\label{orbit}
$\Sigma_\CC(k)$ corresponds one-to-one to the $G_k$-orbits of $ \CC_{\bar{k}}(\bar{k})$.
\end{cor}
%
\subsection{Divisors and Picard groups}\label{functor}

Given a curve $\CC/k$, the group of $k$-rational divisors $\Div_\CC(k)$   is defined as follows. 

\begin{defi}
$\Div_{\CC}(k)=\bigoplus_{\p\in \Sigma_\CC(k)} \Z\cdot\p$, i.e. $\Div_{\CC}(k)$ is the free abelian group with base $\Sigma_\CC(k)$.
\end{defi}

Hence a \textbf{divisor} $D$ of  $\CC$  is a  formal sum 
\[ D=\sum_{\p\in \Sigma_{\CC}(k)} z_\p \, P\]
where   $z_\p \in \Z$ and $z_\p = 0$ for all but finitely many prime divisors $\p$.  The degree of a divisor is defined as 
\[\deg(D) : = \sum_{\p\in \Sigma_{\CC}(k)}z_\p.\]
As we have seen in \cref{orbit} we can interpret divisors as formal 
sums of $G_k$-orbits in $\CC_{\bar{k}}(\bar{k})$.  But we remark that taking points in $\CC(k)$ is in general not enough to get all $k$-rational divisors of $\CC$.
The map $$D\mapsto \deg(D)$$ is a homomorphism from $\Div_{\CC}(k)$ to $\Z$.   Its kernel is  the subgroup $\Div_{\CC}(k)^0$ of divisors of degree $0$.  

\begin{exa}
Let $f\in k(\CC)^*$ be a meromorphic function on $\CC$. For $\p\in \Sigma_\CC(k)$ we have defined the normalized valuation $w_\p$. The \emph{divisor of} $f$ is defined as
\[(f)=\sum_{\Sigma_{\CC}(k)}w_\p\cdot\p.\]
It is not difficult to verify that $(f)$ is a divisor, and that its degree is $0$, see \cite{Sti}. Moreover $(f\cdot g)=(f)+(g)$ for functions $f,g$, and $(f^{-1})=-(f)$. 
The completeness of $\CC$ implies that $(f)=0$ if and only if $f\in k^*$, and so $(f)$ determines $f$ up to scalars $\neq 0$. 
\end{exa}

Thus, the set of principal divisors   $\PDiv_\CC(k)$ consisting of all divisors $(f)$ with  $f\in k(\CC)^*$ is a subgroup of $\Div^0_\CC(k)$.
\begin{defi}The group of divisor classes of $\CC$ is defined by
\[\Pic_\CC(k):=\Div_\CC(k)/\PDiv_\CC(k)\]
and is called the \textbf{divisor class group} of $\CC$.    The group of divisor classes of degree $0$ of $\CC$  is defined by
$$\Pic ^0_\CC(k):=\Div^0_\CC(k)/\PDiv_\CC(k)$$ 
and is called the \textbf{Picard group} (of degree $0$) of $\CC$.
\end{defi}

\noindent \emph{The Picard Functor:}   Let $L$ be a finite algebraic  extension of $k$ and $\CC_L$ the curve obtained from $\CC$ by constant field extension.
Then places of $k(\CC)$ can be extended to places of $L(\CC_L)$.  By the conorm  map  we get an injection of $\Div_\CC(k)$ to $\Div_{\CC_L}(L)$. The well known formulas for the extensions of places yield that 
\[ \conorm_{L/k}(\Div^0_\CC(k))\subset \Div^0_{\CC_L}(L)\]
and that principal divisors are mapped to principal divisors.  Hence we get a homomorphism 
\[\conorm _{L/k}: \Pic^0_\CC(k)\ra \Pic^0_{\CC_L}(L)\]
and therefore a functor   
\[ \Pic^0: L\mapsto \Pic^0_{\CC_L}(L)\]
from the category of algebraic extension fields of $k$ to the category of abelian groups.  Coming  "from above" we have a Galois theoretical description of this functor.  Clearly, 
\[ 
\Div_{\CC_L} (L) = \Div_{\CC_{\bar k} } ( \bar{k} )^{G_L}
\]
and the same is true for functions. With a little bit of more work one sees that an analogue result is true for  $\PDiv_{\CC_L}(L)$ and for $\Pic^0_{\CC_L}(L)$. 

\begin{thm}\label{picard} 
For any curve $\CC_k$ and any  extension $L/k$  with  $k\subset L\subset \bar{k}$  the functor
\[ L\mapsto \Pic^0_{\CC_L}(L) \]
is the same as the functor 
\[ L\mapsto \Pic^0_{\CC_{\bar{k}}}(\bar{k})^{G_L}.\]
In particular, we have 
\[ \Pic^0_{\CC_{\bar{k}}}(\bar{k})=\bigcup_{k\subset L\subset \bar{k}}\Pic^0_{\CC_L}(L), \]
where inclusions are obtained via conorm  maps.
\end{thm}

\begin{rem}
For  a finite extension $L/k$ we  also have the norm map of places of $\CC_L$ to places of $\CC_k$  induces a homomorphism from $\Pic^0_{\CC_L}(L)$ to $\Pic^0_\CC(k)$. In general, this map will be neither injective nor surjective.
\end{rem}

It is one of the most important facts for the theory of curves that the functor $\Pic^0$ can be represented: There is a variety $\J_\CC$ defined over $k$ such that for all extension fields $L$ of $k$ we have a functorial equality
\[\J_\CC(L)=\Pic^0_{\CC_L}(L).\]
$J_\CC$ is the \textbf{Jacobian variety} of $\CC$. This variety will be discussed soon.

\subsection{Riemann-Roch Theorem}
Here we take as guideline the book \cite{Sti} of H. Stichtenoth.

\subsubsection{Riemann-Roch Spaces}
We define a partial ordering of elements in  $\Div_\CC (k)$ as follows;  
\[ D=\sum_{\p\in \Sigma_\CC(k)}z_\p P \]
is \emph{effective} ($D \geq 0$)  if $ z_\p \geq 0$ for every  $\p$,   and $D_1\geq D_2$ if $D_1-D_2\geq 0$.

\begin{defi}
The \textbf{Riemann-Roch space} associated to $D$ is
\[\L(D)=\{f\in k(\CC)^*\mbox{ with } (f)\geq -D\}\cup\{0\}.\]
\end{defi}
So the elements $x\in \L(D)$ are defined by the property that $w_\p(x)\geq -z_\p$ for all $\p\in \Sigma_\CC(k)$.  Basic properties of valuations imply immediately that $\L(D)$ is a vector space over $k$. This vector space has positive dimension  if and only if there is a function $f\in k(\CC)^*$ with $D+(f)\geq 0$, or equivalently,  $D\sim D_1$ with $D_1\geq 0$.

Here are some immediately obtained facts:
$\L(0)=k$ and if $\deg(D)<0$ then $\L(D)=\{0\}$.  If $\deg(D)=0$ then either $D$ is a principal divisor or $\L(D)=\{0\}$.
The following result is easy to prove but fundamental.
\begin{prop}
Let $D=D_1-D_2$ with $D_i\geq 0$. Then 
\[\dim(\L(D))\leq \deg (D_1)+1.\]
\end{prop}
We remark that for $D\sim D'$ we have $\L(D)\sim \L(D')$.   In particular $\L(D)$ is a finite-dimensional $k$-vector space.

\begin{defi}
$\ell(D):=\dim_k(\L(D))$.
\end{defi}

To compute $\ell(D)$ is a fundamental problem in the theory of curves. It is solved by the Theorem of Riemann-Roch.    For all divisors $D$ we have the inequality
\[ \ell(D)\leq \deg(D)+1. \]
For a proof one can assume that $\ell(D) > 0$ and so $D\sim D'>0$.
The important fact is that one can estimate the interval given by the inequality. 

\begin{thm}[\textbf{Riemann}]  \label{riemann}
For given curve $\CC$ there is a minimal number $g_C\in \N\cup \{0\}$ such that for   all $D\in \Div_\CC$ we have 
\[\ell(D)\geq \deg(D)+1-g_\CC.\]
\end{thm}

 For a proof see    \cite[Proposition 1.4.14]{Sti}.  Therefore,  
\[ g_\C=\max \{ \deg{D}-\ell(D)+1;\,\, D\in \Div_\CC(k) \} \]
exists and is a non-negative integer    independent of $D$.

\begin{defi} 
The integer  $g_\CC$ is called the \textbf{genus} of $\CC.$
\end{defi}

We remark that the genus does not change under constant field extensions because we have assumed that $k$ is perfect.  This can be wrong in general if the constant field of $\CC$ has inseparable algebraic extensions. 

\begin{cor} 
There is a number $n_\CC$ such that for $\deg(D) > n_\CC$ we get equality
$$\ell(D)=\deg(D)+1-g_\CC.$$
\end{cor}

\cref{riemann} together with its corollary is the  "Riemann part" of the Theorem of Riemann-Roch for curves. To determine $n_\CC$ and to get more information about the inequality for small degrees one needs canonical divisors. 

\subsubsection{Canonical Divisors}
Let $k(\CC)$ be the function field of a curve $\CC$ defined over $k$. To every $f\in k(\CC)$ we attach a symbol $df$, the \emph{differential}   of $f$ lying in a $k(\CC)$-vector space $\Omega(k(\CC))$ generated by the symbols $df$ modulo the following relations: For $f,g\in k(\CC)$ and $\lambda\in k$ we have:

\begin{itemize}
\item[i)]  $d(\lambda  f+g)=\lambda df +dg$

\item[ii)]  $d(f\cdot g)=f dg+g df$.
\end{itemize}

\noindent The relation between derivations and differentials is given by the

\begin{defi}[Chain rule]
Let $x$ be as above and $f\in k(\CC)$. Then $df = (\partial f/\partial x)  dx$.
\end{defi}
 
As in calculus one shows that the $k(\CC)$-vector space of differentials $\Omega (k(\CC))$ has  dimension $1$ and it is generated by $dx$ for any $x\in k(\CC)$ for which $k(\CC)/k(x)$ is finite and separable. We use a well known fact from the theory of function fields $F$ in one variable.

Let $\p$ be a place of $F$, i.e. an equivalence class of discrete rank one valuations of $F$ trivial on $k$).  Then there exist a function $t_\p \in F$ with $w_\p(t_\p)=1$ and $F/k(t_\p)$ separable. We apply this to $F=k(\CC)$. For all $\p\in \Sigma_\CC(k)$ we choose a function $t_\p$ as above.   For a differential $0\neq \omega\in \Omega (k(\CC))$  we get $\omega=f_\p \cdot dt_\p$.  The divisor $(\omega)$ is given by
\[ (\omega):= \sum_{\p\in \Sigma_\p} w_\p(f_\p)\cdot \p   \]
and  is a called a \textbf{canonical divisor} of $\CC$.

The chain rule implies that this definition is independent of the choices, and the relation to differentials yields that  $(\omega)$ is a divisor.  Since $\Omega(k(\CC) )$ is one-dimensional over $k(\CC)$ it follows that the set of canonical divisors of $\CC$  form a divisor class $k_\CC\in \Pic_\CC(k)$ called the \textbf{canonical class} of $\CC$.   We are now ready to formulate the Riemann-Roch Theorem.

\begin{thm}[\textbf{Riemann-Roch Theorem}]
Let  $W$  be a canonical divisor of $\CC$. For all $D\in \Div_\CC(k)$ we have
\[ \ell(D)=\deg(D)+1-g_\CC+\ell(W-D).\]
\end{thm}
For a proof see    \cite[Section 1.5]{Sti}.  A differential $\omega$ is \emph{holomorphic} if $(\omega)$ is an effective divisor.    The set of holomorphic differentials is a $k$-vector space denoted by $\Omega^0_\CC$  which is  equal to $\L(W)$.    If we take $D=0$ respectively $D=W$ in the theorem of Riemann-Roch we get the following:

\begin{cor} 
$\Omega^0_\CC$ is a $g_\CC$-dimensional $k$-vector space  and $\deg(W)=2g_\CC-2$.
\end{cor}
  
For our applications  there are two further important consequences of the Riemann-Roch theorem.

\begin{cor}  The following are true:
\begin{enumerate}
\item If $\deg(D) > 2 g_\CC-2$ then $\ell(D)=\deg(D)+1-g_\CC.$

\item{In every divisor class of degree $g$ there is a positive divisor.}
\end{enumerate}
\end{cor}

\proof 
Take $D$ with $\deg(D) \geq 2g_\CC -1$. Then $\deg(W-D)\leq -1$ and therefore $\ell(W-D) =0$.  Take $D$ with $\deg(D)=g_\CC$. Then $\ell(D) =1+\ell(W-D)\geq 1$ and so there is a positive divisor in the class of $D$.
\qed

%
\section{Applications of the  Riemann-Roch Theorem}\label{sect-4}

\subsection{The Hurwitz genus formula}
In the theory of curves the notion of a cover is important.

\begin{defi}
Let $\CC,\,\D$ be curves defined over $k$, with $\D$ not necessarily  absolutely irreducible.    A finite surjective morphism $$\eta:\D\ra\CC$$ from  $\D$ to $\CC$ is a \emph{cover morphism}, and if  such a morphism exists we call $\D$ a \textbf{cover} of $\CC$.
\end{defi}

As usual, we denote by 
\[\eta^*: k(\CC)\hookrightarrow k(\D)\]
the induced monomorphism of the function fields and identify $k(\CC)$ with its image. $\eta$ is separable if and only if $k(\D)$ is a separable extension of $k(\CC)$, and $\eta$ is Galois with Galois group $G$ if $k(\D)/k(\CC)$ is  Galois with group $G$.  The cover $\eta$ is geometric if $k$ is algebraically closed in $k(\D)$.

Assume in the following that $\eta$ is separable.  We shall use the well known relations between prime divisors of $k(\CC)$ and those of $k(\D)$ such as extensions, ramifications and sum formulas for the degrees.  In particular we get:

Let $\p\in \Sigma_\CC(k)$ be a prime divisor and   $\PP_1,\dots,\PP_r$  the primes divisors in $\Sigma_\D(k)$ which extend $\p$, written as $\PP/\p$.
Let $t_\p$ be an element in $k(\CC)$ with $w_\p(t_\p)=1.$ The ramification index  $e(\PP_i/\p)=:e_i$ is defined as $e_i={w_{\PP}} _i(t_\p)$, hence there is a function   ${t_\PP} _i$ on $\D$ such that ${t_\PP}_ i^{e_i}=t_\p\cdot u$ with ${w_\PP} _i(u)=0$. The \emph{conorm} of $\p$ is the divisor \[ \conorm  (\p) = \sum_i\PP_i^{e_i}\] and its degree is $[k(\D):k(\CC)]$, the \emph{norm} of $\PP_i$ is $\p$. The cover $\eta$ is  \emph{tamely ramified} if for all $\p\in \Sigma_\CC(k)$ the ramification numbers of all extensions are prime to $\car(k)$.
 
We want to relate the genus of $\D$ to the genus of $\CC$. Let $x\in k(\CC)$ be such that $k(\CC)/k(x)$ is finite separable, and let $dx_\CC$ respectively $dx_\D$ be corresponding  differentials with divisors 
$(dx)_\CC$ and $dx_{D}$. We know that 
\[2g_\C-2=\deg(dx)_\CC\mbox{ and }2g_\D-2=\deg (dx)_\D.\]
We compute the value $z_\p$ respectively ${z_\PP} _i$ of these divisors at $\p\in \Sigma_\CC(k)$ and in the extensions $\PP_1,\dots \PP_r$ with ramification numbers $e_i$.  To ease notation we take $\PP:=\PP_i $, $e_i=e_\PP$ and $t_\PP\in k(\D)$ with $w_\PP=1$.   Then we can choose 
\[
t_\p=u\cdot t_\PP ^{e_\PP}\in k(\CC), 
\]
with $w_\p(u)=0$. By the rules for differentials we get $dt_p= (e_\PP\cdot u\cdot t_\PP^{e_\PP-1}+ u'\cdot t^{e_\PP}_\PP)dt_\PP$ and so
\[
w_\PP(dx)=e_\PP\cdot w_\p(dx) +e_\PP-1.
\]
Summing up over $\PP_1,\dots,\PP_r$ we get that 
\[
\deg \left( \sum_{\PP|\p}z_\PP \right) =\deg \left( \sum_{i=1}^r z_\p\PP_i^e \right)+\sum_{i=1}^r (e_i-1).
\]
Summing up over all $\p\in \Sigma_\CC(k)$ we get the Hurwitz theorem. 

\begin{thm}[Hurwitz]
Any separable, tamely ramified   degree $n$ cover  $\eta: \D\ra \CC$
with  $e_\PP$ the ramification index of $\PP\in \Sigma_\D(k)$ satisfies
\[2g_\D-2=n\cdot(2g_\CC-2)+\sum_{\PP\in \Sigma_{\D}} (e_\PP-1).\] 
\end{thm}

Let us illustrate the theorem with a classical example.

\begin{exa}
Assume that $\CC=\P^1$, the genus  $g_{\P^1}=0$ curve. Let $\D$ be tamely ramified cover of degree $n$ of $\P^1$. Then 
\[ g_\D= 1 - n + \frac 1 2 \,    \sum_{\PP\in \Sigma_\D(k)}(e_\PP-1).\]
In particular $\P^1$ has no unramified extensions.
\end{exa}

The special case $n=2$ will be important for us. Assume that $\car(k)\neq 2$. Then we can apply the Hurwitz formula and get
\[ g_\D=  \frac 1 2 \, r -1,\]
where $r$ is the number of prime divisors of $\P^1$ (or of $\D$ ) which are ramified (i.e. ramification order is larger than $1$) under $\eta$.  

\subsection{Gonality of curves and Hurwitz spaces}
Let $\CC$ be a curve  defined over $k$ and  $\eta:\CC\ra \P^1$  a degree $n$ cover.  We assume that $\CC$ has a $k$-rational point $P_\infty$ and hence a prime divisor $\p_\infty$ of degree $1$.

\begin{defi}
The gonality $\gamma_\CC$ of $\CC$  is defined by 
\[ 
\gamma_\CC=\min \left\{ \deg(\eta):\CC\ra \P^1  \right\}  =  \min  \left\{  [k(\CC):k(x)]    \; | \;     x\in k(\CC)  \right\}   .
\]
\end{defi}
For $x\in k(\CC)^*$,  define the pole divisor $(x)_\infty$ by 
\[ (x)_\infty=\sum_{\p\in \Sigma_\CC(k)}   \max (0,-w_\p(x)) \cdot \p.\]
By the property of conorms of divisors we get  $\deg((x)_\infty) =[k(\CC):k(x)]$ if $x\notin k$ and so
\[ \gamma_\CC = \min \left\{  \deg \, (x)_\infty , \; | \;   x\in  k(\CC) \setminus k.   \right\}
\]
\begin{prop}\label{gon}
For $\g_\CC\geq 2$ we have $\gamma_\CC \leq g$.
\end{prop}

\proof By   Riemann-Roch theorem 
\[ \ell(g_\CC\cdot P_\infty) =1+\ell(W-g_\CC\cdot P_\infty)\]
 and since $g_\CC\geq 2$  the divisor $(W-g_\CC\cdot P_\infty)$ has degree $\geq 0$ and so $\ell(W-g_\CC\cdot P_\infty)\geq 1$. But then $\ell(g_\CC\cdot P_\infty)\geq 2$ and there is a non-constant function $x$ whose pole divisor is a multiple of $\p_\infty$ of order $\leq g_\CC$.

\qed

This proves more than the proposition.

\begin{corollary} For curves $\CC$ of genus $\geq 2$ with prime divisor $\p_\infty$ of degree $1$ there exists a cover
\[ \eta:\CC\ra \P^1\]
with $\deg(\eta)=n \leq g_\CC$ such that $\p_\infty$  is ramified of order $n$ and so the point $P_\infty\in \CC(k)$ attached to $\p_\infty$ is the only point on $\CC$ lying over the infinite point $(0:1)$ of $\P^1$.
\end{corollary}

In general, the inequality in the proposition is not sharp but of size $g/2$ as we shall see below. Curves with smaller gonality are special and so per se interesting.

\subsubsection{Gonality of the generic curve}
Let us assume that $k$ is algebraically closed. We are interested in the classification of isomorphism classes of projective irreducible regular curves of genus  $g\geq 2$.

The moduli scheme  $\M_g$ is a scheme defined over $k$ with the property that it parametrizes these classes (i.e., to every point $P$ there is a unique class of a curve $\CC$ of genus $g$). The coordinates of $\CC$ (chosen in an appropriate affine open neighborhood) are the invariants of $\CC$. It is a classical task to determine such systems of invariants and then to find the curve $\CC$ with these invariants. We shall come back to this in the case of curves of small genus.

\begin{rem}
The scheme $\M_g$ is defined over non-algebraically closed fields $k$ but then it is only a coarse moduli scheme.
\end{rem}

The construction of $\M_g$ is done over $\C$ either by Teichm\"uller theory or, more classically, by Hurwitz spaces (see below), and so over  algebraically closed fields of characteristic $0$ by the so-called Lefschetz principle.  Its existence in the abstract setting of algebraic geometry uses deep methods of geometric invariant theory as developed and studied by Deligne and Mumford in \cite{geominv}.

One knows that $\M_g$ is irreducible and so there exists a generic curve of genus $g$. Moreover the dimension of $\M_g$ is equal to $3g-3$. Curves with special properties (i.e. non-trivial automorphisms or small gonality) define interesting subschemes of $\M_g$. Here is one example.

\begin{defi}
A curve $\CC$ with genus $\geq 2$ is hyperelliptic if and only if it has gonality $2$.
\end{defi}

The subspace of hyperelliptic curves in $\M_g$ is the \emph{hyperelliptic locus} $\M_{g,h}$. We shall see below that this locus has dimension $2g-1$. \newline

\renewcommand\H{\mathcal H}

\noindent \textbf{Hurwitz spaces:}
We continue to assume that $k$ is algebraically closed and consider separable covers $\eta: \CC \to \P^1$ of degree $n$. Then $\eta^*$ allows to identify  $k(\P^1)=:k(x)$ with a subfield of $k(\CC)$.  First, we introduce the equivalence: $\eta \sim \eta'$ if there are isomorphisms $\alpha:\CC\ra \CC'$ and $\beta\in \Aut(\P^1)$ with
\[\beta\circ \eta=\eta'\circ \alpha.\]
The \emph{monodromy group} of $\eta$ is the Galois group of the Galois closure $L$ of $k(\CC)/k(x)$. We embed $G$ into $S_n$, the symmetric group with $n$ letters,  and fix the ramification type of the covers $\eta$. We assume that exactly $r\geq 3$ points in $\P^1(k)$ are ramified (i.e. the corresponding prime divisors have at least one extension to $k(\CC)$  with ramification order $>1$
)
 and that all ramification orders are prime to $\car(k)$. It follows that the ramification groups are cyclic. 

By the classical theory of covers of Riemann surfaces, which can be transferred to the algebraic setting by the results of Grothendieck (here one needs tameness of ramification) it follows that there is a tuple $(\sigma_1, \dots , \sigma_r)$ in $S_n$ such that $\sigma_1 \cdot \cdot \cdot \sigma_r =1$, $\ord(\sigma_i)=e_i$ is the ramification order of the $i$-th ramification point $P_i$ in $L$  and $G:=\< \sigma_1, \dots , \sigma_r\> $ is a transitive group in $S_n$.  We call such a tuple the \textbf{signature} $\sigma$ of the covering $\eta$ and remark that such tuples are determined up to conjugation in $S_n$, and that the genus of $\CC$ is determined by the signature because of the Hurwitz genus formula.

Let $\H_\sigma$ be the set of pairs $([\eta], (p_1, \dots , p_r)$, where $[\eta]$ is an equivalence class of covers of type $\sigma$, and $p_1, \dots , p_r$ is an ordering of the branch points of $\phi$ modulo automorphisms of $\P^1$.
This set carries the structure of a algebraic scheme, in fact it is a quasi-projective variety, the \emph{Hurwitz space}  $\H_\sigma$.   We have the forgetful morphism 
\[\Phi_\sigma: \H_\sigma \to \M_g\]
 mapping $([\eta], (p_1, \dots , p_r)$ to the isomorphic class $[\CC]$ in the moduli space $\M_g$. Each component of $\H_\sigma$ has the same image in $\M_g$.

Define the \textbf{moduli dimension of $\sigma$} (denoted by $\dim (\sigma)$) as the dimension of $\Phi_\sigma(\H_\sigma)$; i.e., the dimension of the locus of genus g curves admitting a cover to $\P^1$ of type $\sigma$. We say $\sigma$ has   \textbf{full moduli dimension} if   $\dim(\sigma)=\dim \M_g$; see \cite{kyoto} for details. 
 
\begin{exa}\label{exa-3}
Take $n=2$, so $G=S_2$, $r\geq 6$ and $\car(k)\neq 2$ and the notations from above. A signature $\sigma$ is completely determined by the $r$ ramification points $P_1,\cdots, P_r$. Hence $\H_\sigma$ consists of classes of hyperelliptic curves of genus $g_r=r/2 \, - 1$ (so $r$ is even). Since we can apply automorphisms of $\P^1$ we can assume that $P_1=(1:0),P_2=(1,1),P_3=(0:1)$ and so we have  $(r-3)$ free parameters modulo a finite permutation group.

So the moduli dimension is $r+3= 2g+2$, and the hyperelliptic locus $\M_{g,h}$ has dimension $2g-1$ and codimension $g-2$. Hence all curves of genus $2$ are hyperelliptic, and for $g\geq 3$ the locus of the hyperelliptic curves has positive codimension.
\end{exa}

For a fixed $g\geq 3$,  we want to find $\sigma$ of full moduli dimension and of minimal degree. This would give a generic covering of minimal degree for a generic curves of genus $g$ and so its gonality. 

A first condition is that $r=3g$. Because of the Hurwitz genus formula this yields conditions for  the ramification cycles, which have to have minimal order.  This is worked out in \cite{generic-3}.

\begin{lem} For any $g\geq 3$ there is a minimal degree $d=\floor{\frac {g+3} 2}$ generic cover
\[\psi_g: \CC_g \to \P^1\]
of full moduli dimension from a genus $g$ curve $\CC_g$ such that it has $r=3g$ branch points and signature:

i) If $g$ is odd, then $ \sigma=(\sigma_1, \dots , \sigma_r)$ such that $\sigma_1, \dots , \sigma_{r-1}\in S_d$ are transpositions and $\sigma_r \in S_d$ is a 3-cycle.

ii) If $g$ is even, then $\sigma=(\sigma_1, \dots , \sigma_r)$ such that $\sigma_1, \dots , \sigma_r \in S_d $ are transpositions.
\end{lem}

\subsubsection{Equations for Curves}
There is a one-to-one correspondence between function fields $F$ of transcendence degree $1$  over the field of constants $k$  with $k$ algebraically closed in $F$  and isomorphism classes of projective regular absolutely irreducible curves $\CC$ with $k(\CC)=F$.  The natural question is: Given $F$, how can one find $\CC$ as embedded projective curve in an appropriate $\P^n$?

The main tool to solve this question are Riemann-Roch systems.  Let $D$ with $\ell(D)=d+1>0$ and $(f_0, f_1, \dots, f_d)$ a base of $\L(D)$. Then
\[
\begin{split}
\Phi_D : \CC(\bar{k}) & \ra \P^d ({\bar k}) \\
P & \mapsto (f_0(P):f_1(p):\dots:f_d(P))\\
\end{split}
\]
is a rational map defined in all points for which $f_0,\dots, f_d$ do not vanish simultaneously.  $\L(D)$ is without base points if this set is empty, and then $\Phi_D$ is a morphism from $\CC$ in $\P^d$.

\begin{lem}
For $g\geq 3$ and $D=W_\CC$ the space $\L(W)=\Omega^0_\CC$ is without base points, and so $\Phi_W$ is a morphism from $\CC$ to $\P^{g_\CC-1}$.
\end{lem}

$\Phi_W$ may not be an embedding but the only exception is that $\CC$ is hyperelliptic, and than the image of $\Phi_W$ is the projective line.

\begin{thm}
Let $\CC$ be a curve of genus $g_\CC>2$ and assume that $\CC$ is not hyperelliptic. Then $\Phi_W$ is an embedding of $\CC$ into $\P^{g_\CC-1}$ and the image is a projective regular  curve of degree $2g_\CC-2$ (i.e. the intersection with a generic hyperplane has $2g_\CC-2$ points).
\end{thm}

So having determined a base of the canonical class of $\CC$ one gets a parameter representation of $\CC$ and then one can determine the prime ideal in $k[Y_0, \dots , Y_{g_\CC}]$ vanishing on $\Phi_W(\CC)$.   $\Phi_W$ is the \textbf{canonical embedding} of $\CC$.

\begin{exa}
Take $g_\CC=3$ and assume that $\CC$ is not hyperelliptic. Then the canonical embedding maps $\CC$ to a regular projective plane curve of degree $4$. In other words: All non-hyperelliptic curves of genus $3$ are isomorphic to non-singular quartics in $\P^2$.
\end{exa}

\noindent \textbf{Plane Curves:} \,  Only very special values of the genus of $\CC$ allow to find plane regular projective curves isomorphic to $\CC$. We have just seen that $g=3$ is such a value. The reason behind is the Pl\"ucker formula, which relates degree, genus and singularities of plane curves. But of course, there are many projective plane curves which are birationally equivalent to $\CC$. 

Take $x\in k(\CC)\setminus k$ with $k(\CC)/k(x)$ separable. Then there is an element $y\in k(\CC)$ with $k(x,y)=k(\CC)$, and by clearing denominators we find a polynomial $G(x,y)\in k[X,Y]$ with $G(x,y)=0$. Then the curve $\CC'$ given by the homogenized polynomial $$G_h(X,Y,Z)=0$$ is a plane projective curve birationally equivalent to $\CC$ but, in general, with singularities. Using the gonality results we can chose $G(X,Y)$ such that the degree in $Y$ is $\floor{\frac {g+3} 2}$.   Using the canonical embedding for non hyperelliptic curves and general projections we can chose $G_h(X,Y,Z)$ as homogeneous polynomial of degree $2g_\CC-2$.

In the next subsection we shall describe a systematic way to find plane equations for hyperelliptic curves.

\subsubsection{Plane equations for elliptic and hyperelliptic curves, Weierstrass normal forms}\label{section-4-2-3}

We first focus on elliptic curves. \newline

\noindent \textbf{Elliptic Curves:}  \, We assume that $\E$ is a curve of genus $1$ with a $k$-rational point $P_\infty$ and corresponding prime divisor $\p_\infty.$ By definition, $\E$ is an \emph{elliptic curve defined over $k$}. We look at the Riemann-Roch spaces $\L_i:= \L(i \cdot \p_\infty)$ and denote their dimension by $\ell_i$. Since $2g_\E-2= 0 $ we can use the theorem of Riemann-Roch to get that  $\ell_i=i$.  Hence $\L_1=\<1\>$, $\L_2=\<1,x\>$ with a function $x\in k(\E)$ with $(x)_\infty=2\p_\infty$, $\L_3= \< 1, x, y \>$ with $(y)_\infty=3\p_\infty$ and $\L_5= \<1, x, x^2, y, x y \>$ with $5$ linearly independent functions.

Now look at $\L_6$. This is a vector space of dimension $6$ over $k$. It contains the seven elements $\{1, x, x^2, x^3, y, xy, y^2\}$ and hence there
is a non-trivial linear relation 
\[
\sum_{0\leq i\leq 3;\,0\leq j\leq 2}a_{i,j} x^i y^2.
\]
Because of the linear independence of $(1, x, x^2, y, xy)$ we get that either $a_{3,0}$ or $a_{0,2}$ are not equal $0$, and since $x^3$ and $y^2$ have a pole of order $6$ in $\p_\infty$ it follows that $a_{0,2}\cdot a_{3,0}\neq 0$. By normalizing we get  $x$ and $y$ satisfy the equation
\[ Y^2 + a_1 X \cdot Y +a_3 Y  =  a_0  X^3 + a_2 X^2 + a_4 X + a_6.\]
By multiplying with $a_0^2$ and substituting $(X,Y)$ by $(a_0 X, a_0 Y)$ we get an \textbf{affine Weierstrass equation} for $\E$:
\[ {W_\E}_{aff}: \;   Y^2 + a_1X \cdot Y + a_3Y = X^3 + a_2 X^2 + a_4X +a_6.\]
The homogenization give the cubic equation 
\[W_\E  : \;  Y^2 \cdot Z+ a_1X\cdot Y\cdot Z +a_3Y\cdot Z^2=a_0X^3+a_2X^2\cdot Z+a_4X\cdot Z^2+a_6\cdot Z^3,\]
which defines a plane projective curve. 

The infinite points of this curve  have $Z=0$, and so the   only infinite  point is $P_\infty =(0,1,0)$ corresponding to the chosen $\p_\infty$. Looking at the partial derivatives one verifies that $\E$ has no singularities if and only if the  discriminant  of the affine equation  ${W_\E}_{aff}$ as polynomial in  $X$ is different from $0$, and that this is equivalent with the condition that $k(\E)$ is not a rational function field.

\begin{thm}
Elliptic curves defined over $k$ correspond one-to-one the isomorphism classes of plane projective curves without singularities given by Weierstrass equations 
\[ W_\E: Y^2\cdot Z+ a_1X\cdot Y\cdot Z +a_3Y\cdot Z^2=X^3+a_2X^2\cdot Z+a_4X\cdot Z^2+a_6\cdot Z^3\]
with non-vanishing   $X$-discriminant.
\end{thm}

Since we are dealing with  isomorphism classes of such curves we can further normalize the equations and finally find invariants for the  isomorphism class of a given $\E$. This is a bit tedious if $\car (k) \, | \, 6$. In this case we refer to \cite{Silverman}. 

Assume that $\car(k) \neq 2, 3$. Then we can use Tschirnhausen transformations to get an equation
\[ W_\E: Y^2\cdot Z =X^3 -g_2X\cdot Z^2-g_3\cdot Z^3\]
and the reader should compare this equation with the differential equation satisfied by the Weierstrass $\wp$-function. 

We use this analogy and define $\Delta(\E)=4g_2^3-27 g_3^2$ and this is, because of the regularity of $\E$, an element $\neq 0$, and
\[ j_E=12^3\frac{4g_2^3}{\Delta_\E}.\]
If $k$ is algebraically closed then $j_\E$ determines the isomorphic class of $\E$.  For an arbitrary $k$ the curve  $E$ is determined up to a \emph{twist}, which is  quadratic if $\car(k)$ is prime to $6$ (see \cite{Silverman}).  \newline

\noindent \textbf{Weierstrass equations for hyperelliptic curves:}  
Let $\CC$ be a curve over $k$ of genus $g\geq 2$  with a degree 2  cover 
\[ \eta : \; \CC \ra \P^1.\]
We assume that there is a point  $P_\infty\in \CC(k)$ corresponding to a prime divisor $\p_\infty$ of $\CC$ of degree $1$.  Take $Q_\infty=\eta(P_\infty)\in \P^1(k)$ and $x\in k(\P^1)$ with $(x)_\infty=\p_{0,\infty}$ with $\p_{0,\infty}$ a prime divisor of degree $1$ of $\P^1$.  Thus,  $\conorm  ( \p_{0,\infty}) = 2\cdot \p_\infty$ and so $\eta$ is ramified in $Q_0$, or  $\conorm  ( \p_{0,\infty}) =  \p_\infty\cdot \p'_\infty$. In any case $\conorm ( \p_{0,\infty})=:D$ is a positive divisor of degree $2$. We define the  Riemann-Roch spaces  $\L_i=\L(i\cdot D)$ and $\ell_i=\dim_k(\L_i)$. 

By assumption $\L_1$ has as base $(1,x)$ and so $\ell_1 =2$.
Since $\deg(g+1)\cdot D >2g-2$ the theorem of Riemann- Roch implies that $\ell_{g+1}=2(g+1)-g+1=g+3$. Hence there is a function $y\in \L_{g+1}$ linearly independent from  powers of $x$.  So $y\notin k[x]$.  The space $\L_{2(g+1)}$ has dimension $3g+3$ and contains the $3g+4$ functions 
\[ 
\{1,x,x^{g+1}, y,x^{g+2},xy,\dots,x^{2(g+1)},x^{g+1}y,y^2\}.
\]
So there is a nontrivial $k$-linear relation between these functions, in which $y^2$ has to have a non-trivial coefficient. We can normalize and get the equation 
\[ y^2+h(x)y=f(x) \; \; \textit{ with }  \; \;  h(x),   f(x) \in k[x] \]
and $\deg  h(x) \leq g+1, \,\deg(f)\leq 2g+2$.   So 
\[ {W_\CC}_{aff}: Y^2+h(X)Y=f(X)\]
is the equation for an affine part $\CC_{aff}$ of a curve birationally equivalent to $\CC$. It is called an \emph{affine  Weierstrass equation}  for $\CC$, and its homogenization is the equation of a projective plane curve  $\CC'$ birationally  equivalent to $\CC$.

The prime divisors of $\CC$ are extensions of prime divisors of $k(x)$ and hence correspond (over $\bar{k}$) to points $(x,y)$ in $\mathbb{A}^2$ or the points lying over $\p_{0,\infty}$. To get more information we use the Hurwitz genus formula and assume for simplicity that $\car(k)\neq 2$ and so $\eta$ is tamely ramified and separable; for the general case see \cite[Section~14.5.1]{book}.

Then we can apply the Tschirnhausen transformation and can assume that $h(x)=0.$  We know that $\eta$ has to have $2g+2$ ramification points. Ramification points of $\eta$ are fixed points under the hyperelliptic involution $\omega$ which generates $\Gal (k(\CC)/k(x)$. Since  $\omega$ acts on points $(x,y)$ by sending it to $(x,-y)$ 
%
%
 the affine ramification points correspond to the zeros of $f(X)$. If $\p_{0,\infty}$ is unramified  then it follows that $f(X)$ has to have $2g+2$ zeros, and so $\deg(f(X))=2g+2$ and all zeros are simple.

Assume that $\p_{0,\infty}$ is ramified. Then there have to be $2g-1$ places with norm $\neq \p_{0,\infty}$ and so $\deg(f(X) )=2g+1$ and again all zeros are different.
Hence in both cases we have that $\CC_{aff}$ is without singularities. This is not true for the point $(0,1,0)$, the only point at  infinity of $\CC'$. It is a singular point, and it corresponds to two points (over $\bar{k}$) on $\CC$ if $\p_{0,\infty}$ is unramified, and to one point on $\CC(k)$ if $\p_{0,\infty}$ is ramified.
For computational purposes the latter case is more accessible. The arithmetic in $k(\CC)$ is analogue to the arithmetic in imaginary quadratic fields. \\

\renewcommand\G{\bar G}

\noindent \textbf{The automorphism group of $\CC$:}\,    We assume that $k$ is algebraically closed.   We identify the places of $k(x)$ with the points of $\P^1= k \cup \{\infty\}$ by their $X$-coordinate. As seen $k(\CC)$ is a quadratic extension field of $k(x)$ ramified exactly at $2g+2$ places $\a_1, \dots , \a_{2g+2}$ of $k(x)$. The corresponding places of $k(\CC)$ are called the {\it Weierstrass points} $Q_1,\dots,Q_{2g+2}$ of $k(\CC)$, the set formed by these points is denoted by  $\mathcal{P}$. 

Weierstrass points play a very important role in the arithmetic of curves. For a detailed discussion see \cite{Sti}. In particular, Weierstrass points of $\CC$ are uniquely determined up to permutations.  So, every automorphism  of $\CC$ and equivalently, of $k(\CC)/k$, fixes $\mathcal{P}$ and  so fixes $k(x)$, and therefore $k(x)$ is the unique subfield of index $2$ in $k(\CC)$.

It follows that $\< \omega \>$ is central in $\Aut (k(\CC)/k)$, and  $\G:=G/\< \tau \>$ is naturally isomorphic to the subgroup of $\Aut(k(x)/k)$ induced by $G$. We have a natural isomorphism $\Gamma:=PGL_2(k) {\overset \iso \to } \Aut(k(x)/k)$. The action of $\Gamma$ on the places of $k(x)$ corresponds under the above identification to the usual action on $\P^1$ by fractional linear transformations $t \mapsto \frac {at+b} {ct+d}$. Since $G$ permutes the Weierstrass points and $2g+2\geq 6$ its action determines  $G$ and so we get an  embedding $\G \to S_n$. 

Since $k(\CC)$ is the unique degree 2 extension of $k(x)$ ramified exactly at $a_1, \dots $, $ a_{2g+2}$, each automorphism of $k(x)$ permuting these $2g+2$ places extends to an automorphism of $k(\CC)$.  Hence under the isomorphism $\Gamma \mapsto \Aut(k(x)/k)$, $\G$ corresponds to the stabilizer $\Gamma_{\mathcal P}$ in $\Gamma$ of the $2g+2$-set $\mathcal P$.   By a theorem of Klein, $\G$ is isomorphic to a cyclic group, a dihedral group, or  $A_4$, $S_4$ or $A_5$. Hence, we can determine $\Aut (\CC)$ as a degree 2 central extension of $\G$ for any fixed genus $g\geq 2$.   \newline


\noindent \textbf{Minimal Degrees:} \; 
We have seen above that non-hyperelliptic curves of genus $\geq 3$ are birational equivalent to plane projective curves of degree $\leq 2g+2.$ 
But in general, this is not the minimal degree one can achieve.   On the other side one has an estimate from below  for the degree of plane curves birational equivalent to a hyperelliptic curve of genus $g\geq 3$; see   \cite{CM} for details.

\begin{prop}
Let $\CC$ be a hyperelliptic curve of genus $g$ and let $\CC'$ be a plane projective curve birationally equivalent to $\CC$.  Then the degree of the equation of $\CC'$ is $\geq g+2$.
\end{prop}

\subsubsection{Addition in Picard groups over $\Fq$}
We take $k=\Fq$ and $\CC$ a curve of genus $g$ defined over $\F_q$.  By a  result of F. K. Schmidt (proved by using Zeta-functions)   curves over finite fields
have a rational divisor $D_0$ of degree $1$ (Caution: Only for curves of genus $\leq 1$ this implies that they have a rational point.) It is not difficult to show that this divisor can be computed effectively.  We use this to present divisor classes $c$ of degree $0$ of $\CC$.

Let $\D_\CC(\Fq)_{>0}^g$ denote the positive divisors of degree $g$ of $\CC$. A consequence of the theorem of Riemann-Roch is that the  map 
\[
\begin{split}
\varphi :  \, \D_\CC(\Fq)_{>0}^g   &   \ra \Pic^0_\CC(\Fq) \\
D     & \mapsto \varphi(D) =D-g\cdot D_1\\
\end{split}
\]
is surjective. A first consequence is that $\Pic^0_\CC(\Fq)$ is a finite  abelian group since there are only finitely many positive divisors of degree $D$ rational over $\Fq$.     Our  aim is to find an algorithm, which computes the addition in $\Pic^0_\CC(\Fq)$ fast.   The main task is the following \emph{reduction}:

Given $D, D'\in \D_\CC(\Fq)_{>0}^g$ find a divisor $S\in\D_\CC(\Fq)_{>0}^g$ with 
\[D+D'-2D_1\sim S-D_1.\]
Then $S-\D_1$ lies in the divisor class that is the sum of the divisor class of $D-D_1$ with the class of $D'-D_1$.  An analogue reduction is well-known from computational number theory and ideal classes of orders. There one uses Minkowski's theorem instead of the Riemann-Roch theorem.

The idea of  F. He{\ss} in \cite{hess-1}  and worked out with many additional details in \cite{Diemhabil} is to use the fact the holomorphic functions in affine open parts of $\CC$  are Dedekind domains and that divisors with support on these parts  can be identified with ideals of these rings. As first step compute (e.g. from the function field $k(\CC)$) a plane curve $\CC'$ birationally equivalent to $\CC$ of  a degree $d$ of size $\O(g)$ (see our arguments above).

The next step is to go to an affine part of $\CC'$ which is without singularities and for which divisors can be identified with ideals in its coordinate ring (approximation properties of functions in function fields can be used since we are only interested in  divisor classes).  Now  the algorithms known from number theory are applicable. The result is given by the following theorem.

\begin{thm}\label{add}[He{\ss}, Diem]
Let $\CC$ be a curve of genus $g$ over $\F_q$.   The addition in the degree 0 class group of $C$ can then be performed in an expected time which is
polynomially bounded in $g$ and $\log(q)$.
\end{thm}

This result is a highlight in algorithmic arithmetic geometry and it opens the access to the Picard groups as abelian groups for arbitrary curves.
Of course, it will be a challenge to implement it. In our context, namely to construct crypto systems, its importance is the \emph{existence} of the algorithm which make certain attacks thinkable!

In the next sections we shall see how we can find explicit algorithms and even formulas to perform group operations in Picard groups of hyperelliptic curve  very rapidly. 

\subsubsection{The Jacobian Variety of a Curve}
In  \cref{functor} we defined the \emph{Picard functor}  $\Pic^0_\CC$ from the category of extension fields $L/k$ into the category of abelian groups given by 
\[ L\mapsto {\Pic^0_{\CC_L}}  (L).\]
In addition we  stated that $\Pic^0_\CC$ is a Galois functor, i.e. that if $k\subset L\subset \bar{k}$ then ${\Pic^0_{\CC}}_L(L)={\Pic^0_{\CC}}_{\bar{k}}(\bar{k})^{G_L}$.
We also announced that this functor is \emph{representable} in terms of algebraic geometry.

More precisely: Let $\CC$ be a curve of positive genus and assume that there exists a $k$-rational point  $P_0\in \CC(k)$ with attached prime divisor $\p_0$.
There exists an abelian variety $\J_\CC$ defined over $k$ and a uniquely determined  embedding 
\[\phi_{P_0}: \CC\ra \J_\CC \; \text{ with } \;  \phi_{P_0}(P_0)={0_\J}_\CC \]
such that
\begin{enumerate}
\item for  all extension fields  $L$ of $k$ we get $J_\CC(L)=\Pic^0_{\CC_L}(L)$  where this equality is given in a functorial way and
\item if $\AA$ is an Abelian variety and $\eta:\CC\ra \AA$ is a morphism sending $P_0$ to $0_\AA$ then there exists a uniquely determined  homomorphism $\psi:\J_\CC\ra \AA$ with $\psi\circ \phi_{P_0}= \eta$.
\end{enumerate}
$\J_\CC$ is uniquely determined by these conditions and is called the \textbf{Jacobian variety} of $\CC$. The map $\phi_{P_0}$ is given by sending a prime divisor $\p$ of degree $1$ of $\CC_\L$ to the class of $\p-\p_0$ in $\Pic^0_{\CC_L}(L)$. \newline

\noindent \textbf{Properties of Jacobian varieties:} \,   From functoriality and universality of the Jacobian it follows that we can introduce coordinates for divisor classes of degree $0$ such that the group law in $\Pic^0_{\CC_L}(L)$ is given by rational functions defined over $k$ and depending only on $\CC$ (and not on $L$). Moreover, we can interpret the norm and conorm  maps on divisor classes geometrically.

Let $L/k$ be a finite algebraic extension. Then the Jacobian variety $\J_{\CC_L}$ of $\CC_L$ is the scalar extension of $\J_\CC$ with $L$, hence a fiber product with projection $p$ to $\J_\CC$. The norm map is $p_*$, and the conorm  map is $p^*$.
%
%

\begin{prop}
If  $f:\CC\ra \D$ is a surjective morphism of curves sending $P_0$ to $Q_0$,  then there is a uniquely determined surjective homomorphism 
\[ f_*:\J_\CC\ra J_\D \]
such that $f_*\circ \phi_{P_0}=\phi_{Q_0}$.
\end{prop}

\proof Apply the universal property to the morphism $\phi_{Q_0}\circ f$ to get $f_*$. The surjectivity follows from the fact that for $k=\bar{k}$
the sums of divisor classes of the form $\p-\p_0$ with $\p\in \Sigma_\CC(k)$ generate  $\Pic^0_\CC(\bar{k})$.
\qed

A useful observation is
\begin{corollary}
Assume that  $\CC$ is a curve of genus $\geq 2$ such that  $\J_\CC$ is a simple abelian variety, and that $\eta: \CC\ra \D$ is a 
separable	cover of degree $>1$. Then $\D$ is the projective line.
\end{corollary}
For the proof use the Hurwitz genus formula and the universal properties of Jacobians.

What about the \textbf{existence} of Jacobian varieties?   Over the complex numbers the classical theory of curves (key words: Riemann surfaces and the Theorem of Abel-Jacobi) is used to prove the existence of Jacobian varieties  already  in the 19-th century. In fact, this notion is historically earlier than the notion "Abelian variety" introduced by A. Weil as most important tool for his proof of the geometric Riemann hypothesis.  By the Lefschetz principle the existence of Jacobian varieties follows for algebraically closed fields of characteristic $0$. 

For a proof  in the framework of Algebraic Geometry (and so over arbitrary ground fields $k$) see Lang~\cite{lang}.  The important fact is that we "know" a birational affine model of $J_\CC$. 

By the Theorem of Riemann-Roch we have a surjective map from $\Sigma^g_\CC(L)$ to $\Pic^0_\CC(L)$ by sending any positive divisor $D$ of degree $g$ to $D- g\cdot \p_0$. We can interpret such positive divisors geometrically. Take the g-fold cartesian product $\CC^g$ of the curve $\CC$ of genus $g$ and embed it (via Segre's map) into a projective space. On this variety we can permute the factors and so have an action of $S_g$, the symmetric group with $g$ letters.  Define the $g$-fold symmetric product $\CC^{(g)}$ by $\CC^g/ S_g$.  Then we can identify $\CC^{(g)}(L)$ with  $\Sigma^g_\CC(L)$ and so define a birational map from $\CC^{(g)}$ to $\J_\CC$.  Taking an affine part of $\CC$ (e.g. found as a regular part of a plane model of $\CC$)  we get an affine variety which is birational equivalent to $\J_\CC$. 

The Jacobian varieties connect the arithmetic in divisor classes of curves  (which is very accessible to algorithms) with the very rich geometric structure of abelian varieties (e.g. isogenies, endomorphisms and $\ell$-adic representations).
 
%
\subsubsection{Construction of curves by period matrices}
It is convenient to assume in the following that $k$ is algebraically closed.
We look at the following task: Assume that a point $P$ in the moduli scheme $\M_g(k)$ is given by coordinates in a certain coordinate system.
How can we find an equation for a curve $\CC$ corresponding to $P$?

It is useful to look at the case that $k=\C$ and at the parametrization  of isomorphism classes principally polarized abelian varieties by period matrices. 
We reformulate the question and ask whether we can find a curve such that the Jacobian has a given period matrix. Of course, the first problem is that not every principally polarized abelian variety is the Jacobian of a curve, and the decision for this is the well-known \emph{Schottky Problem}
which is unsolved till now.

There are two cases where the situation is better: If the dimension of the Abelian variety is $\leq 3$ then such a curve exists, and if 
we are looking for hyperelliptic curves we can solve the Schottky Problem and determine a Weierstrass equation if we know the period matrix.

This latter result is based on \emph{invariant theory}. Details are worked out in the thesis of H.J. Weber \cite{Web} (explicitly for curves up to genus $5$). Important cases for our applications are curves of genus $1$ (use the $j$-invariant), genus $2$ and genus $3$.
We remark that this method works very well over number fields and by reduction, over finite fields, too.

We  shall give more details in the interesting case that the genus of $\CC$ is equal to $2$.

%
\subsubsection{Example: Curves of genus 2}\label{gen-2}
Let $\CC$ be a genus 2 curve defined over a field $k$. By \cref{gon} we have that its gonality is $\gamma_\CC = 2$. Hence, genus 2 curves are hyperelliptic and we denote the hyperelliptic projection by 
$\pi : \CC \to \P^1$.   By the Hurwitz's formula this covering has $r=6$ branch points which are images of the Weierstrass points of $\CC$.  The moduli space has dimension $r-3=3$; see \cref{exa-3}. 
  
The arithmetic of the moduli space of genus two curves was studied by Igusa in his seminal paper \cite{Ig} expanding on the work of Clebsch, Bolza, and others.  Arithmetic invariants by $J_2, J_4, J_6, J_8, J_{10}$ determine uniquely the isomorphism class of a genus two curve. Two genus two curves $\CC$ and $\CC^\prime$ are isomorphic over $\bar k$  if and only if there exists $\l \in \bar k^\star$ such that $J_{2i} (\CC) = \l^{2i} J_{2i} (\CC^\prime)$, for $i=1,\dots, 5$. If $\char k \neq 2$ then the invariant $J_8$ is not needed. 

From now on we assume  $\char k \neq 2$.   Then   $\CC$ has an affine Weierstrass equation 
\begin{equation}\label{eq-g-2}
y^2=f(x)= a_6 x^6 + \dots + a_1 x + a_0, 
\end{equation}
over $\bar k$, with discriminant $\Delta_f = J_{10} \neq 0$.   The moduli space $\M_2$  of genus 2 curves, via the Torelli morphism, can be identified with the moduli space of the principally polarized abelian surfaces $\mathbb A_2$ which are not products of elliptic curves. Its compactification $\mathbb A_2^\star$ is the weighted projective space $\WP (k)$  via the Igusa invariants $J_2, J_4, J_6,  J_{10}$.   Hence, 
\[ A_2 \iso \WP (k) \setminus \{ J_{10} = 0 \}.\]
A point $\p \in \WP$ for $J_2\neq 0$ can be written as 
\[ 
\left[  1 :    \frac 1 {2^4 3^2} \x_1 : \frac 1 {2^6 3^4} \x_2 + \frac 1 {2^4 3^3} \x_1 : \frac 1 {2 \cdot 3^5} \x_3   \right] 
\]
where $\x_1, \x_2, \x_3$ are given as ratios of Siegel modular forms and are called \textit{absolute invariants} and denoted by $i_1$, $i_2$, $i_3$ by other authors; see \cite{Ig-2}.  Two genus 2 curves are isomorphic over $\bar k$ if and only if they have the same absolute invariants. Notice that the absolute invariants are not defined for $J_2=0$. There are different sets of absolute invariants used by many authors, but all of them are not defined over $J_2=0$.  

\paragraph{Recovering the curve from invariants}
Given a moduli point $\p \in \M_2$, with automorphism group of order 2, we can recover the equation of the corresponding curve over a minimal field of definition following Mestre's approach \cite{mestre}, where the point is given in terms of the absolute invariants. The case of automorphism group of order $> 2$ was done in \cite{cardona} and   \cite{ants}.   In all these papers the case when the absolute invariants are not defined had to be treated differently, introducing a new set of invariants.  In \cite{univ-g-2} it is given an equation of the curve in terms of $J_2, J_4, J_6, J_{10}$ without using any absolute invariants. 

In \cite{w-height}, for any number field $K$,  a \textit{height} on the  moduli space  $\WP (K)$ is introduced. This makes it possible to store the \textit{smallest} tuple of invariants in a unique way. This is used  in \cite{beshaj-guest} to create a database of all genus 2 curves with  small height and defined over $K$ including all the twists of minimal moduli height.

\subsubsection{Example: Elliptic Curves}
Let $\E$ be an elliptic curve over $k$, i.e. a curve of genus $1$ with a $k$-rational point. Its isomorphism class
over $\bar{k}$ is uniquely determined by the $j$-invariant. As seen above, $\E$ is isomorphic to a plane curve $\E'$ given by a Weierstrass equation. 

We choose one $k$-rational point $P_\infty$ with prime divisor $\p_\infty$ and projective coordinates such that  $P_\infty=(0:1:0)$ is the infinite point of the curve $\E'$ with equation 
\[ 
Y^2Z+a_1XYZ+a_3YZ^2=X^3+a_2X^2Z+a_4XZ^2+a_6Z^3
\]
and identify $\E$ with $\E'$.

Let $\J_\E$ be the Jacobian variety of $\E$. We look at  
$\phi_{P_\infty} : \E   \ra \J_\E $
given by 
\[  P   \mapsto [\p-\p_\infty] \]
where $[.]$ means divisor class.   Since $2g_\E-2=0$ the  Riemann-Roch theorem implies that for all extension fields $L$ of $k$ in each $L$-rational divisor class of degree $1$ there is exactly one prime divisor $\p$ of degree $1$ corresponding to a point $P\in \E(L)$, and to each divisor class $c$ of degree $0$ there is exactly one prime divisor $\p$ of degree $1$ with $c=[\p-\p_\infty]$. So $\phi_{P_\infty}$ is injective and surjective and hence an isomorphism of projective varieties.  By transport of structure we endow $\E$ with a group structure:\\

\textit{For extension fields $L$ of $k$ and $P_1,P_2\in \E(L)$ define $P_1\oplus P_2$ as the  point belonging to the prime divisor 
in the class $\p_1+\p_2-2\p_\infty$.}   

It is obvious that this makes $\E(L)$ to an abelian group with neutral element $P_\infty$.  We conclude: Three points $P_1$, $P_2$, $P_3$ sum up to $0$ if $\p_1+\p_2+\p_3-3\p_\infty =(f)$ with $f\in L(\E)$.

Now recall that $\E$ has degree $3$ and so lines intersect with $\E$ in $3$ points (counted with multiplicities) and so $f$ defines a line in $\P^2$. Hence $P_1+P_2+P_3=0$ if and only if the three points are collinear, and then $P_1\oplus P_2=\ominus P_3$.   Using coordinates we get an algebraic recipe for addition:

\textit{For given $P\neq P\infty$ take the line through $P$ and $P_\infty$ to get: $\ominus(P)$ is the third intersection point of the line with $\E$ (if this point is equal to $P$ the line is a tangent and $P=\ominus P$ is an element of order $2$). Given two points $P_1\neq P_2$ compute the line through these two points, take its third intersection point $P_3$ with $\E$ to get $P_1\oplus P_2=\ominus P_3$.
}

\smallskip

By elementary algebra one can perform this recipe by writing down formulas in rational functions in $(X,Y,Z)$   and so we get

\begin{thm} 
After the choice of a base point  $P_\infty$ the elliptic curve $\E$ is an Abelian variety of dimension $1$  with neutral element  $P_\infty$ which is equal to $\J_\E$.
\end{thm}

\noindent \textbf{Division polynomials for elliptic curves:}  \, To simplify we assume that 
 $\char k\neq 2, 3$ and so  we can take the affine Weierstrass equation of $E$ as  
\[E : \quad  Y^2=X^3+aX+b, \]
 for $a, b\in k$. 
Recall that for $m\in \N$ the endomorphism $[m]$ of $\E$ is induced by the scalar multiplication by $m$. We want to give formulas for this endomorphism.
\begin{lem}
For any integer $m$ and point $P(x, y) \neq \O$ in $E$, the point $[m]P$ has coordinates
\[
[m]P = \left(   \frac{\phi_m(x,y)}{\psi_m(x,y)^2}, \frac{\om_m(x,y)}{\psi_m(x,y)^3}     \right)
\]
where  the polynomials $\phi_m, \psi_m, \om_ m$ are given by the  recurrences
\begin{equation}
\begin{split}
    \psi_1 &= 1,\\
    \psi_2 &= 2Y^2,\\
    \psi_3 &= 3X^4 + 6aX^2 + 12bX - a^2,\\
    \psi_4 &= (2X^6 + 10aX^4 + 40bX^3 - 10a^2X^2 - 8abX - 2a^3 - 16b^2)2Y^2,  \\
    \ & \dots \\
    \psi_{2m+1}  &= \psi_{m+2}\psi_m^3 - \psi_{m-1}\psi_{m+1}^3     				\quad \text{for $m \geq 2$,}\\
    \psi_2\psi_{2m} &= (\psi_{m+2}\psi_{m-1}^2 - \psi_{m-2}\psi_{m+1}^2)\psi_m 		\quad    \text{for $m\geq 3$.}
\end{split}
\end{equation}
and    
\begin{equation}  
\begin{split}
 \phi_m &= x\psi_m^2 - \psi_{m+1}\psi_{m-1},\\
\om_m &= \psi_{m-1}^2\psi_{m+2} + \psi_{m-2}\psi_{m+1}^2.
\end{split}
\end{equation}
\end{lem}

The proof follows from classical identities of the Weierstrass function $\wp$ if $k=\C$ and is then transferred to arbitrary perfect fields
(see \cite{Silverman}).
The polynomial $\psi_m$ is called the   \textbf{$m$-th division polynomial} and it vanishes in $E[m]$

\begin{cor} All $m$-torsion points $P(x, y)$ of $E$ have coordinates satisfying $\psi_m (x, y) =0$
\end{cor}

This provides a computational approach on how to determine the $m$-torsion points for any given $m\geq 2$. 

\subsection{Cantor's Algorithm}\label{Cantor}
Inspired by the group law on elliptic curves and its geometric interpretation we give an \emph{explicit}  algorithm for the group operations on Jacobian varieties of hyperelliptic curves. 

Take  a genus $g \geq 2$ hyperelliptic curve $\CC$ with at least one rational Weierstrass point  given by the affine Weierstrass  equation
\begin{equation}\label{hyp}
W_\CC : \; y^2 + h(x) \, y = x^{2g+1} + a_{2g} x^{2g} + \dots + a_1 x + a_0 
\end{equation}
over $k$. We  denote the prime divisor corresponding to  $P_\infty =(0:1:0)$  by $\mathfrak {p}_\infty$.   The affine coordinate ring of $W_\CC$ is  
\[  \O   =  k[X,Y]/(Y^2 + h(X) \, \<Y -( X^{2g+1} + a_{2g} X^{2g} + \dots + a_1 X + a_0)\> )\]
and so prime divisors $\p$ of degree $d$ of  $\CC$ correspond to prime ideals $P\neq 0$ with $[\O/P:k]=d$.    Let $\omega $ be the hyperelliptic involution of $\CC$. It operates on $\O$ and on $\Spec(\O)$ and fixes exactly the prime ideals which   "belong" to Weierstrass points, i.e. split up in such points over $\bar{k}$.
 
Following Mumford \cite{Mum} we introduce polynomial coordinates for points in $J_{\CC}(k)$. The first step is to normalize representations of divisor classes. In each divisor class $c\in \Pic^0(k)$ we find a unique \emph{reduced} divisor
\[ D=n_1\p_1+\cdots +n_r \p_r - d\, \p_\infty \]
with $\sum_{i=1}^r n_i\deg (\p_i) = d\leq g$, $\p_i \neq \omega(\p_j)$ for $i\neq j$ and $\p_i\neq \p_\infty$ (we use Riemann-Roch and the fact that $\omega$ induces $-id_{J_\CC}$).

Using the relation between divisors and ideals in coordinate rings we get that $n_1\p_1+\cdots +n_r \p_r$ corresponds to an ideal $I\subset \O$ of degree $d$ and the property that if the prime ideal $P_i$ is such that  both $P$ and $\omega(P)$ divide $I$ then it belongs to a Weierstrass point.  The  ideal $I$ is a free $\O$-module of rank $2$ and so
\[I=k[X]u(X)+k[x](v(X)-Y).\]

\noindent \textbf{Fact}:   $u(X),v(X)\in k[X]$, $u$ monic of degree $d$, $\deg(v)<d$ and $u$ divides $v^2+h(X)v-f(X)$;  see   \cite[Theorem 4.143]{book}.

\smallskip

Moreover, $c$ is uniquely determined by $I$, $I$ is uniquely determined by $(u,v)$ and so we can \emph{take $(u,v)$ as coordinates for $c$.}

\begin{thm}[Mumford representation] \label{Mum-rep}
Let $\CC$ be a hyperelliptic curve of genus $g\geq 2$ with affine equation 
\[ y^2 + h(x)\, y  \, = \,  f(x), \]
where $h, f \in k[x]$, $\deg f = 2g +1 $, $\deg h \leq g$.

 Every non-trivial group element $c \in \Pic^0_\CC (k)$ can be represented in a unique way by a pair of polynomials $u, v \in k[x]$, such that 

i) $u$ is a monic

ii) $\deg v < \deg u \leq g$

iii) $u \, | \, v^2+ vh -f$
\end{thm}

How to find the polynomials $u,v$?   We can assume without loss of generality that $k=\bar{k}$ and identify prime divisors $\p_i$ with points $P_i=(x_i,y_i)\in k\times k$. Take the reduced divisor $D =n_1\p_1+\cdots + n_r \p_r - d \p_\infty$ now with $r=d\leq g$.  Then 
\[  u(X)= \prod_{i=1}^r(X-x_i)^{n_i}. \]
Since $(X-x_i)$ occurs with multiplicity $n_i$ in $u(X)$ we must have for $v(X)$ that 
\[  \left( \frac d {dx}  \right)^j \left[ v(x)^2 + v(x) \, h(x) - f(x)  \right]_{x=x_i} =0, \]
and one determines $v(X)$ by solving this system of equations.\\

\noindent \textbf{Addition:} Take the divisor classes represented by $[(u_1,v_1)]$  and $[(u_2,v_2)]$ and  in "general position". Then the product is represented by the ideal $I\in\O$ 
given by 
\[ \< u_1u_2, u_1(y-v_2), u_2(y-v_1), (y-v_1)(y-v_2) \>.\]
We have to determine a base, and this is done by Hermite reduction.  The resulting ideal is of the form $\<u'_3(X), v'_3(X)+w'_3(X)Y\>$ but not necessarily reduced. To reduce it one uses recursively the fact that $u \, \mid \, (v^2-hv-f)$.

The formalization of this procedure and the treatment of special cases is called \textbf{Cantor's algorithm}. For readers acquainted with algorithmic number theory it may be enlightening to compare this algorithm with the well known method to  add in class groups of imaginary quadratic number fields, going back to Gauss and based on the theory of definite quadratic forms with fixed discriminant. The very explicit and efficient "generic" algorithm can be found in \cite[Algorithm 14.7]{book}.   For curves of genus $2$ a detailed analysis including all special cases is done in \cite[Section 14.3.2]{book}, including a determination of complexity (see Table 14.2 and Table 14.13). For curves of genus $3$ we refer to Section 14.6 in\cite{book}.  \newline

\noindent \textbf{Addition by Interpolation} \, Another approach to describe  addition in the Jacobians of hyperelliptic curves is to use approximation by rational functions;  see \cite{Leitenberger}. This is analogous to the geometric method used for elliptic curves. 

For simplicity we assume that $k=\bar{k}$.  Let $D_1$ and $D_2$ be reduced divisors on $\Jac_k \CC $ given by 
\begin{equation} 
\begin{split}
D_1  & = \p_1 + \p_2 + \dots + \p_{h_1} - h_1 \p_\infty,  \\
D_2  & = \mathfrak{q}_1 + \mathfrak{q}_2 + \dots + \mathfrak{q}_{h_2} - h_2\p_\infty,  \\
\end{split}
\end{equation}
where  $\p_i$ and $\mathfrak{q}_j$ can occur with multiplicities,  and $0 \leq h_i \leq g$, $i=1, 2$. 
As usual we denote by $P_i$ respectively $Q_j$ the points on $\CC$ corresponding to $\p_i$ and $\mathfrak{q}_j$.

Let $g(X)= \frac {b(X)} {c(X)} $ be the unique  rational function going through the points $P_i$, $Q_j$. In other words we are determining $b(X)$ and $c(X)$ such that $h_1+h_2-2r$ points $P_i$, $Q_j$ lie on the curve
\[ Y \, c(X) - b(X) \, =  \, 0.\] 
This rational function is uniquely determined and has the form
\begin{equation}\label{pol} 
Y= \frac {b(X)} {c(X)}    = \frac     {b_0 X^p + \dots b_{p-1} X + b_p} {c_0 X^q + c_1 X^{q-1} + \dots + c_q} 
\end{equation}
where 
\[
p=\frac {h_1+h_2+g-2r-\epsilon} 2, \; \; q = \frac {h_1+h_2-g-2r-2+\epsilon} 2,
\]
$\epsilon$ is the parity of $h_1+h_2+g$.    By replacing $Y$ from  \cref{pol} in  \cref{hyp} we get a polynomial of degree $\max \{ 2p, \, 2q (2g-1) \}$, which gives $h_3\leq g$ new roots apart from the $X$-coordinates of $P_i, Q_j$.   Denote the corresponding points on $\CC$ by $R_1, \dots , R_{h_3}$ and   $\bar R_1, \dots , \bar R_{h_3}$ are the   corresponding symmetric points with respect to the $y=0$ line. 
Then, we define 
\[ D_1 + D_2 = \bar R_1 + \dots \bar R_{h_3} - h_3 \O.\]
For details we refer the reader to \cite{Leitenberger}.  

\begin{rem}
For $g=1, 2$ we can take $g(X)$ to be a cubic polynomial. 
\end{rem}

\begin{exa}[Curves of genus $2$]\label{gen-2-add}
Let $\CC$ be a genus 2 curve defined over a field $k$ with a rational Weierstrass point.  If $\char k \neq 2, 3$ the $\CC$ is  birationally isomorphic  to an affine plane curve with equation
\begin{equation} 
Y^2 = a_5 X^5 + a_4 X^4+ a_3 X^3 + a_2 X^2 + a_1 X + a_0.
\end{equation}
Let $\p_\infty$ be the prime divisor corresponding to the point at infinity. Reduced divisors in generic position  are  given by 
\[ D = \p_1 + \p_2 - 2 \p_\infty \]
where $P_1(x_1, y_1) $, $P_2(x_2, y_2)$ are points in $\CC(k)$ (since $k$ is algebraically closed) and $x_1\neq x_2$.  For any two divisors $D_1=\p_1 + \p_2 - 2\p_\infty$ and $D_2 = \mathfrak{q}_1 + \mathfrak{q}_2 - 2\p_\infty$ in reduced form, we determine the cubic polynomial 
\begin{equation} Y= g(X) = b_0 X^3 + b_1 X^2 + b_2 X + b_3,\end{equation}
going through the points $P_1 (x_1, y_1)$, $P_2(x_2, y_2)$, $Q_1(x_3, y_3)$, and $Q_2(x_4, y_4)$.   This cubic will intersect the curve $\CC$ at exactly two other points $R_1$ and $R_2$ with coordinates 
\begin{equation}
R_1  = \left( x_5,   g(x_5)    \right) \; \text{ and } \; R_2  =  \left( x_6,  g(x_6)    \right),  
\end{equation}
where $x_5$, $x_6$ are roots of the quadratic equation
\begin{equation}  
x^2 + \left( \sum_{i=1}^4 x_i   \right) x + \frac  {b_3^2-a_5} { b_0^2 \prod_{i=1}^4 x_i} = 0. 
\end{equation}
Let us denote by $\overline R_1 = (x_5, - g(x_5) )$ and $\overline R_2 = (x_6, - g(x_6) )$. 
Then,
\begin{equation} 
[D_1] \oplus [D_2] = [\overline R_1  + \overline R_2 - 2 \p_\infty].
\end{equation}
\end{exa}

After having defined explicitly the addition in $\Jac_ \CC$ it is a natural problem that  given a reduced divisor $D \in \Jac_ \CC$, determine explicitly the formulas for $[n]D$,  at least in generic cases similarly as in the case of elliptic curves. 
Hence, one wants to determine explicitly division polynomials (i.e polynomials that have torsion points of order $n$ as zeroes) or  more generally, ideals which define zero-dimensional schemes containing $J_\CC[n]$. There has been a lot of activity on this area lately; see    \cite{onishi},   \cite{Ka-1}, \cite{Ka-2}.

\subsection{Automorphisms  of curves   and their Jacobians} 

Let $\CC$ be an algebraic curve defined over $k$ and $\Jac_k (\CC)$ its Jacobian.  What are the automorphism groups of $\CC$ and $\Jac_k (\CC)$?

\subsubsection{Automorphisms of curves}
Let $\Aut(\CC)$ be the automorphism of the curve $\CC$. If $\CC=\P^1$ the automorphism group over an infinite field $k$ is infinite, as seen in \cref{section-4-2-3}. The same is true for elliptic curves $\E$, if $E(k)$ is infinite, i.e. $k=\bar{k}$, for then there are infinitely many  translations.  But caution: If we look only at automorphisms with a fixed point (which are automatically homomorphisms with respect to the group structure with the fixed point as origin) then the group is finite, well understood, and for generic elliptic curves $\E$ equal to $\{id_\E,-id_\E\}$.
 
For curves of genus $\geq 2$ the picture changes completely. The automorphism group  $\Aut(\CC)$ is a finite group.  The reason is the existence of $2g+2$  Weierstrass points and the faithful action of $\Aut(\CC)$ on these points, which gives an injection of $\Aut(\CC)$ into $S_{2g+2}$. In fact, one can describe all occurring groups. If $\car(k)$ = $0$ one has the \emph{Hurwitz bound}
\[ \# \Aut(\CC)\leq 84( g-1).\]
One gets a stratification of $\M_g$ by strata of curves with the same automorphism group, and the generic curve of genus $G>2$ has trivial automorphism group.

It is very interesting to study curves with large automorphism group; for instance a curve of genus $3$ with group $PSl(2,7)$ (168 automorphisms) is the famous Klein Quartic, which also occurs as modular curve (see below). For more details and a full account of automorphisms of curves see \cite{kyoto}.   

\subsubsection{Automorphisms of Jacobian varieties}
By functoriality it follows that automorphisms of $\CC$ induce automorphisms of $\Jac_\CC$, or, to be more precise, of $(\Jac_\CC,\iota)$ where $\iota$ is  the principal polarization of $\Jac_\CC$ attached to $\CC$. 

\begin{thm}
Let $\CC$ be an algebraic curve and $\AA := \Jac (\CC)$ with canonical principal polarization $\iota$.  Then,
\[
\Aut \CC \iso 
\left\{
\begin{split}
& \Aut (\AA, \iota), \quad \text{if } \CC \text{is hyperelliptic} \\
& \Aut (\AA, \iota)/ \{\pm 1\}, \quad \text{if } \CC \text{is non-hyperelliptic} \\
\end{split}
\right.
\]
\end{thm}

The above result can be used to find Jacobians of genus 3 hyperelliptic curves; see \cite{Weng} and \cref{gen-3}.

\subsubsection{Endomorphism of Abelian varieties}

The ring of endomorphisms of generic Abelian varieties is "as small as  possible". For instance, if $\car(k)=0$ $\End(\AA)=\Z$ in general. If $k$ is a finite field, the Frobenius endomorphism  will generate a larger ring, but again, this will be all in the generic case. A concrete result is  the following \cite{zarhin-1}:
\begin{thm}[Zarhin]\label{zarhin-thm}
Let $\CC$ be a hyperelliptic curves with affine equation $y^2=f(x)$, $n=\deg f$, and $f \in \Q[x]$. If $\Gal (f) $ is isomorphic to $A_n$ or $S_n$  then $\End_{\overline \Q} \, (\Jac \CC ) \iso \Z$.
\end{thm}

The theorem is actually true over any number field $K$. See \cite{zarhin-2} for detailed results on endomorphisms of Jacobians of hyperelliptic and superelliptic curves. 

From this point of view it will be interesting to find Abelian varieties with larger endomorphism rings. This leads to the theory of real and complex multiplication.  For instance, the endomorphism ring of the Jacobian of the Klein quartic contains an order in a totally real field of degree $3$ over $\Q$.  We shall see in \cref{sect-correspondences} that the Jacobians of modular curves have real multiplication.


%
\subsection{Endomorphism ring of an abelian surface}

For $\char k \neq 2$, a point $\p $ in the moduli space $\M_2$ is determined by the tuple $(J_2, J_4, J_6, J_{10})$, for discriminant $D:=J_{10} \neq 0$. In the case of $\char k =2$ another invariant $J_8$ is needed; see \cite{Ig}. 

\medskip

\noindent \textbf{Humbert surfaces:}   For every $D:=J_{10} > 0$ there is a Humbert hypersurface $H_D$ in $\M_2$ which parametrizes curves $\CC$ whose Jacobians admit an optimal action on $\O_D$; see \cite{HM95}.   Points on $H_{n^2}$ parametrize curves whose Jacobian admits an $(n, n)$-isogeny to a product of two elliptic curves.

\medskip

\noindent \textbf{Shimura curves:}   For every quaternion ring $R$  there are irreducible curves $S_{R, 1}$, $\dots$, $S_{R, s}$ in $\M_2$ that parametrize curves whose Jacobians admit an optimal action of $R$. Those $S_{R, 1}$, $\dots$ , $S_{R, s}$ are called  \emph{Shimura curves}.      

\medskip

\noindent \textbf{Curves with complex multiplication:}   Curves whose Jacobians admit complex multiplication correspond to isolated points in $\M_2$.   We have the following:

\begin{prop}   
$\Jac (\CC)$ is a geometrically simple Abelian variety if and only if it is not $(n, n)$-decomposable for some $n>1$. 
\end{prop}

A more detailed discussion is given in \cite[Section 2.5]{lombardo}. The endomorphism rings of Abelian surfaces can be determined by the Albert's classification and results in \cite{Oort}. We summarize in the following:   

\begin{prop}  
The endomorphism ring $\End_{\overline \Q}^0  \, (\Jac \CC )$ of an abelian surface is either 
$\Q$, a real quadratic field, a CM field of degree 4, a non-split quaternion algebra over $\Q$, $F_1 \oplus F_2$,  where each $F_i$ is either $\Q$ or an imaginary quadratic field, the Mumford-Tate group $F$,  where $F$ is either $\Q$ or an imaginary quadratic field.   
\end{prop}

\begin{rem}  Genus 2 curves with extra involutions have endomorphism ring larger than $\Z$. 
Let $\CC$ be a   genus 2 curve defined over $\Q$. If $\Aut (\CC)$ is isomorphic to the Klein 4-group $V_4$, then  $\CC$ is isomorphic to a curve $\CC^\prime$ with equation 
\[ y^2 = f(x)= x^6 - a x^4+b x^2-1.\]
We denote  $u=a^3+b^3$ and $v=ab$. The    discriminant   
\[\Delta_f = - 2^6 \cdot \left(  27-18v+4u-u^2   \right)^2,\]  
is not a complete square in $\Q$ for any values of $a, b \in \Q$.  In this case $\Gal_{\Q} (f) $ has order 24.  There is a twist of this curve, namely $y^2=f(x)=x^6+a^\prime x^4+ b^\prime x^2 +1$, in which case $\Delta_f$ is a complete square in $\Q$ and $\Gal_{\Q} (f) $ has order 48. In both cases, from \cref{zarhin-thm} we have that   $\End_{\overline \Q} (\Jac \CC^\prime ) \neq  \Z$.
\end{rem}


Next we turn our attention to determining the endomorphism ring of abelian surfaces. Let us first recall a few facts on characteristic polynomials of Frobenius for abelian surfaces. The Weil $q$-polynomial arising in genus 2 have the form 
\begin{equation}\label{weil-g-2}
 f (T ) = T^4 - aT^3 + (b+2q) T^2 - aq T + q^2,
\end{equation}
for $a, b \in \Z$ satisfying the inequalities 
\[ 2 |a| \sqrt{q} - 4q \leq b \leq \frac 1 4 a^2 \leq 4q .\]
We follow the terminology from  \cite{bhls}.  Let $\CC$ be a curve of genus 2 over $\F_q$ and $\J=\Jac_ \CC$.  Let $f$ be the Weil polynomial of $J$ as in \cref{weil-g-2}. We have that $\# \CC (\F_q)=q+1-a$, $\#J (\F_q)=f(1)$ and it lies in the genus-2 Hasse interval 
\[ \H_q^{(2)} = \left[  (\sqrt{q}-1)^4, (\sqrt{q}+1)^4    \right] \]
In \cite{bhls} are constructed decomposable $(3, 3)$-jacobians with a given number of rational points by glueing two elliptic curves together.

Next we describe some of the results obtained in \cite{lombardo} for $\End_K (\AA)$ in terms of the characteristic polynomial of the Frobenius. We let $K$ be a number field and $M_K$ the set of norms of $K$.   Let $\AA$ be an abelian surface defined over $K$  and $f_v$ the characteristic Frobenius for every norm $v \in M_K$. 


\begin{lem}
Let $v$ be a place of characteristic $p$ such that $\AA$ has good reduction.  Then $\AA_v$ is ordinary if and only if the characteristic polynomial of the Frobenius 
\[ f_v (x) = x^4 + a x^3 + b x^2 + ap x + p^2, \]
satisfies $b \not\equiv 0 \mod p$.   
\end{lem}

Then from   \cite[Lemma~4.3]{lombardo}   we have  the following.

\begin{lem}
Let $\AA$ be an absolutely simple abelian surface. The endomorphism algebra $\End_{\bar K}^0 (\AA)$ is non-commutative (thus a division quaternion algebra) if and only if for every $v \in M_K$, the polynomial $f_v (x^{12})$ is a square in $\Z[x]$. 
\end{lem}

The following gives a condition for geometrically reducible abelian surfaces. 


\begin{prop}[\cite{lombardo}]
i) If $\AA/K$ is geometrically reducible then for all $v\in M_k$ for which $\AA$ has good reduction the polynomial $f_v (x^{12})$ is reducible in $\Z[x]$. 

ii) If $\CC$ is a  smooth, irreducible genus 2 curve with affine equation $y^2=f(x)$  such that $f(x) \in K[x]$ is an irreducible polynomial of degree 5   then $\Jac \CC$ is absolutely irreducible. 
\end{prop}

In \cite{lombardo} is given a detailed account of all the cases and an algorithm how to compute $\End_K \AA$.  

%
\renewcommand\H{\mathbb H}

\section{Modular curves}\label{sect-5}
As stated in \cref{sect-2} we are interested in isogenies between Abelian varieties. Of special interest is the case of elliptic curves, and it turns out that their isogenies are fairly well accessible both from the theoretical and algorithmic point of view.  The reason for this is is the very rich and well understood structure of \emph{modular curves}, which parametrize isogenies of elliptic curves.

\subsection{Modular curves over $\C$}\label{mod-curves}
As we have seen in \cref{sect-1} the isomorphism classes of elliptic curves $\E$ over $\C$ correspond one-to-one to isomorphism classes of lattices  $\Lambda_\tau$ with $\tau$ in $\H =\left \{ z = x+ i y \in \C \, \left | \frac{}{}\right. \, y > 0\right \} \subset \C$.
Moreover, $\Lambda_\tau\cong \Lambda_{\tau'}$ if and only if $\tau'=\frac{a\tau+b}{c\tau+d}$ with  $\begin{pmatrix} a & b \\ c& d  \end{pmatrix}$   in $SL_2(\Z)$.

So the moduli space of isomorphism classes of elliptic curves over $\C$ is in a natural way equal to $\H/\SL_2(\Z)$, which is, via the $j$-function, identified with $\mathbb{A}^1$.
It is more convenient to work with compact Riemann surfaces  and so, with projective curves.   Hence we compactify by adjoining \emph{cusps}.  

Define $\H^*=\H\cup \Q\cup \{ i\cdot\infty \}$ and extend the action of $\SL_2(\Z)$ in an obvious way (e.g.: $1/z$  interchanges the cusp $0$ with $i\cdot \infty$).  Then $\H^*/\SL_2(\Z)=\P^1$, the $j$-function extends to a meromorphic function on $\H^*$ with a simple pole at the cusps (which are all equivalent) and all points  $P \in \P^1 \setminus \{ j (i\cdot \infty) \}$ have a \emph{modular interpretation}: to $P=j(\tau)$ there corresponds the isomorphism class of the elliptic curve $\E_\tau$ with lattice $\Z+\tau \Z$. 
In this interpretation we call $\P^1 = X(1)$ and $\A^1=Y(1)$ modular curves of level $1$.

We introduce now a special family of congruence subgroups of $\SL_2(\Z)$, which are linked to isogenies. For $N\in \N$ define 
\[
\begin{aligned}
\Gamma_0 (N)  & := \left\{ \begin{pmatrix} a & b\\ c & d \end{pmatrix} \in SL(2, \Z) \;  | \; c \equiv  0 \mod N     \right\} \\
\end{aligned}
\]
This is a subgroup of $SL_2(\Z)$  of finite index, and so $\H^*/\Gamma_0 (N)$ is a compact Riemann surface $X_0(N)$ with a natural cover morphism
\[ \eta_N : X_0(N) \ra X(1) \]
of degree $\varphi(N)$, the value of the Euler function of $N$.   As affine part in $X_0(N)$ we find $Y_0(N)=X_0(N)\setminus \eta^{-1}(j(i\cdot \infty))$ and points $P$ in $Y_0(N)$ have the following  modular interpretation.

Let $P\in Y(1)$ corresponding to an isomorphism class of the lattice $\Lambda_\tau=\Z+\tau \Z$. The inverse image  $\eta^{-1}(P)$ consists of the equivalence classes of lattices isomorphic to $\Lambda_\tau$ \emph{with the additional information} consisting of a lattice $\Lambda_\tau(N)\supset\Lambda_\tau$ with
$\Lambda_\tau(N)/\Lambda_\tau\cong \Z/N$. Hence a point $P_N\in Y_0(N)$ is the $\C$-isomorphism class of a \emph{pair} $(\E_\tau,C_N)$ where $C_N$ is a cyclic group of order $N$ in $\E[N]$. 

In other words, $Y_0(N)$ parametrizes \emph{isomorphism classes of elliptic curves together with cyclic isogenies $\eta_N$ of order $N$}
and so they are moduli spaces for the pairs  $(\E,\eta_N)$ over $\C$.
Using Hurwitz genus formula and the well known fixed points of $\Gamma_0(N)$ one can compute the genus $g_N$ of $X_0(N)$.   We have that $g_N\sim N$. 
For $N=p$ a prime we get
\[
g_N={
\left \{
\begin{aligned}
& 0  \quad if  \quad  p=2,3 \\
& {\frac {(p-13)} {12}} \quad  if  \quad   p\equiv 1 \mod 12        \\
& {\frac {(p-5)} {12} }  \quad if   \quad  p\equiv 5 \mod 12  \\
& {\frac {(p-7)} {12}}  \quad if   \quad  p\equiv 7  \mod 12  \\
& {\frac {(p-11)} { 12}}  \quad if   \quad  p\equiv 11  \mod 12  \\
\end{aligned} 
\right.
}
\]

\subsubsection{Modular functions and forms}
Because of the construction as quotient of $\H^*$ we can identify algebraic-geometric objects on $X_0$ with analytic objects attached to
$\H^*$ with specific symmetry properties.  For instance, rational functions on $X_0 (N)$ come from those meromorphic functions on $\H^*$  that 
are invariant under the action of $\Gamma_0(N)$, ie.
\[f \left( \frac{az+b}{cz+d} \right)=f(z) \; \mbox{ for all} \;  \begin{pmatrix} a & b\\ c & d \end{pmatrix} \in \Gamma_0(N)\]
and are called \emph{modular functions} of level $N$. A basic example is $j(z)$, which is modular of level $1$ and hence of all levels.
%
%
Since for all $\tau$ the elliptic curve the lattice $\Z+\frac{\tau}{N} \, \Z$ gives rise to an isogeny of $\E_\tau$ with cyclic kernel of order $N$ the function $j \left(\frac{z}{N}\right)=:j_n(z)$ is a modular function of level $N$, and one checks that $\C (X_0(N) )$ is isomorphic to $\C ( j(z), j (\frac{z}{N}) )$.

Differentials $\omega$ on $X_0(N)$ are of the form $f(z)dz$ with $f(z)$ meromorphic on $\H^*$ and, because of their invariance under $\Gamma_0(N)$,
they  satisfy a functional equation
\[f \left( \frac{az+b}{cz+d} \right)=(cz+d)^2f(z) \; \mbox{ for all} \; \begin{pmatrix} a & b\\ c & d \end{pmatrix} \in \Gamma_0(N).\]
Such functions $f$ are called \emph{modular forms of weight $2$}.  If $\omega$ is holomorphic then $f(z)$ is holomorphic on $\H$ and vanishes on the cusps: $f(z)$ is a cusp form of weight $2$ and level $N$.  We know from the Riemann-Roch theorem that the cusp forms $\mathcal{S}_0(N)(\C)$  are a $\C$-vector space of dimension $g_N$. This  space plays a very important role in the arithmetic of modular forms.

\subsubsection{$q$-expansion}
Modular forms and functions are invariant under the transformation $z\mapsto z+1$.  Define as new variable
 $q=e^{2\pi iz}$ \, (hence transform neighborhoods of $i\cdot\infty$ to  neighborhoods of $0$). Then we can form the Fourier series of modular forms and functions $f(q)$ and get the  \textbf{q-expansion} 
\[ f(q) = \sum_{j=k}^\infty a_j q^j \]
 with $k\in \Z$. If $f$ is a cusp form then $k\geq 1$. Modular forms and functions are uniquely determined by their Fourier expansion.

\subsection{Modular Polynomials}\label{mod-pol-ell}
The important fact is that in many cases one can compute the coefficients of the $q$-expansion effectively. This is so for $j(q)$ (and hence for $j_N$ in the variable $q^{1/N}$). In particular, $j(q)$ has a pole of order $1$ at $0$ and has coefficients in $\Z$, and for the exact expression see   \cite[Proposition 5.41]{book}.

We know that $j_N$ satisfies a polynomial over $\C(j)$ of degree $\varphi(N)$. To find this polynomial (and so an equation for $X_0(N)$) one compares sufficiently many coefficients of the Fourier expansion and finds the \emph{modular polynomials} $\phi(j,j_N)$ which are monic and symmetric in $(j,j_N)$ and have degree $\varphi(N)$ in $j$ and $j_N$.

So two  elliptic curves with $j$-invariants $j_1$ and $j_2$ are isogenous under a cyclic isogeny of degree $N$ if and only if $\phi_N (j_1,j_2)=0$.    The equation $\phi_N (X, Y)=0$ is the canonical equation of the affine modular curve $Y_0 (N)$, and $X_0(N)$ is given by the closure in $\P^2$.  We display $\phi_N (x,y)$  for $N=2, 3$. 
\begin{small}
\[
\begin{split}
\phi_2 & =x^3-x^2y^2+y^3+1488xy(x+y)+40773375xy-162000(x^2+y^2)\\
& +8748000000(x+y) -157464000000000 \\
\phi_3  & = - x^{3}  y^{3}+2232 x^{3} y^{2}+2232 y^{3} x^{2}+ x^{4}- 1069956 x^{3}y+2587918086 x^{2} y^{2} \\
& -1069956 y^{3}x+ y^{4} +36864000 x^{3} +8900222976000 x^{2}y +8900222976000 y^{2}x \\
& + 36864000 y^{3}  +452984832000000 x^{2}-770845966336000000 xy+ 452984832000000 y^{2} \\
& +1855425871872000000000 x+ 1855425871872000000000 y  \\
\end{split}
\]
\end{small}
There are tables of modular polynomials for large $N$.  

\begin{rem}
There have been some attempts in the last decade to generalize modular polynomials for higher dimensional varieties. The interested reader should consult \cite{lauter} for abelian surfaces.
\end{rem} 

\subsection{The arithmetical theory of modular curves}
A general reference for the following discussions is \cite{mazur}. One observes that the modular polynomials have coefficients in  $\Z$. This is no accident.

Looking at the modular interpretation one sees that the curves $X_0(N)$ represent over $\C$ a moduli problem: "Parametrize the isomorphism classes of pairs $(\E,\eta_N)$ of elliptic curves with cyclic isogeny  of degree $N$".

This problem makes sense over any ring $R$ with invertible element $N$ and so over $\Z[1/N]$.  Because of the existence of twists $\chi$ and the fact that the pair $(\chi(\E),\chi(\eta_N))$ is isomorphic to $(\E,\eta_N)$ over $\bar{k}$ but in general not over $k$ we cannot expect to find a fine moduli scheme for our problem. But a coarse moduli scheme exists and so $\H^*/\Gamma_0(N)$ is the constant field extension of a curve $X_0(N)$ (i.e a scheme of relative dimension $1$) over $\Z[1/N]$ with affine part given by  the equation $\phi_N(X,Y)$.

A deep analysis  of reductions of modulo prime numbers dividing  $N$ due to Deligne-Rapoport (\cite{De-Ra})  and Katz-Mazur (\cite{Ka-Ma}) shows that  one can define $X_0(N)$ over $\Z$. Hence all geometric objects like holomorphic differentials (and so cusp forms) have a $\Z$-structure and so it makes sense to speak of $\mathcal{S}_0(N)(R)$, the $R$-module of cusp forms over the commutative ring $R$ (to avoid complications assume that $N$ is invertible in $R$).

How to find these cusp forms?  Following a ground breaking idea of \textbf{J. Tate} one replaces the theory of elliptic curves over $\C$ by  the theory  of elliptic curves in the realm of rigid $p$-adic spaces and finds the elliptic Tate curve in $\Z((q))$, the  ring of  formal Laurent series in one variable in $q$.  Using this curve and exploiting functions and differentials on this curve one finds the  $q$-expansion over all rings $R$ and gets for instance:

\begin{thm}[\textbf{$q$-expansion principle}]
Let $R$ be a commutative ring with $N\in R^*$. Then the $q$-expansion determines uniquely holomorphic differentials on $X_0(N)$, and in particular we have:

\begin{enumerate}
\item $f\in \mathcal{S}_0(N)(R)$ if and only if the Fourier expansion of $f$ has coefficients in $R$, and

\item there exists a $\Z$-base $f_1,\dots,f_{g_N}$ of $\mathcal{S} _0(N)(\Z)$ such that 
\[ \mathcal{S}_0(N)(R)= \<f_1,\dots,f_{g_N} \>\otimes R.\]
\end{enumerate}
\end{thm}

We can apply these results to the Jacobian variety  and its minimal model over $\Z$, denoted by $\J_0(N)_\Z$. Endomorphism of this variety are uniquely determined by their action on $\mathcal{S} _0(N)(\C)$ and so accessible to computation after we have computed a base of the cusp forms over $\Z$. This computation is done in an effective way by using \emph{modular symbols}; see \cite{Me}.   
It turns out that $\End_\Z( \J_0(N)_\Z$ is large, it contains the Hecke algebra $\mathcal{T}_N$ generated by Hecke operators constructed in \cref{sect-correspondences}. In particular, it follows that simple factors of  $\J_0(N)$ have real multiplication.

Galois representations of $G_\Q$ attached to the Tate modules of these factors split up in a sum of $2$-dimensional representations. One of the greatest results in Arithmetic Geometry is the theorem of Khare-Kisin-Wintenberger (see \cite{Kh-W}, \cite{Kisin})  which confirms the conjecture of Serre and states that all odd two-dimensional  Galois representations are obtained (up to twists by characters) by the operations on appropriate factors of $\J_0(N)$  with an explicit recipe how to find $N$ and the twist character. From this, the Fermat conjecture FLT is an easy consequence (a "five line proof").

\part{Cryptography}\label{part-2}
 For communication in the modern world it is crucial that one can use open communication channels to
\begin{itemize}
\item  {exchange keys,}
\item  {sign messages}
\item authenticate entities, and
\item encrypt and decrypt (not too large) messages
\end{itemize}
with simple protocols, clear and easy to follow implementation rules based on  \emph{cryptographic primitives}, which rely on (hopefully) hard mathematical tasks. The part of cryptology, which is devoted to solve these challenges is \emph{public key cryptography} and relies on the ground breaking ideas of Diffie and Hellman \cite{DH}.
In this paper we shall concentrate to the first item in its simplest form (and not using a secure protocol). The reader is encouraged to look at the other aspects, too, e.g. in \cite{book}. 

\section{Diffie-Hellman Key Exchange}\label{Diffie-Hellman}

\subsection{The Classical Case}
The task to solve is: find a protocol such that two partners $P_1$ and $P_2$ can agree on a common secret by using public channels and algorithms.

A groundbreaking solution was found by W. Diffie and H. E. Hellman in \cite{DH} with the idea to use (computational) one-way functions.  They suggest to use the multiplicative group of finite fields $\Fq$ and, for a prime number  $\ell | q-1$ choose a primitive root $\zeta_\ell$, private keys $k_i\in \Z$ and public keys $p_i=\zeta_\ell^{k_i}$. The common secret is 
\[s_{1,2}=\zeta_\ell ^{k_1\cdot k_2}.\]
All computations are very fast (polynomial in $\log (q)$).    The security is measured by the hardness of the \textbf{Diffie-Hellman} computational problem (CDH):
\[ \mbox{For random elements }  a,   b\in \{0, \dots,    \ell-1\} \mbox{ and given }  \zeta_\ell^a,\,\,\zeta_\ell^b\mbox{ compute } \zeta_\ell^{a\cdot b}.\]
Let $\zeta$ be a primitive root of unity in $\Fq^*$.  Define the (classical) discrete logarithm (DL) of an element $x\in \Fq^*$ with respect to the base $\zeta$ by 
\[ \log_\zeta(x)=  \min \left\{    n\in \N \mbox{ such that } \zeta^n=x.   \right\} \]
It is obvious that an algorithm that computes discrete logarithms (e.g. in $\zeta_\ell$) solves (CDH).   This problem is rather old (going back at least to the 19-th century). C.F. Gauss  introduced the term "index" in the Disquisitiones Arithmeticae (1801) for the discrete logarithm modulo $p$, and there are tables for primes  up to $1000$ by C. G. Jacobi (1839).

A systematic algorithm is given in the book of Kraichik (1922) \cite{Kra}; in fact this is the index-calculus algorithm reinvented and refined in cryptography  from 1980  till today \cite{JOP}.  As result one gets algorithms of subexponential complexity (with relatively small constants, see \cite{JOP}), which are even dramatically faster if $q$ is not a prime. 
%
%
The reason for these fast algorithms is the fact that it is easy to lift elements in $\Fq$ to elements in rings of integers of number fields.

\subsection{A First Abstraction}
Obviously we can use the Diffie-Hellman key exchange scheme if we have
\begin{itemize}
\item a finite cyclic group $(C,\circ)$ with a generator $g_0$,

\item a numeration, i.e. an injective map 
\[ f : C \ra \N, \]
\item an addition law $\oplus$ on $f(C)$ with
\[ f(f^{-1}(a)\circ f^{-1}(b)) =a\oplus b\mbox{ for all }a,b\in f(C).\]
\end{itemize}
$f(C)$ becomes a $\Z$-module with the usual scalar multiplication: $0\cdot a= f(0_C)$, $n\cdot a=(n-1)$-fold addition of $a$ to itself, $(-n)\cdot a=n\cdot(\ominus a)$ for $a\in f(C),\,n\in \N$.

The private keys are again $k_i \in \Z$, the public keys are $k_i \cdot f(g_0)$, and the common secret is  $k_1\cdot(k_2\cdot f(g_0))$. The CDH problem is: for random $a_1$, $a_2\in f(C)$ with publicly unknown $k_1, k_2$ such that $k_i\cdot f(g_0)=a_i$ compute $c=(k_1\cdot k_2)\cdot f(g_0)$.  

Define the  discrete logarithm (DL) by 
\[\log_{g_0}(a) := \min \, \{n\in \N \, \mid \,   \text{ such that } n\cdot f(g_0)=a.  \} \] 
Again, the computation of the (DL) solves CDH.  By elementary number theory (CRT and p-adic expansion) one sees immediately that the computation of (DL) is reduced to the computation of the discrete logarithms in all $f(C_\ell)$ with $C_\ell$ the subgroup of $C$ of elements of order dividing $\ell$ and $\ell$ dividing $|C|$.  Hence we can and will assume from now on that $C$ is cyclic of prime order $\ell$.
For shortness, we denote the task to compute the discrete logarithm  by (DLP). 
%
\subsubsection{Black Box Groups}
A "generic" object of the situation above is given by a black box group $C$ of prime order $\ell$.
\begin{enumerate}
\item There are algorithms   that compute (DL) (probabilistically) with 
$\O(\sqrt{\ell})$
group operations  in $C$ (e.g. Shank's bay-step giant step algorithm, 
Pollard's
%
 $\rho$ algorithm et.al.),  these algorithms are applicable for all finite cyclic groups, and one cannot do better.

\item Up to algorithms with subexponential complexity, the computation of (DL) in $C$ is equivalent with (CDH)(Maurer-Wolf).
\end{enumerate}
%
\subsection{Mathematical Task}\label{group}
In order that we can use  (a family of) groups $C$ for crypto systems based on discrete logarithms  they have to satisfy  four crucial conditions:
\begin{enumerate}
\item $C$ has a known large prime order $\ell$ and a numeration $f:  C\ra \N$. 
\item Condition for the numeration: The elements in $C$ can be stored in a computer in a {compact} way (e.g. $\O (\log \ell)$ bits needed).
\item The group composition $\oplus$ induced by $f$ is given by an algorithm  that is easily and efficiently implemented and  very {fast}.
\item The computation of the DL in $f(C) $ (for random elements) is  {very hard} and so infeasible in practice (ideally the bit-complexity should be exponential in $\log\ell$).
\end{enumerate}
It is surprisingly hard to construct such groups. All known examples today are related with subgroups of Picard groups of hyperelliptic curves of genus $\leq 3$ over prime fields $\F_p$.   It will be one of the main aims of the paper to explain this statement.
%
\subsection{Q-bit Security}
As said, we shall describe below DL-systems  for which we have good reasons to believe that the bit-complexity is exponential and so the task in  \cref{group} is solved.  But the possibility that 
\emph{quantum computing} may be realizable in foreseeable time yields new aspects for the discussion of security of crypto systems.  By Shor's algorithm it follows  that the q-bit complexity of discrete logarithms in \textbf{all} finite groups is  polynomial!. 

So it is challenging to find key exchange systems that are not based on discrete logarithms in groups but still are near to the original idea of Diffie and Hellman. 
In the quantum world new relations between crypto primitives arise, and it seems that hidden subgroup problem  and connected to it, the hidden shift problem related to groups $G$ are central (\cite{Reg} and \cite{Kup}).  Here the state of the art is that for abelian groups $G$ the problems can be solved in subexponential time and space, for dihedral groups there is  "hope". 
%
\subsection{Key Exchange with $G$-sets} 
The DL-system in  \cref{group} can be seen in the following way:  By scalar multiplication 
$f(C)$ becomes a $\Z$-set, and so elements of $\Z$ induce commuting endomorphisms on $f(C)$.

Denote by $\Z'$ the semigroup of elements in $\Z$ prime to $|f(C)|$. Then the set $A$ of generators of $f(C)$ becomes a $\Z'$-set,
and elements in $\Z'$ induce commuting endomorphisms of $A$.
%
A next step to generalize the Diffie-Hellman key exchange is to replace  $\Z'$
 by a (semi-) group $ G $ and the set of generators  
of $f(C)$
by a $G$-set $ A \subset \N$ on which $ G $ operates  transitively.  For $g \in G $, define $t_g \in \End_{set}(A)$ by 
\[ a \mapsto t_g(a):=g \cdot a. \]
Let $G_1$ be a semi-subgroup of $G$ and $G_2=Z(G_1)$ the centralizer of $G_1$ in $G$   (if $G$ is abelian then $G=G_1=G_2$).
Because of 
\[ g_1 \cdot (g_2 \cdot a_0) = (t_{g_1} \circ t_{g_2}) \cdot a_0 = (t_{g_2}\circ t_{g_1}) \cdot a_0 \]
for  $g_1\in G_1$ and $g_2\in G_2$ we can use $(A, a_0, G_1, G_2)$ for key exchange by defining an obvious analogue of the scheme in   \cref{group}.

The security of this exchange depends on the difficulty to find the translations $t_{g_i}$.  We remark that though the security of such systems is, in general, not related to discrete logarithms, it may happen that the generic algorithms from \cref{group} can still be applied. 

What about quantum security? One breaks the system if one can determine $t_{g_1}$.  This is a typical problem for the hidden shift.  Take the maps
\[
\begin{split}
f_0 :  G_1    &  \ra A,   \text{ such that  }      f_0 (g)=t_g \cdot a_0 \\
f_1 : G_1     &   \ra A,  \text{ such that  }      f_1 (g)=t_g \cdot(t_{g_1} \cdot a_0) \\
\end{split}
\]
and find the shift.   For $G_1$ abelian and finite there is an algorithm of Kuperberg \cite{Kup}, which solves this task in subexponential time.   In particular we see that every Diffie-Hellman key exchange based on $\Z$-sets has at best subexponential security.

\subsection{Abstract Setting of Key Exchange}

On our way to generalization we get rid of the  algebraic structures.   Assume  $A\subset\mathbb N$ and let  $B_1, B_2 \subset \End_{set}(A)$. Choose $a_0\in A$.
We need the \textbf{centralizing condition}.   The elements of $B_1$ {commute} with the elements of $B_2$   on $B_i\{a_0\}$.  Then 
\[  b_1 (b_2 (a_0))= b_2 (b_1(a_0)) \]
and this is all we need for key exchange. 

The effectiveness of this exchange is given if for $b_i\in B_i, b_j\in B_j$  the value $b_i(b_j(a_0))$ can be {quickly evaluated} (i.e., calculated and represented). The analogue of the Computational Diffie-Hellman problem is 
\begin{center}
 \mbox{\textbf{CDH}}: For randomly given $a_1, a_2 \in A$,      compute  (if exists)  $a_3$   with  $a_3 = b_{a_1} \cdot ( b_{a_2} \cdot a_0)$, 
\end{center} 
where $b_{a_i}\in B_i$ such that  $b_{a_i}\cdot a_0=a_i$. It is clear that {CDH} can be solved if one  can calculate {for random $ a\in B_i\cdot\{a_0\}$} an endomorphism  $ b_a\in B_i $ with $b_a(a_0)=a$. 
We remark that $b_a$ may  not be uniquely determined by $a$.

\medskip

\noindent \textbf{Problem:} 
\begin{enumerate}
\item Find a "genuine" usable instance for the abstract setting!
\item What can one say about quantum computing security?
\end{enumerate}

\subsection{Key Exchange in Categories}
We make a final step of abstraction.
As always we assume that we have two partners $P_1$ and $P_2$  who want to have a common secret key.

Let $\CC_i$, $i=1, 2$ be two categories whose objects are the same  sets $A_j$ and with morphisms $B^i_{j,k}=\Mor^i(A_j,A_k)$.  We fix a "base" object $A_0$ and assume that $\CC_1,  \CC_2,  A_0$ satisfies the following conditions:
\begin{enumerate}
\item For every $\varphi\in B^1(A_0,A_j)$ and every $\psi\in B^2(A_0, A_k)$  the  pushout exists,  i.e. there is a uniquely (up to isomorphisms) determined triple 
\[(A_l,\,\, \gamma_1\in B^1(A_k,A_l),\,\,\gamma_2\in B^2(A_j, A_l))\]
with 
\[ \gamma_2\circ \varphi=\gamma_1\circ \psi \]
such that 
%
 this triple is minimal (universality condition).

\item $P_1$ can determine $A_l$ if he  knows $\varphi$, $A_k$ and an additional (publicly known) information $P (\psi)$ (which is often a subset of $A_k$), and an analogue  fact holds for $P_2$. \\
\end{enumerate}

\paragraph{ \textbf{Key Exchange}}  Given such categories $\CC_1,  \CC_2$ the partners can chose $\varphi,  \psi$, send $A_j$, $A_k$ and $P(\psi)$ respectively $P(\varphi)$ and compute the \textbf{\emph{common secret}} $A_l$. \\

\paragraph{ \textbf{Effectiveness}} We assume that all the objects concerning $\CC_i$ can be handled by computers in a fast and compact way, in particular, for chosen $\varphi$, $\psi$  the objects $A_j$, $A_k$ as well as the additional information can be computed rapidly.  Moreover, using the given information, $P^i$ can compute of $A_l$  quickly. \\

\paragraph{\textbf{Security}} The scheme is broken if (CDH) is weak: For randomly given $A_j$, $A_k$ determine $A_l$, which is the pushout of 
\[ 
A_0\stackrel{\varphi}{\ra}A_j
\]
and 
\[ A_0\stackrel{\psi}{\ra}A_k.
\]
For this, it is allowed to use the additional information. We shall see an example for this categorial key exchange in \cref{supersingular}, and till now all algorithms for breaking this system have exponential complexity.

\section{Index calculus in Picard groups}\label{index-calculus}
We want to use systems based on discrete logarithms in groups $G$ and so find groups which satisfy the conditions formulated in  \cref{group}. Motivated by ideas of V. Miller and N. Koblitz we want to use subgroups of Picard groups of curves over finite fields.    \cref{add} of  Hess-Diem  implies that at least in principle for such groups the conditions 2 and 3 are satisfied.
For finding subgroups of large prime order one has to be able to determine the order of Picard groups rapidly. Here the key word is point counting, and again there is, in principle,  a solution by a polynomial time algorithm due to Pila generalizing the Schoof algorithm for elliptic curves. But these algorithms are much too slow, and an acceleration is only known for elliptic curves (AES-algorithm), for curves of genus $2$ (Gaudry, Schost) and for curves with special endomorphism rings (complex multiplication or real multiplication). 

But before investing a lot of work in point counting it is useful to look at the security aspect. We want to compare the hardness  of the computation of the DL
in the specific groups with the generic hardness, i.e. $\sim |G|^{1/2}$. We recall that a main reason against the classical DL was the index-calculus algorithm, which is based on the (easy) lifting of finite fields to integers in number fields or function fields over finite fields. This kind of attack is not possible in Picard groups of curves of positive genus as pointed out by Miller and Koblitz:  The "golden shield" of the N{\'eron-Tate quadratic form prevents a (easy) lifting of elements in Abelian varieties over finite fields to  global fields.   But unfortunately there are very effective variants of the index- calculus attack to Picard groups. 

\subsection{Introduction to index calculus}
Let $(G,\oplus)$  be a cyclic group of order $N$ with generator $g_0$.

\textbf{First step:}\;  Find a "\emph{factor base}" consisting of relatively few elements and  compute $G$ as $\mathbb Z-$module given by the free abelian group
generated by the base elements modulo relations. 

So choose a subset $\mathcal{B}= \{g_1,\dots,g_r\} $ of $G$  generating $G$ and look for relations  
\begin{equation}\label{index_calc} 
R_j:\oplus_{i=1}^r [n_i] g_i =0_G.
\end{equation}
Obviously $R_j$ yields the relation 
\begin{equation}\label{system} 
\sum_{i=1}^r n_i \log_{g_0} (g_i) \equiv 0 \mod N
\end{equation}
for discrete logarithm.

We assume that we can find sufficiently many independent relations  as in Eq.~\eqref{index_calc}  for solving the system in Eq.~\eqref{system} via linear algebra for $\log_g g_i$, $i=1, \dots , r$. Then we have an explicit presentation of $G$ as $\Z$-module by
\[G\cong\Z^r/ \< \dots ,R_j, \dots \>.\]

\textbf{Second step:} \;  Take $g\in G$ randomly and chose a "random walk" with steps $g^0=g,\dots, g^j=[k_j]g^{j-1}$ and assume that after a few steps $j$
we find a tuple $e_1, \dots ,e_r$ with $e_i$ small and 
\[ g^j=[e_1]g_1+\cdots [e_r]g_r.\]
"To find" means: There is a fast algorithm to decide whether such $e_i$ exist, and then the computation of these $e_i$ is also fast. 
This boils down to a smoothness condition. (Recall: A number $n \in \N$ is $B$-smooth if all prime divisors of $n$ are $\leq B$, and results from analytic number theory by Canfield, Erd\"os, Pomerance  determine  the probability for $n$ being smooth.) The second step is usually done by an appropriate sieving method.  

The important task in this method is to balance the number of elements in the factor base to make the linear algebra over $\Z$  manageable and to guarantee "smoothness" of arbitrary elements
with respect to this base. Usually one finds a kind of \emph{size} in $G$ (size of lifted elements in $\Z$ or degree in polynomial rings, degree of reduced divisors , \dots) to define factor bases. Typically, successful index-calculus approaches give rise to algorithms for the computation of the DL in $G$  which have {\emph{subexponential}} complexity and so, for large enough order of $G$, the DL-system  has a poor
security.

For an axiomatic approach of index-calculus algorithms we refer to \cite{EG}. This principle is refined in concrete situations with enormous effect as we shall see below.

\subsection{Index calculus for hyperelliptic Jacobians}
Index calculus can be applied to a DL in Jacobians of hyperelliptic curves.   Let $\CC$ be a hyperelliptic curve of genus $g\geq 2$ over a finite field $\F_q$ of characteristic $p$ and 
$G$ a cyclic subgroup in $\Pic^0_{\CC}$.


We can represent every element in $G$ in a unique way by the Mumford representation $[u (x), v (x)]$, where $u(x)$ is a polynomial in $\Fq[X]$ of degree $\leq g$. 
%
%
%
%
%

%
As factor base we choose points in $\Pic^0_\CC$ with $u(X)$ irreducible of degree at most $B$, a chosen smoothness bound. A divisor is said to be \textit{$B$-smooth} if all the prime divisors in its decomposition have degree at most $B$.  This leads to the historically first algorithm to compute discrete logarithms in Picard groups of hyperelliptic curves. It is due to  Adleman, Demarrais, and Huang \cite{adh}.   For an explicit description of the algorithm see \cite[pg. 525]{book}. For $N\in \Z^{>0}$, $s, c \in \R$, with $0 \leq s \leq 1$ denote 
\[ L_N (s, c) = \exp \left( (c + o(1) \right) \left( \log N)^s \left( \log \log N\right)^{1-s} \right), \]
as $N \to \infty$. Then we have 

\begin{thm}
For $\log q \leq (2g+1)^{1-\epsilon}$,  there exists a constant $c \leq 2.18$ such that the discrete logarithms in $\Jac_{\CC} (\Fq)$ can be computed in expected time $L_{q^{2g+1}} (1/2, c)$. 
\end{thm}
%
%
%
 This remarkable result  gives an subexponential  algorithm for "large" genus. 
%
 %
But much more important for practical applications are \emph{exponential} algorithms, which weaken the DLP for small but realistic genus.
The first groundbreaking result is

\begin{thm}[Gaudry]
Let $\CC$ be a  a hyperelliptic curve of genus $g\geq 2$ defined over a finite field $\F_q$. If $q > g!$ then discrete logarithms in $\Jac_{\F_q} (\CC)$ can be computed in expected time $O (g^3 q^{2+\epsilon})$.
\end{thm}
Since the expected size of $\J_\CC(\Fq)$ is $q^g$ (see Weil's result,  \cref{Hasse-Weil})  we are, for $g>4$, far away from the generic security bound, and so we have to exclude hyperelliptic curves of genus $\geq 5$ if we want a DL-system in Picard groups.
But Gaudry's  result can be sharpened. {N. Th\'{e}riault} suggested to use "large primes" as well as the original elements of the factor base consisting  of points on the curve of small degree.  With many more refinements \cite{DGTT} one gets

\begin{thm}\label{double}
There exists a (probabilistic) algorithm which computes the DL, up to $\log$-factors, in the divisor class group of hyperelliptic curves of genus $g$
in expected time of $\O(q^{(2-2/g)})$.
\end{thm}
This rules out $g=4$ for hyperelliptic curves since the ratio of the expected group order to time complexity, $\O (q^g)/\O (q^{2- 2/g})$, gets too big. 

\subsection{Index-calculus in Picard groups in curves with plane models of small degree}
The following is mainly work of C. Diem. He  gives an algorithm for computing discrete logarithms in $J_\CC(\Fq)$ assuming that one has a plane curve $\CC'$ of degree $d$. We recall that  for non-hyperelliptic curves $d=2g_\CC-2$ is possible but that for hyperelliptic  curves $d\geq g_\CC+1$.

So the minimal degree of plane models of hyperelliptic curves of genus $\geq 3$ is larger than the degree of such models for non-hyperelliptic curves.   Using factor bases constructed with the help of Semaev polynomials and using a large amount of ingredients from abstract algebraic geometry (e.g.  membership tests for zero-dimensional schemes) Diem succeeds to prove:
\begin{thm} Fix $d \geq 4$.    Then the DLP in $\Pic^0_\CC$ of curves birationally equivalent to plane curves of degree $d$ can be solved, up to $\log$-factors, in expected time $\O(q^{2- \frac{2}{d-2}})$.
\end{thm}

For genus $4$ and non-hyperelliptic curve $\CC$ we get $d=6$ and so the hardness of $D$ is bounded, up to $\log$-factors,  by $\O(q^{3/2} )$. Since the expected group size is $q^4$ this is too far away from the generic complexity, and it is not advisable to use (hyperelliptic or not hyperelliptic) curves of genus $4$ for DL-systems.

For non-hyperelliptic curves of genus 3 we get $d=4$ and so the complexity of the DL is $\O(g)$ and again such curves can not be used for DL-systems.  
Hence our discussion for fulfilling conditions $3$ and $4$ in   \cref{group} can be restricted to hyperelliptic and elliptic curves of genus $1, 2, 3$. Before doing this in detail we have one more general section, interesting both from theoretical and practical point of view. 

\section{Isogenies of Jacobians via correspondences and applications to discrete logarithms}\label{sect-correspondences}
We describe  a general construction of isogenies between abelian varieties closely attached to Jacobians of curves.  The crypto-graphical relevance of these constructions is that every computable isogeny yields a transfer of the (DLP), and it may be easier to  solve the problem after the application of the isogeny.

As always, $k$ is assumed to be a perfect field. Let  $L$ be a finite algebraic extension field of $k$. Let $\mathcal{D}_1$ be a regular projective curve over $L$ and $\mathcal{D}_2$ a regular projective curve defined over $k$.   We recall some properties of cover morphisms of curves and attached norm and conorm homomorphisms of Jacobians. Let $\mathcal{H}$ be a curve over $L$ and 
\[ \varphi_1:\mathcal{H}\ra \mathcal{D}_1, \]
respectively 
\[ \varphi_2: \mathcal{H}\ra\mathcal{D}_2\times _{\Spec(k)}\Spec (L)=:\mathcal{D}_{2,L},\]
be  $L$-rational morphisms.   The morphism $\varphi_1$ induces the $L$-rational   \textbf{conorm  morphism}
\[ \varphi_{1}^*: \mathcal{J}_{\mathcal{D}_1}\ra\mathcal{J}_\mathcal{H}\]
and the morphism $\varphi_1$ induces the   \textbf{norm morphism}
\[ \varphi_{2,*}:\mathcal{J}_\mathcal{H}\ra \mathcal{J}_{\mathcal{D}_{2,L}}.\]
By composition we get a homomorphism
\[\eta_L: \mathcal{J}_{ \mathcal{D}_{1}} \ra \mathcal{J}_{\mathcal{D}_{2,L}}\]
defined over $L$.

Let $\mathcal{W}_{L/k}$ be the Weil restriction of the Jacobian of $\mathcal{D}_1$ to $k$. This is an abelian variety defined  over $k$ with $\mathcal{W}_{L/k}(k)=\Pic^0_{\mathcal{D}_1}$. Applying the norm map from $L$ to $k$ and using the functorial properties of the Weil restriction we get a homomorphism
\[\eta:  \mathcal{W}_{L/k}\ra \mathcal{J}_{\mathcal{D}_2}.\]
In general, neither the kernel nor the cokernel of $\eta$ will be finite. But under, usually mild, conditions one can assure that that  $\eta$  has a finite kernel, and so it induces an isogeny of  $\mathcal{W}_{L/k}$   to an abelian subvariety of $\mathcal{J}_{\mathcal{D}_2}$.

As application we get a transfer of the discrete logarithm problem from $\Pic^0_{\mathcal{D}_1}$ (defined over $L$) to the DL-problem in a subvariety of $\mathcal{J}_{\mathcal{D}_2}$ (defined over $k$). Of course, the efficiency of this transfer depends on the complexity of the algorithms computing the norm and  conorm  maps (hence $\varphi_i$ and $[L:k]$ must have  reasonably small degrees), and an attack makes sense only if the DL-problem after the transfer is easier than before.

\subsection{Weil Descent}
Take $k=\Fq$ and $L=\mathbb F_{q^d}$ with $d>1$ and  $\mathcal{H} = \mathcal{D}_{2,L}$, i.e.  a given curve  $\CC$ defined over $\F_{q{^d}}$ is covered by a curve $\mathcal{D}_{\F_{q{^d}}}$, which is the scalar extension of a curve $\mathcal{D}$  defined over $k$.

This yields a $k$-rational homomorphism from the Weil restriction $\mathcal{W}_{L/k}$ of $\mathcal{J}_\CC$ to $\mathcal{J}_\mathcal{D}$.  Then $\mathcal{D}$ will (in all non-trivial cases) be a curve of a genus larger than the genus of $\CC$ but since it is defined over the smaller field $\Fq$  one can hope that one can apply fast algorithms to compute the discrete logarithm  in   $\mathcal{J}_\mathcal{D}(\Fq)$, e.g. by methods of index-calculus in \cref{index-calculus}.  Indeed, if $\CC$ is not defined over a proper subfield of $\F_{q^d}$ this is the principle of the so-called GHS-attack in (see \cite{GHS} and \cite[Section 22.3.2]{book}), which is successful in remarkably many cases.

If $\CC$ is already defined over $\Fq$ one is lead to the so-called trace-zero varieties in $\Jac_{\CC} ( \F_{q^d})$ (see \cite[Section 7.4.2]{book})  and again correspondences induced by covers of curves can be used for attacks on crypto systems based on discrete logarithms on these varieties  by   work of Diem  \cite[22.3.4]{book}. These results already indicate that the use of  Picard groups of curves (e.g. elliptic curves) over non-prime fields $\Fq{_{^d}}$ with $d\geq 4$ is not advisable for cryptographic use.

By more recent work of C. Diem this "feeling" is reinforced for instance for  families of elliptic curves in towers of finite fields. The methods used in these papers use the Weil restriction method explained above only as a "guideline" and  sometimes as tools for proof. The real heart of the methods of Diem is the use of Semaev's summation polynomials.  In this context and in particularly because of suggestions of pairing based cryptography using (supersingular) elliptic curves it is important to mention the  enormous progress made in the computation of discrete logarithm in the multiplicative group of finite non-prime fields \cite{JOP}.

\subsection{Modular Correspondences}\label{Correspondences}
We recall from \cref{sect-5} that  for $N$ prime to $\char (k)$  the modular curve $X_0(N)$ is a regular projective curve, defined over $\Z[1/N]$   and so in particular over $\Q$ and over $\Fp$ with $p$ prime to $N$.   As explained in \cref{sect-5} there is an affine part $Y_0(N)$, which is a (coarse) moduli scheme for the   isomorphism classes of pairs $(E,\eta_N)$ of elliptic curves with cyclic isogeny of degree $N$. This means that for every point $P=(j_E,j_\eta)$ in $Y_0(N)(k)$ there is an elliptic curve $E$ defined over $k$ and an isogeny  $\eta_N:E\ra E'$ with $\ker(\eta_N)$ invariant under the action of $G_k$ and as abelian group isomorphic to $\Z/N$ such that the invariants of $E$ and $E'$ are $(j_E,j_\eta)$.

The points in $X_0(N)\setminus Y_0(N)$ are the cusps, and it is important that these points have a modular interpretation, too. For example, if $N$ is squarefree, then there is one  cusp point at $\infty$ (in the upper half plane) which  corresponds to the pair (N{\'e}ron polygon with $N$ vertices,   $\< \zeta_N \>$) where $\zeta_N$ is a primitive $N$-th root of unity.   

Let $\ell$ be a prime not dividing  $\car(k)\cdot N$. By the  splitting $\Z/\ell\cdot N\cong\Z/\ell\times \Z/N$ and an analogous splitting of the kernel of a cyclic isogeny of degree $\ell\cdot N$ in $C_\ell\times C_N$ we get a natural $k$-morphism  
\[ 
\varphi_\ell: X_0(\ell\cdot N)\ra X_0(N).
\]
Let $\omega_\ell$ be the involution of $X_0(\ell\cdot N)$ induced by the map that sends the pair $(E,\eta)$ with  $\ker(\eta)=C_\ell\times C_N$ to the pair $(E,\eta')$ where the kernel of $\eta'$ is $E[\ell]/C_\ell\times C_N$.    Define 
\[   
\psi_\ell := \varphi_\ell\circ \omega_\ell:  X_0(\ell\cdot N)\ra X_0(N).
\]
We are in the situation described above (with $k=L$) and can define the Hecke correspondence
\[ 
T_\ell: \J_0(N)\ra \J_0(N), 
\]
by  $T_\ell:=\varphi_{\ell*}\circ \psi_{\ell}^*$.

The Hecke ring of $X_0(N)$ is $\mathcal{T}_N= \<T_\ell$   with  $\ell$   prime to $ N \>$, the ring generated by the endomorphisms $T_\ell$.   It is a commutative ring, which is very near to $\End (\J_0(N))$; see  \cite{mazur}.   It acts on the vector space of holomorphic differentials of $X_0(N)$ which can be identified with the $k$-vector space of $k$-rational cusp forms $\Ss_0(k)$ of level $N$ (and trivial nebentype). By classical theory one knows that $\mathcal{T}_N$ is endowed with an Hermitian structure due to the Peterson scalar product, and so the eigenvalues of the operators $T_\ell$ are totally real numbers.

\begin{rem}
 Assume that $\AA$ is a simple factor of $\J_0(N)$. Then $\End^0(\AA)$ contains a totally real field of degree    $\dim(\AA)$.
\end{rem}

This means that factors of $\J_0(N)$ have very special and large endomorphism rings. As consequence there is a splitting of Galois representations of  $G_\Q$ constructed by the action on Tate modules of $\J_0(N)$ into a sum of two-dimensional representations with real eigenvalues, and these "modular representations" play a most important role in number theory, e.g. for the proof of Fermat's Last Theorem. The narrow relation to arithmetic is reflected by the \textbf{Eichler-Shimura} congruence
\[ T_\ell =\mathrm{Frob}_\ell +\ell/\mathrm{Frob}_\ell, \]
where $\mathrm{Frob}$ is the Frobenius endomorphism on $\mathcal{J}_0(N) \bigotimes \mathbb{F}_\ell$.    In particular, $\mathrm{Frob}_\ell$ satisfies the \textbf{Eichler-Shimura} equation
\[
X^2-T_\ell\cdot X+\ell =0.
\]
A curve $\CC$ whose Jacobian is a factor of $\mathcal{J}_0(N)$ is called \textbf{modular of level $N$}.

Using cusps forms it is possible to determine its period matrix, decide whether it is hyperelliptic, and then compute its   Weierstrass equation (see \cite{Weng}, \cite{Web}). The \emph{importance for cryptography} is the fact that $\mathrm{Frob}$ satisfies the quadratic Eichler-Shimura equation over a totally real number field, and this    can be used for point counting for curves of genus $\geq 2$ as in \cite{Gaudry}. Hence \emph{modular curves of genus $2$ are potentially usable for $DL$-systems}.

\subsection{Correspondences via Monodromy Groups}
We assume that we have a cover morphism 
\[ f:\CC  \ra  \P^1\]
defined over $k$ of degree $n$, satisfying some fixed ramification  conditions and having a fixed monodromy group $G_f:=\mbox{Mon } (f)$.    We have morphisms
\[ \tilde{f} : {\tilde{\mathcal{H}}}\stackrel{h}{\ra} \CC\stackrel{f}{\ra}\P^1 \]
with $\tilde{f}$ a Galois cover of $f$ with Galois group $G_f$.   For simplicity, we assume that the field of constants of ${\tilde{\mathcal{H}}}$ is $k$.  This setting is motivated by the theory of \emph{Hurwitz spaces} and it is hoped that one can exploit their rich and, over $\C$, well understood theory (\cite{FK1} and \cite{FK2}). 

Next we choose  subgroups $H_1\subset G_f$ fixing $\CC$ and $H_2$ containing $H_1$.   Let $\mathcal{H}$ be the curve fixed by $H_1$ and $\mathcal{D}$  the fixed curve under $H_2$. So $\mathcal{H}$ covers both $\CC$ and $\mathcal{D}$.  Let 
\[ h : \mathcal{H}\ra \CC  \; \; \text{ and } \; \; g:\mathcal{H}\ra \mathcal{D} \]
with morphisms induced by the Galois action.   Hence the degree of $h$ is equal to $\deg (h) = \frac{|G_f|}{|H_1|\cdot n}$ and the degree of $g$ is equal to   $\deg(g)=\frac{|H_2|}{|H_1|}$.   We get a correspondence
\[ \eta: \J_\CC\ra \J_\mathcal{D} \]
by applying $g_*\circ h^*$ to the Picard groups. In general, $\eta$ will  be neither injective nor surjective. 

\begin{lem}
Assume   that $\J_\mathcal{D}$ is a simple abelian variety, $\dim \J_\mathcal{D} = g_{\CC}$,    and 
that
  there is a prime divisor  $\p_\infty$ of $\CC$ which is totally ramified under $h$, i.e. there is exactly one prime divisor $\PP_\infty$ of $\mathcal{H}$ with norm $\p$,  and that there is no non-constant morphism of degree $\leq \deg(h)$ from $\mathcal{D}$ to the projective line.
 Then  $\eta$ is an isogeny. 
\end{lem}

\proof Since $J_{\mathcal{D}}$ is simple, it is enough to show that $\eta$ is not the zero map.  Let $\p'_\infty $ be the norm of $\PP_\infty$ under $g$. Without loss of generality we can assume that $k$ is algebraically closed.  So we find a prime divisor $\PP$ of $\mathcal{H}$ which is different from all prime divisors in  $g^{-1}(\p'_\infty)$.

Let $c$ be the class of $\p-\p_\infty$, where $\p = h_*(\PP)$. Then $\eta(c)$ is the class of the divisor 
\[ D_\p:= \sum_{\PP\in h^{-1}(\p)}g_*(\PP) - \deg(h)\cdot g_*(\PP_\infty).\]
Note that $D_\p \neq 0$ (as  divisor). If the class of $D_\p$ would be trivial, then there would be a non-constant function on $\mathcal{D}$ with pole order $\leq\deg(h)$ and hence a non-constant map of $\mathcal{D}$ to the projective line of degree $\leq \deg(h)$, which is a contradiction. 

\qed

We shall see in \cref{gen-3} that we can realize the situation (over $\bar{k}$) of the lemma for hyperelliptic curves of genus $3$ with non-decomposable Jacobian, $f$ a polynomial of degree $6$, $G_f=S_4$, $H_1$ a subgroup of order $2$ and $H_2$ a subgroup of order $6$. This leads to  isogenies of degree $8$ discussed by B. Smith, and generically maps hyperelliptic curves to non-hyperelliptic curves.
The \emph{importance for cryptography} is that generic hyperelliptic curves of genus $3$ are not usable for DL-systems.

It is an open and challenging problem to find  other interesting correspondences  of low degree between Jacobian varieties induced by correspondences between curves and (possibly) attached to Hurwitz spaces.

\section{Genus 3 curves and cryptography}\label{gen-3}
\begin{quest}
Can one use curves $\CC$ of genus $3$ for DL-systems?
\end{quest}

To find equations for random curves of genus $3$ is easy: Either take a regular plane quartic (non-hyperelliptic curve) or a curve with  equation $Y^2=f(X)$ with $\deg(f)=7$. 
In both cases the addition law is easily implemented and fast. If $\CC$ is hyperelliptic, the Cantor algorithm is well-studied and fast, moreover one can transform it into formulas (involving, alas, many special cases), which  are sometimes more convenient for implementations near to specialized hardware. The generic cases for addition and doubling are explicitly given by Algorithms 14.52 and 14.53 in \cite{book}. The timings  are not too far away from additions in groups in elliptic curves (with comparable size order); see \cite[Table 14.13]{book}.
For non-hyperelliptic curves see \cite{Oyono-1}. 
But we have already discussed the security problem: One can only use hyperelliptic curves of genus $3$ for which the Jacobians do not possess an easily computed isogeny to another principally polarized abelian variety which has a non-hyperelliptic polarization. As we shall see next this will endanger the DL in  "generic" hyperelliptic curves of genus $3$.

\subsection{Isogenies via $S_4$-Covers}\label{trigonal}
As observed by B. Smith \cite{Smith}  "many" hyperelliptic curves are isogenous to  non-hyperelliptic curves via an isogeny  with degree dividing $8$. This fact is interpreted in terms of Hurwitz spaces and connected modular spaces in  \cite{FK1, FK2}.  We refer for details and refinements to these papers.

For our purposes it will be enough to look at the case that the base field $k$ is algebraically closed, which we shall assume from now on. For applications in cryptography one has to study rationality problems;  see \cite{Smith} and \cite{FK2}. The construction relies on the so-called trigonal construction of Donagi-Livn\'{e} \cite{Donagi}.
We begin with a hyperelliptic curve $\CC$ of genus $3$ and its uniquely determined hyperelliptic projection  $f_1:  \CC \ra \P^1$   with $8$ ramification points $P_1,\cdots, P_8$, which extend to the Weierstra{\ss} points of $\CC$.  By linear algebra we show that there is a map
\[ f_2:\P^1 \ra \P^1\]
of degree $3$  with the following properties:
\begin{itemize}
\item $f_2$ is unramified in $P_1,\cdots P_8$, its ramification points are denoted by $Q_1,\cdots Q_4$ on the base line  $\P^1$.    The ramification order in $Q_i$ is $2$, and so each $Q_i$ has exactly one unramified extension under $f_2$ denoted by $Q'_i$.

\item $f_2 (\{P_1,\cdots P_8\} ) =\{S_1,\cdots S_4\}$ such that, after a suitable numeration, $f_2(P_i)=f_2(P_{4+i})$ for $1\leq i\leq 4$. 
\end{itemize}
Now we can use Galois theory.
\subsubsection{The monodromy group of $f_2$} Obviously,  the Galois closure $\tilde{f_2}=f_2\circ h_2$ of $f_2$ has as Galois group the symmetric group $S_3$ (since $f_2$ is not Galois because of the ramification type), and $h_2$ is  a  degree $2$ cover $\mathcal{\E'}\stackrel{h_2}{\ra} \CC$.   From Galois theory we get that  $ \tilde{f_2}= \pi\circ \eta $,  where 
\[ \eta:\mathcal{E}'\ra \E \]
is a cyclic cover of degree $3$ with Galois group equal to the alternating subgroup $A_3$.   Then, $\E$ is a quadratic cover of $\P^1$ ramified exactly at the discriminant 
\[ \Delta_1=Q_1 + \cdots + Q_4\]
of $f_2$. Therefore $\E$ is an elliptic curve with cover map $\pi$ to $\P^1$.  From construction and Abhyankar's lemma it  follows that $\eta$ is unramified. Hence $\E'$ is an elliptic curve, too, and $\eta$ is an isogeny of degree $3$ (after applying a suitable translation). 

\subsubsection{The monodromy group of $f=f_2\circ f_1$}
Since $f$ is a cover of degree $6$,  its Galois group can be embedded into $S_6$.  But a closer analysis using the specific ramification situation shows; see  \cite[Thm.~3]{FK1}.

\begin{lem} 
The monodromy group of $f$ is isomorphic to $S_4$.
\end{lem}

Let  $\tilde{f} : \tilde{\CC} \ra \P^1$ 
be the Galois cover of curves factoring over $f$ with Galois group $S_4$.  Let $\CC'$ be the subcover of $\tilde{\CC}$ with function field equal to the composite of the function fields of $\CC$ and $\E'$, i.e. the normalization of the fiber product of $\CC$ with $\E'$.   Let 
\[ \pi_\CC : \CC' \ra \CC \]
the projection to $\CC$, which is a cover of degree $2$.   The Galois group of $\tilde{\CC}/\CC$ contains two transpositions. Let $\sigma$ be one of them, chosen such that with $G_2= \< \sigma \>$ we get  $\CC':=\tilde{\CC}/G_2$.  Hence, $\sigma$ is contained in precisely two of the stabilizers $T_1, \ldots, T_4$ of the elements $\{1,2,3,4\}$ on which $S_4$ acts. 
Let 
\[ \pi_T : \tilde{\CC} \ra \mathcal{D}:= \tilde{\CC}/T \]
be the quotient map. Then $\tilde{f}$ factors over $\pi_T$ as $\tilde{f} = g\circ\pi_T$, where $g: \mathcal{D} \ra\P^1$ has $\deg(g) = 4$. Note that $g$ is primitive
(does not factor over a quadratic subcover).   We can use the Hurwitz genus formula to compute the genus of $\mathcal{D}$.  For this we have to determine the ramification of $\mathcal{D}/\P^1$ under $g$.  

\begin{lem}  
The genus of $\mathcal{D}$ is equal to $3$, and so is equal to the genus of $\CC$. 
\end{lem}

We are interested in the case that $\J(\CC)$ is simple. Then we get from \cref{sect-correspondences}  that:

\begin{prop}
Let $\J_\CC$ be a simple abelian variety and  $\mathcal{D}$ be non-hyperelliptic.  The pair of cover maps $(\pi_\CC,\pi_T)$ from $\CC'$ to $(\CC,\mathcal{D})$ induces an isogeny
\[ \eta:\J_\CC\ra\J_\mathcal{D},\]
whose kernel is elementary-abelian and has degree $\leq 8$.
\end{prop} 
 
A more detailed analysis due to E. Kani shows that the proposition is true without the assumption that $\mathcal{D}$ is non-hyperelliptic.  Thus we have the following:
 
\begin{cor} 
The notations are as above. Let $k$ be equal to $\Fq$ and assume that $\mathcal{D}$ is non-hyperelliptic.  Then the computation of the Discrete Logarithm in $\Pic^0_\CC$ has complexity $\O(q)$.
\end{cor}

This result motivates the question whether the assumptions of  the Corollary are often satisfied.   Empirically, B. Smith has given a positive answer.  A rigorous answer is given in \cite{FK2}. 
We have already explained that by the construction of a $(2,3)$-cover as above we have found a  generically finite and dominant morphism from a Hurwitz space $\mathcal{H_1}$ to the hyperelliptic locus in the moduli space $\M_3$ of curves of genus $3$. Hence $\mathcal{H_1}$ is a scheme of dimension $5$.
Via the trigonal construction we have, to each hyperelliptic curve $\CC$, found a curve $\mathcal{D}$ of genus $3$  with a cover map 
\[ g:\mathcal{D}\ra \P^1\]
 with $\deg(g)=4$ and the monodromy group of $g$ equal to $S_4$. Moreover, a detailed study of the construction allows to determine the ramification type of $g$ in the generic case: 

There are $8$ ramification points of $g$, exactly $4$ points $P_1,\dots,P_4$ amongst them are  of type $(2,2)$ (i.e. $g^*(P_i)=2(Q_{i,1}+Q_{i,2})$, and the other $4$ ramification points are of type $(2,1,1)$.   Hence $(\mathcal{D},g)$ yields a point in a Hurwitz space $\mathcal{H}_2$ of dimension $5$.

In \cite{FK2} one discusses the hyperelliptic locus $\mathcal{H}_{hyp}$ in $\mathcal{H}_2 $. The computational part of this discussion determines conditions for the coefficients of Weierstra{\ss} equations for curves $\mathcal{D}$ lying in  $\mathcal{H}_{hyp}$. This is rather complicated, but one sees that generically these coefficients are parametrized by a $4$-dimensional space.  Rather deep and involved geometric methods have to be used to transfer these computations into scheme-theoretical results
and to get

\begin{thm}
The Hurwitz space $\mathcal{H}_{hyp}$ is a unirational, irreducible variety of dimension $4$, provided that $\char k > 5$. Moreover, the natural forgetful map
\[
\mu: \mathcal{H}_{hyp} \ra \M_3
\]
to the moduli space $\M_3$ of genus $3$ curves has finite fibers and so its image is also irreducible of dimension $4$.
\end{thm}

\begin{cor}
Assume that $k$ is algebraically closed.  There is a 1-codimensional subscheme $U$ of $\M_{3,hyp}$ such that for $\CC\notin U$ the isogeny $\eta$  maps $\J_\CC$ to the Jacobian of a non-hyperelliptic curve $\mathcal{D}$.
\end{cor}

Replacing the algebraically closed field $k$ by a finite field $\Fq$ one has to study rationality conditions for $\eta$. This is done in \cite{Smith} and \cite{FK2}. As result we get the following:

\begin{cor}
There are $\O(q^5)$ isomorphism classes of hyperelliptic curves of genus $3$ defined over $\Fq$ for which the discrete logarithm in the divisor class group of degree $0$ has complexity ${\O}(q)$, up to log-factors.  Since $|\Pic^0(C)|\sim q^3$, the DL system of these hyperelliptic curves of genus $3$ is weak.
\end{cor}

\subsection{Point Counting}
In general, not much is known about fast point counting algorithms on curves of genus $3$ (aside of the general fact that for all abelian varieties there is a polynomial time algorithm due to Pila).  But as we have seen above,  for applications in cryptography we have to restrict ourselves to special hyperelliptic curves (where it is not at all clear what "special" means for a concrete curve), and so we do not lose much by restricting to  hyperelliptic curves $\CC$ whose Jacobian $\J_\CC=:\J$ has a special endomorphism ring $\O_\J$.

\subsubsection{Real Multiplication}

A first possibility is to assume that $\J$ has real multiplication. This means that $\O_\J$ contains an order $\mathcal{R}$ of a totally real field of degree $3$.   An immediate consequence is that there are many isogenies at hand, and in the case of genus $2$ this situation has accelerated the point counting dramatically \cite{GKS}.  So there is hope that the same could happen for Jacobians of dimension $3$.
Hence
 it is interesting to construct hyperelliptic curves $\CC$ such that $\J_\CC$ has real multiplication. In view of the results about Jacobians of the modular curves $X_0(N)$ in \cref{Correspondences} it is natural to look for curves whose Jacobian is a quotient of $J_0(N)$ for some $N$.  This was successfully done by H. J. Weber \cite{Web}. The procedure is as follows: 

First, one computes eigenspaces of dimension $3$ of the space of cusp forms of level $N$ under the Hecke operators. Using the attached differentials one can compute  (over $\C$) the period matrix of the corresponding factor $\J$ of $J_0(N)$ and decides whether it is principally polarized and hence is the Jacobian of a curve $\CC$. Using theta-null values one decides whether $\CC$ is hyperelliptic.

If so, one can compute invariants of the curve, and (e.g by  using Rosenhain models compute a Weierstrass equation (at the end over $\Z[1/N]$) of $\CC$  whose Jacobian is  a simple factor of $\J_0(N)$. Reduction modulo $p$ gives hyperelliptic curves over $\Z/p$ of genus $3$ with (known) real multiplication.

The method works quite well but has one disadvantage: Since there are many non-hyperelliptic curves of genus $3$ with real multiplication we are not sure whether the constructed  curve is isogenous to a non-hyperelliptic curve under the trigonal construction described above.

\subsubsection{Complex Multiplication}

We strengthen the condition on $\End (\J)$ and assume that $\J$ has complex multiplication and 
so is the reduction of a curve 
 defined over a number field. Recall that this means that there is an embedding of $\End(\J)$ as order $\O$  into a CM-field $K$, i.e. $K$ is a totally imaginary quadratic extension of a totally real field $K_0$ of degree $3$ over $\Q$.

The arithmetic of $\CC$ and $\J$ is reflected by the arithmetic of orders in $K$. In particular, one finds the Frobenius endomorphism of reductions of $\CC$ modulo prime ideals $\p$ of $K$ as element in $\O$. This solves the problem of point counting on $\CC$ modulo $\p$ immediately. Moreover, class field theory of $K$ gives both a classification of isomorphism classes of curves $\CC$ with  CM-field $K$ and methods to find period matrices of $\J$ and so equations of $\CC$.   Details and more references can be found in \cite{book} sections 5.1 and 18.3.

But trying to find examples for hyperelliptic curves attached to CM fields of degree $6$ one runs into trouble since these examples seem to be very rare. 
(Recall that the hyperelliptic locus in $\M_3$ has codimension 1.)
This was one of the results of the thesis of A. Weng (Essen 2001). So one has to use some force: If $\J$ has an automorphism of order $4$ the curve $\CC$ has an automorphism of order at least $2$ 
and if $\J$ is simple the quotient of $\CC$ by this automorphism has to be $\P^1$, and so $\CC$ is hyperelliptic and has an   automorphism of order $4$.
The existence of $J$ with automorphism $\varphi$ of order $4$ is obtained by a special choice of the CM-field $K$:

Let $K_0$ be a totally real field of degree $3$ with class number $1$ (there are many fields with these properties) and take $K=K_0(\sqrt{-1})$, and for $\O$ take the maximal order of $K$. In \cite{Weng} one finds in detail how these choices lead to many examples of hyperelliptic curves over finite fields suitable for cryptography.

%
%

There is a bit of hope that the following question may have an affirmative answer.

\begin{quest}
We assume now that $\CC$ is a hyperelliptic curve with an automorphism   of order $4$.  Is $\CC$   resistant against the trigonal attack?
\end{quest}

If the answer would be yes and since automorphisms of degree $4$ survive under isogenies of degree prime to $2$  one could hope 
to have a positive answer to the following. 
   
\begin{quest}
Let $\CC$ be a hyperelliptic curve with an automorphism of order $4$ and with simple Jacobian variety $\J$.   Let 
\[ \eta:  \J   \ra \J' \]
be an isogeny with $\J'$ principally polarized.   Is $\J'$  the Jacobian variety of a hyperelliptic curve?
\end{quest}

In the case of a positive answer to the question the CM-curves with CM-field $K$ of degree $4$ containing $\sqrt{-1}$  would deliver  nice and easy to handle candidates for cryptographically usable DL-systems

\section{Genus 2 curves and cryptography}\label{genus-2}
Curves $\CC$ of genus $2$ with at least one  rational Weierstra{"s} point $P_\infty$ are very interesting objects for creating DL-systems and in most aspects they can very well compete with elliptic curves (see \cref{ell-curves}). The research area around these curves is attractive since there is a lot of activity but also a lot of unsolved problems till now.

\paragraph{\textbf{Security}} The hardness of the DL in the Picard groups of randomly chosen curves over prime fields of order $q$ is comparable with the hardness on elliptic curves over prime fields of order $q^2$, in particular, all known versions of index-calculus attacks have a complexity equal to $\O(q)$ and hence are not more efficient than generic algorithms. (Recall that because of  
\cref{Hasse-Weil} we can expect that $|\Pic^0_\CC|\sim q^2$.)\\

\textbf{Addition:}   By our assumption we can assume that $\CC$ is given by a Weierstra{"s} equation 
\[ Y^2Z^3=f(X,Z)\]
 with $f(X,Z)$ homogenous and of degree $5$ in $X$.

Hence we can use Mumford representations of reduced divisors and the Cantor algorithm (see \cref{Cantor}).  A detailed analysis including all special cases is done in \cite[Section 14.3.2]{book}, including a determination of complexity (see Table 14.2 and Table 14.13).
 
Alternatively we  use the interpolation formulas given explicitly in \cref{gen-2-add}, and we have the choice to chose coordinates taylor-made to soft- and hardware environments.  As result we can state that the efficiency of group operations in $\Pic^0_\CC$ is on the same level as it is for elliptic curves.
 
If we are only interested in scalar multiplication (e.g. for key exchange) we can use, as in the case of elliptic curves, a "Montgomery ladder" to compute this multiplication. The role of $x$-coordinates of points on elliptic curves is played by coordinates on the Kummer surface related to the Abelian surface $\J_\CC$.\\
 
\textbf{Kummer Surfaces} The following  is due to P. Gaudry \cite{G2}  and Gaudry/Lubicz  \cite{GL}.  We embed $\CC$ into $J_\CC$ by using $P_\infty$ as base point, and continue the  hyperelliptic involution $\omega$ of $\CC$ to $J_{\J_\CC}$. Then $J_\CC/w=:K$ is the Kummer variety of $\CC$, and we have an
embedding of $\mathbb{P}^1\cong \CC/w$ into $K$.

On $K$  the action of $\mathbb{Z}$ is induced by the group structure on the Jacobian. One checks that one has a scalar multiplication but no group structure (compare the case of elliptic curves). Hence the usual add-and double algorithm to get a fast
scalar multiplication does not work. To repair this one uses the \textbf{Montgomery  ladder} (see \cite{book} ) which is  well known for
elliptic curves.

To make the ladder very fast one uses a remarkable tool: classical modular forms in an abstract setting!

More concretely,  P. Gaudry uses in \cite{G2} classical theory of theta functions, their p-adic interpretation and reduction, exploits "classical" doubling formulas and gets \textbf{extremely simple doubling formulas}.

One drawback is that  the model used for $\mathbb{\CC}$ based on theta functions has   bad reduction modulo 2. So  in \cite{G2} Gaudry had to exclude the important case that the ground field has even characteristic. More arithmetic geometry, namely  the theory of minimal models enabled him together with D. Lubicz to remove this restriction \cite{GL}.

The third necessary aspect important for the construction of DL-systems is point counting, which has to be  so effective that in reasonable time one finds by a random search curves $\CC$ and fields $\F_p$ with the property that a large prime number $\ell$ with size $\sim q^2$ divides
$|\J_\CC(\F_p)|$.

\subsection{Point counting on curves of genus $2$}
\subsubsection{Point counting on random curves}\label{random}
 A generic method to determine the order of a finite group is given by a variant of Shank's baby-step giant-step method, whose efficiency depends on the size of the interval in which on can place $|\J_\CC(\F_p)|$. This is used if one knows $|\J_\CC(\F_p)|$ modulo a rather big number $N$ together with the information given by the Hasse estimate.
 
To get such a congruence one tries to determine the characteristic  polynomial of the Frobenius endomorphism $\phi_p$ modulo "enough" small numbers by 
 its action on torsion points. So a first step for counting algorithms is to determine polynomials or ideals which vanish on torsion points of a given order. This procedure was already the key part of Schoof's algorithm for elliptic curves (\cref{schoof}).
 
But for really fast algorithms for elliptic curves one needs  one more ingredient: isogenies and corresponding modular polynomials respectively ideals.
 
 For curves of genus $2$ division polynomials and modular polynomials are  not so well understood as for elliptic curves but as we have announced already in \cref{gen-2-add} this is  an active area of research. 
 
 The starting point is the Mumford representation of points on $\J_\CC(\F_p)$.
We assume that $\CC$ is given by  $y^2=f(x)$ as before and $D = \< u, v\>$ is a reduced divisor.  Most reduced divisors have weight 2, i.e the degree of $u$ is $2$. The set of those divisors with strictly lower weight is called $\Theta$. A divisor of weight 1 i.e., with a single point $P = (x_P , y_P )$, is represented by 
\begin{equation}\label{weight-1}
 \<  u(x), v(x) \> = \<    x - x_P , y_P  \>. 
\end{equation} 
The unique divisor of weight 0, is the identity $\O$ given as 
$ \O = \< u(x), v(x) \>  = \< 1, 0 \>$. 
Any divisor of weight 2 is given as
\[ \< u(x), v(x) \> = \< x^2+ u_1 x +u_2, v_0 x +v_1 \>. \]
The following algorithm using \textbf{division polynomials} is due to Gaudry and  Harley.

For a divisor of weight 1, i.e. given by an ideal $P= \< x- x_P, y_P\>$ in general position we have 
\[ [l]  P = \left\<   x^2 + \frac { d_1^{(l)} (x_p)} {d_0^{(l)} (x_p)  } x + \frac {d_2^{(l)} (x_p) } {d_0^{(l)} (x_p) }, 
y_P \left( \frac { e_1^{(l)} (x_p)} {e_0^{(l)} (x_p) } x + \frac { e_2^{(l)} (x_p)} { e_0^{(l)} (x_p)}      \right)   \right\> 
\]
where $e_i^{(l)}, d_i^{(l)}$ are polynomials with degrees $\deg d_i^{(l)}= 2 l^2 - (i+1)$ and $\deg e_i^{(l)}= 3 l^2 - (i+1)$, for $i=0, 1, 2$; see \cite{GH} for details.   

For a divisor of weight two, we consider it as a sum of two divisors of weight 1, say $D=P_1 + P_2$ where $P_1 = \< x-x_1, y_1\>$ and $P_2 = \< x-y_2, y_2\>$, where $x_1$ and $x_2$ are roots of $u(x)$ and $y_i = v(x_i)$ and $u,v$ come from the Mumford presentation. Then, 
\[ [l] D = [l] P_1 + [l]P_2.\] 
With these formulas one computes the order of  $\Jac_\CC{\F_p} $ modulo $l$; see \cite[Section~5.4]{GH}. The cost of the algorithm is 
\[ O(l^2) M(l^2) + O (d \log q) M(l^4) + O (l^2+ \log q ) M(d), \]
where $M(n)$ is the number of field operations required to multiply two polynomials of degree $n$,  and $d$ is the smallest degree of resultants of $u(x)$ and $ v(x)$; see \cite{GH}.

Next one tries to determine "modular equations" for finite subschemes of $\J_\CC$. Here a paper of P.Gaudry and E. Schost  is a remarkable beginning
(see \cite{Gaudry-2005}.  Mixing the results and methods and using many tricks  Gaudry and Schost succeed in  \cite{GS} to count points on the Jacobian of some hundreds of random curves of genus $2$ and finally found one having the security level of AES 128.  The development still goes further. The interested reader should have a look at the paper  \cite{AG}.


\subsection{Modular curves of genus two}
In \cref{random}  we have seen that in principle one can find cryptographically relevant curves of genus $2$  by a search on random curves using an analogue of Schoof's algorithm. But the necessary input of computing capacity (and implementation art) is rather heavy.

So it may useful to look for classes of special curves for which point counting is easier.
Again work of Gaudry, together with David Kohel and Benjamin Smith shows that one can accelerate the algorithm dramatically if the Jacobian of $\CC$ 
has real multiplication (see  \cite{GKS}.

The  new algorithm has, for large $p$, complexity $\widetilde{\O}(\log^5 q)$. This is used to compute a 256-bit prime-order Jacobian, suitable for cryptographic applications, and also the order of a 1024-bit Jacobian.

 Hence it is interesting to construct curves of genus $2$ with real multiplication.
 We describe how this can be done by using factors of Jacobians of modular curves. 
 The basic reference is  Wang \cite{wang}.

 Let $N$ be a positive integer and $X_0 (N)$ the modular curve as described in \cref{mod-curves}.  
Let $S_2 (N)$ be the space of cusp forms of weight 2 for $\Gamma_0 (N)$ and $f= \sum_{n=1}^\infty a_n q^n \in S_2 (N)$ be a newform. This newform determines a simple abelian variety $\AA_f$ which is a factor of $J_0 (N):= \Jac X_0 (N)$. 
The knowledge of the newform is equivalent with the knowledge of holomorphic differentials on the factor and this is used by 
Wang to compute the period matrix of $\AA_f$ by computing the complex integrals of a symplectic basis.  Moreover, he determined conditions when this period matrix $\Omega_f$ corresponds to a principally polarized abelian variety. He did this for factors of dimension $\geq 2$  but we  focus on the case when this variety has dimension two. 

Having found the period matrix one has to construct the curve.  We have mentioned this task already in before and cited \cite{Web}. 

%
%

\subsubsection{Thetanulls}  
Once the period matrix $\Omega_f$ is determined, one can compute the theta functions.   For any genus 2 curve we have six odd theta characteristics and ten even theta characteristics. The following are the sixteen theta characteristics, where the first ten are even and the last six are odd. For simplicity, we denote them by $\T_i = \ch{a} {b}$ instead of $\T_i \ch{a} {b} (z , \t)$ where $i=1,\dots ,10$ for the
even theta functions.
\begin{small}
\[
\begin{split}
\T_1 = \chr {0}{0}{0}{0} , \,  \T_2 = \chr {0}{0}{\frac{1}{2}} {\frac{1}{2}} ,  \T_3 =\chr
{0}{0}{\frac{1}{2}}{0} , \, \, \T_4 = \chr {0}{0}{0}{\frac{1}{2}} , \, \,    \T_5 = \chr{\frac{1}{2}}{0}
{0}{0} ,\\
  \T_6  = \chr {\frac{1}{2}}{0}{0}{\frac{1}{2}} , \, \,
  \T_7 = \chr{0}{\frac{1}{2}} {0}{0} , \, \,
  \T_8 = \chr{\frac{1}{2}}{\frac{1}{2}} {0}{0} , \, \,
  \T_9 = \chr{0}{\frac{1}{2}} {\frac{1}{2}}{0} , \, \,
  \T_{10} = \chr{\frac{1}{2}}{\frac{1}{2}} {\frac{1}{2}}{\frac{1}{2}} ,\\
\end{split}
\]
\end{small}
and the odd theta functions  correspond to the following characteristics
\[   \chr{0}{\frac{1}{2}} {0}{\frac{1}{2}} , \,
   \chr{0}{\frac{1}{2}} {\frac{1}{2}}{\frac{1}{2}} , \,
    \chr{\frac{1}{2}}{0} {\frac{1}{2}}{0} , \, \,
    \chr{\frac{1}{2}}{\frac{1}{2}} {\frac{1}{2}}{0} , \,
    \chr{\frac{1}{2}}{0} {\frac{1}{2}}{\frac{1}{2}} , \,
    \chr{\frac{1}{2}}{\frac{1}{2}} {0}{\frac{1}{2}}  \]

We call   fundamental theta constants $\T_1, \, \T_2, \, \T_3, \, \T_4$.  All the other theta constants can be expressed in terms of these four; see \cite{wijesiri} for details.   The classical result of Picard determines   the relation between theta characteristics and branch points of a genus two curve.
\begin{lem}[Picard] Let $\X$ be a genus 2 curve be given by affine equation
\begin{equation} \label{Rosen2}
Y^2=X(X-1)(X-\lambda)(X-\mu)(X-\nu).
\end{equation} Then, $\lambda, \mu, \nu$   can be written as follows:
%
\begin{equation}\label{Picard}
\l = \frac{\T_1^2\T_3^2}{\T_2^2\T_4^2}, \quad \mu = \frac{\T_3^2\T_8^2}{\T_4^2\T_{10}^2}, \quad \nu = \frac{\T_1^2\T_8^2}{\T_2^2\T_{10}^2}.
\end{equation}
\end{lem}

Such branch points were expressed in terms of  the fundamental theta constants. 
%
\begin{lem}[\cite{wijesiri}]  \label{possibleCurve}
Every genus two curve $\X$ can be written in the form:
\[
y^2 = x \, (x-1) \, \left(x - \frac {\T_1^2 \T_3^2} {\T_2^2  \T_4^2}\right)\, \left(x^2 \, -   \frac{\T_2^2 \, \T_3^2 +
\T_1^2 \, \T_4^2} { \T_2^2 \, \T_4^2} \cdot    \a  \, x + \frac {\T_1^2 \T_3^2} {\T_2^2 \T_4^2} \, \a^2 \right),
\]
where $\a = \frac {\T_8^2} {\T_{10}^2}$ and in terms of $\, \, \T_1, \dots , \T_4$ is given by
\[
  \a^2 + \frac {\T_1^4 + \T_2^4 - \T_3^4 - \T_4^4}{\T_1^2 \T_2^2 - \T_3^2 \T_4^2 } \, \a + 1 =0
\]
Furthermore,  if $\alpha = {\pm} 1$ then $V_4 \embd \Aut(\X)$.
\end{lem}

From the above we have that $\T_8^4=\T_{10}^4$ implies that $V_4 \embd \Aut(\X)$.  in \cite[Lemma~15]{wijesiri} 
determines a necessary and equivalent statement when $V_4 \embd \Aut(\X)$ in terms of thetanulls. 

The last part of the lemma above shows that if $\T_8^4=\T_{10}^4$ then all coefficients of the genus 2
curve are given as rational functions of powers of these 4 fundamental theta functions. Such fundamental theta functions
determine the field of moduli of the given curve. Hence, the curve is defined over its field of moduli.

%

Once the fundamental thetanulls are computed, the arithmetic invariants $J_2$, $J_4$, $J_6$, $J_{10}$ can be computed via formulas given in \cite{Ig-2}.  

Till now the computations were made over $\R$ with large enough precision. Now we use the $\Z$-structure of $X_0(N)$ to identify the invariants with integers. Reducing modulo primes one finds invariants of curves defined over finite fields.
%
%
 Then by \cite{univ-g-2} we compute an equation of the genus 2 curves over a minimal field of definition for these invariants and so over  a finite field.

\begin{rem}
In \cite{wang} there are used absolute invariants instead of the above arithmetic invariants. Moreover, the case of curves with automorphism group of order $> 2$ doesn't seem to have been considered. Nevertheless, it seems as this was not a problem for $N\leq 200$, which seems to suggest that no such  genus two curves appear for such $N$. 
\end{rem}

\subsection{CM curves}
We further specialize and want to use curves of genus $2$ whose Jacobian has complex multiplication. We shall use class field theory and the theory of Taniyama-Shimura of CM-fields to find such curves over number fields.By reduction we find curves with CM over finite fields, and again class field theory of CM-fields reduces point counting modulo $p$ to the computation of the trace of an element in the CM-field with norm $p$.
 
Choose a squarefree integer $d \in \N$ such that $K_0 :=\Q (\sqrt{d})$ has class number one.  Let $\alpha= a + b\sqrt{d}$ be squarefree and $\alpha > 0$. Thus $K= K_0 (i \sqrt{\alpha})$ is a CM field of degree 4.  We choose $d$ and $\alpha$ such that $K/\Q$ is Galois with group $V_4$ (i.e. Klein four-group). Since $[K:\Q]=4$ and $K$ is CM field we have four distinct embeddings $\varphi_i$, $i=1, \dots , 4$ of $K$ into $\C$. A tuple $(K, \Phi) = (K, \{ \varphi_1, \varphi_2 \})$ is called CM-type.  For an ideal $ I \subset \O_K$ we define
\[ \Phi (I)= \{ (\varphi_1 (x), \varphi_2 (x))^t, x \in I \}. \]
Then  $\C^2/\Phi(I)$ is an Abelian variety with complex multiplication by $\O_K$. Conversely every abelian variety $\AA$ of CM-type $(K, \Phi)$ with complex multiplication by $\O_K$ is isomorphic to an abelian variety $\AA_{I, \Phi}$; see Shimura-Taniyama (1961) \cite{ShTa}. 

The period matrix of $\AA_{I, \Phi}$ lies in the Siegel upper half plane $\mathbb H_2$ and therefore we can equip $\AA_{I, \Phi}$ with a principal polarization determined by an element $\gamma \in K$.      



\subsubsection{Class polynomials}

For elliptic curves with complex multiplication by $\O_K$ the $j$-invariant lies in the Hilbert class field of the imaginary quadratic field $K$. The case of $g=1$  is simpler  due to the fact that the \textit{reflex CM-field} $\hat K$ is equal to $K$ (see \cite{Shi98}), which is not true for higher genus. 
The following is mostly due to A. Weng \cite{Weng-2}. 

\begin{thm}      
Let $K$ be a CM-field such that $[K:\Q]=4$.  

\begin{itemize}

\item[i)] For every genus 2 curve $\CC$ with CM-type by $\O_K$, the absolute invariants $\x_1, \x_2, \x_3$ are algebraic numbers that lie in a class field over the reflex CM-field $\hat K$. 

\item[ii)] For two genus 2 curves $\CC$ and $\CC^\prime$ with CM with $\O_K$ we have that $\x_i (\CC)$ and $\x_i (\CC^\prime)$, for $i= 1, 2, 3$,  are Galois conjugates.

\item[iii)] Let $\{ \CC_1, \CC_2, \dots , \CC_s\}$ be a set of representatives of isomorphism classes of genus 2 curves whose Jacobians have CM with endomorphism ring $\O_K$.  Denote by $\x_i^{(j)}$, the $i$-th absolute invariant of $\CC_j$. The polynomials
\[ H_{K, i} (X) := \prod_{j=1}^s (X- \x_i^{(j)}), \]
for $i=1, 2, 3$, have coefficients in $\Q$.  
\end{itemize}
\end{thm}

Polynomials $H_{K, 1}$, $H_{K, 2}$, $H_{K, 3}$, are called the \textbf{class polynomials}. 


\begin{thm}   
Let $K$ be a CM-field of degree 4 and $p\geq 7$ a prime which does not divide the denominators of the class polynomials $H_i (X):=H_{K, i} (X)$, $i=1, 2, 3$.  Then the following hold:

\begin{itemize}

\item For all $w\in \O_K$ with $w \bar w = p$, $H_i (X)$ have a linear factor over $\F_p$ corresponding to $w$.

\item For each $\alpha \in \F_p$ there are two $\F_p$-isomorphism classes $\AA_{p,1}$ and $\AA_{p,2}$ of principally polarized abelian varieties over $\F_p$ with absolute invariants $\x_i = \alpha$, for $i=1, 2, 3$. 

\item The principally polarized abelian varieties  $\AA_{p,1}$ and $\AA_{p,2}$ have CM by $\O_K$. 

\item The number of $\F_p$-rational points of $\AA_{p, j}$, $j=1, 2$, is given by 
\[ \prod_{r=1}^4 \left(   1 + (-1)^j w_r \right) \]

\item The equation $w \bar w=p$ for $w\in \O_K$ has (up to conjugacy and sign) at most two different solutions. Hence, for every $CM$-field of degree 4 there are at most four different possible orders of groups of $\F_p$-rational points of principally polarized abelian varieties defined over $\F_p$ with CM by $\O_K$. 
\end{itemize}

\end{thm}

Once we compute the class polynomials $H_{K, i}$ we can reduce them module $p$  (for large enough $p$) and get $H_{K, i} (X) \mod p$.  The roots of $H_{K, i} (X) \mod p$ are the absolute invariants of genus 2 curves $\CC$ modulo $p$. Now that we know invariants of the curve we can determine its equation as in \cite{univ-g-2}.  Then the reduced curve is defined over $\F_p$ or a quadratic extension. 

For example, if we are in the first case of the above theorem, say we find elements $w_1, \bar w_1 \in \O_K$ such that $w_1 \bar w_1 = p$ then there exists at most one second solution (up to conjugation) such that $w_2 \bar w_2 = p$. We set $W:=\{ \pm w_1, \pm w_2\}$.   Then the order of $\Jac (\CC \mod p)$, over $\F_p$, is $\{ \chi_w (1) \, | \, w \in W\}$, where $\chi_w (T)$ is the characteristic polynomial of $w$. 

By using the CM-method for curves of genus $2$ one gets a very efficient way to construct cryptographically strong DL-systems as extensive tables as shown in  the thesis of A. Weng.

One could hope to use these results not only for DL-systems but also for isogeny graphs of Jacobians of dimension $2$ and it could be worthwhile to investigate whether they could be used in the way we shall see in  \cref{ell-curves} below.

\section{Elliptic curve cryptography}\label{ell-curves}

Finally we come to the most interesting and well-understood case.   We shall use isogenies between elliptic curves and their computation quite often, and so we begin, for the convenience of the reader, 
with a fundamental result of J. V{\'{e}lu}; see \cite{Vel}.

\begin{prop}[V\'{e}lu's formula]  Let $E$ be an elliptic curve,  defined over a field $k$, with equation
\[E : y^2=x^3+ax+b\]
and  $G\subset E(\bar{k})$ be a finite subgroup invariant under $G_k$. The separable isogeny $\phi : E\to E/G$, of kernel $G$, can be written as follows: 

For any $P(x, y)$ we get
\begin{equation}
\phi(P) = \left( x + \sum_{Q\in G\setminus\{\O\}}x(P+Q)-x(Q), \;      y + \sum_{Q\in G\setminus\{\O\}}y(P+Q)-y(Q)
    \right)
\end{equation} 
and the curve $E/G$ has equation $y^2=x^3+a'x+b'$, where
\begin{align*}
    a' &= a - 5\sum_{Q\in G\setminus\{\O\}}(3x(Q)^2+a),\\
    b' &= b - 7\sum_{Q\in G\setminus\{\O\}}(5x(Q)^3+3ax(Q)+b).
\end{align*}
 
\end{prop}

Thus, knowing a finite subgroup $G$ of $E$ we can explicitly construct the corresponding isogeny $E\to E^\prime:=E/G$.  

\subsection{Endomorphism ring of $E$}
The following results are mostly due to \textbf{M. Deuring} and mainly contained in the beautiful paper \cite{Deu}.

\begin{defi} 
Let $\E$ be an elliptic curve over $k$. $\E$ is ordinary if and only if $\End(\E)$ is commutative. $\E$ is supersingular if and only if $\End(\E)$ is not commutative.
%
\end{defi}

\begin{thm}[Deuring]\label{lift}
Let $\E$ be an elliptic curve defined over a field $k$. The following hold:
\begin{enumerate}
\item[i)] If  $\car (k) =0$,    then $\E$ is ordinary and
\begin{itemize}
\item $\End_{\ov{k}}(\E)=\Z$ (generic case) or {$\End_{\ov{k}}(\E)$ is an order $O_\E\subset \Q(\sqrt{-d_\E})$,\, $d_\E>0$ (CM-case)}.

\item  Take $\E$ with $CM$ with order $O_\E$.  Let $\mathcal{S}_\E$ be the set of $\C $-isomorphism classes of elliptic curves with endomorphism ring $O_\E$. 
Then $\Pic (O_\E)$ acts in a natural and simply transitive  way on $\cS_\E$, hence $\cS_\E$ is a principally homogeneous  space with translation group $\Pic(O_\E)$: 
For $c\in \Pic(O_\E)$, $\frak A\in c$ and $\C/O_\E=\E_0$ we get  $c\cdot [\E_0] $ is the class of $\C/\frak A$.

\end{itemize}

\item[ii)] (\textbf{Deuring's Lifting Theorem})\label{lift}
Let $\E$ be an elliptic curve over $\Fq$ which is ordinary over $\ov{\Fq}$. Then there is, up to $\C$-isomorphisms, exactly one 
elliptic curve $\tilde{\E}$ with CM over a number field $K$ such that
\begin{itemize}
\item there is a prime $\frak p$ of $K$ with $\tilde{\E}_\frak p \cong \E$ with $\tilde{\E}_\frak p$ the reduction of $\tilde{\E}$ modulo $\frak{p}$, and

\item $\End (\tilde{\E}) =\End(\E)_\frak p = O_\E$,  with $O_\E $ an order in an imaginary quadratic field.
\end{itemize}

\item[iii)] If $\E$ is  supersingular,  then 
\begin{itemize}
\item up to twists, all supersingular elliptic curves in characteristic $p$ are defined over $\mathbb F_{p^2}$, i.e. their $j$-invariant lies in $\mathbb F_{p^2}$.

\item  $|\E(\mathbb F_{p^2})|=(p\pm1)^2$, and the sign depends on the twist class of $\E$.  

\item $\End_{\ov{\mathbb{ F}_p}}(\E)$ is a maximal order in the quaternion algebra $\Q_p$, which is unramified outside of $\infty$ and $p$.
\end{itemize}

\end{enumerate}
\end{thm}
%
We remark that the endomorphism ring of an elliptic curve over a finite field $\Fq$ is never equal to $\Z$ since there is the  Frobenius endomorphism $\phi_{\Fq,\E}$ induced by the Frobenius automorphism of $\Fq$ which has degree $q$. We give a first application of the lifting theorem.

\begin{cor}[Hasse]
Let $\E$ be an ordinary elliptic curve over $\Fq$. Then the Frobenius endomorphism $\phi_{\Fq,\E}$ is an integer in an imaginary quadratic field with norm $q$, and hence has a minimal polynomial
\[\chi_{\E,q}(T) = T^2 - \tr  (\phi_{\Fq,\E}) \cdot T + q \]
 with 
\[ | (\tr (\phi_{\Fq,\E})^2 -4q |< 0.\]
\end{cor}

Recall that the number of $\Fq$-rational points of $\E$  is
\[|\E(\Fq)|=:n_{\Fq,\E}=\chi_{\E,q}(1).\]

\begin{cor}\label{hasse}
$| n_{\Fq,\E}-q-1|< 2\sqrt q.$
\end{cor}

Using the result iii) in \cref{lift} and the observation that the eigenvalues of $\phi_{\mathbb{F}{q^d},\E}$ are the $d$-th power of the eigenvalues of $\phi_{\Fq,\E}$ we get that 
\[ | n_{\Fq,\E}-q-1|\leq 2\sqrt q \]
for all elliptic curves of $\Fq$.  This is the  \emph{Hasse bound} for  elliptic curves, a special case of the  Hasse-Weil bound for  points on Jacobian varieties  over finite fields  (\cref{Hasse-Weil}).
%
\subsection{Point Counting}\label{schoof}

Corollary~\ref{hasse} is the key fact for a polynomial time algorithm for computing the order of $\E (\Fq)$ for elliptic curves $\E$ defined over the field $\Fq$, which is called  \textbf{Schoof's Algorithm}.

The idea is to compute $\chi_{\E,q}(T)\mod n$ for small numbers $n$ by computing the action of $\phi_{\Fq,\E}$ on $\E[n]$ (take for instance $n=\ell$ as small prime number or $n=2^k$ with $k$ small) and then to use CRT and the Hasse bound for  the trace of $\phi_{\Fq,\E}$ to determine $\chi_{\E,q}(T)$. To do this use  the classical n-division polynomials $\Psi _n$. The disadvantage is that $\deg (\Psi_n) \sim n^2/2$ and therefore the Schoof algorithm is too {slow}.

The way out of this problem is to use \'{e}tale isogenies with cyclic kernel of order $n$ and the fact (see \cref{sect-5}) that we can interpret these isogenies with the help of points on an explicitly known curve, namely the modular curve $X_0(n)$. An explicit equation for an affine model of $X_0(N)$ is given by the classical modular polynomial $\phi ( j, j_N)$. It allows an effective computation of isogenies (as functions including the determination of the image curve) at least if $n$ is of moderate size).

\begin{thm}[V\'{e}lu, Couveignes, Lercier, Elkies, Kohel, and many other contributors:]
The cost for the computation of an isogeny of degree $\ell$ of an elliptic curve $\E$ over $\Fq$ is  
\[\O(\ell^2+\ell\log(\ell)\log(q)).\]
\end{thm}

The \textbf{Idea of Atkin-Elkies} is: Use  \emph{\'{e}tale isogenies} of  small degree of  $\E$ instead of points, and use the modular polynomial $\phi_n$ of degree $\sim n$. The resulting \emph{Schoof-Atkin-Elkies algorithm} is very fast, in particular if one assumes  as "standard conjecture" the generalized Riemann hypothesis (GRH). 

\begin{cor}[SAE]  
$|\E(\Fq)|$
can be computed (probabilistically, with GRH) with complexity $\O((\log q)^4)$.  Therefore we can construct, for primes $p$ sufficiently large, (many) elliptic curves with  $|\E(\F_p)|=k\cdot \ell$  with $k$ small (e.g. $k=1$ if we want) and $\ell$ a prime so large that  (using classical computers and according to our best knowledge) the security level of the discrete logarithm in $\E(\F_p)$ is matching  AES 128 (or larger).
\end{cor}

\subsection{Looking for Post-Quantum Security}
As we have seen in  \cref{schoof} we can construct elliptic curves over prime fields such that the resulting  DL-systems are secure under the known attacks. But the situation changes totally if we allow algorithms based on quantum computers.  We shall discuss now how we can use isogenies of elliptic curves  (and maybe, of curves of larger genus with convenient endomorphism rings) to find key exchange schemes staying in the frame of Diffie-Hellman type protocols as described in \cref{Diffie-Hellman}.

\subsubsection{The isogeny graph}

\begin{defi}
Let $\E$ be an elliptic curve over $\Fq$.  The  \textbf{isogeny graph} of $\E$ is a graph where nodes are $j$-invariants of elliptic curves isogenous to $\E$ and edges are  isogenies between the curves attached to the nodes.
\end{defi}

For some applications one restricts the degree of the isogenies defining edges.  Isogeny graphs of ordinary elliptic curves are discussed by using the lifting theorem of Deuring, identifying these graphs with graphs coming from ideal classes of orders in imaginary quadratic fields , and then using analytic number theory and properties of modular forms. The study of graphs of supersingular elliptic curves uses properties of maximal orders of quaternion algebras. In both cases one gets the following result due to D. Jao, S. D. Miller and R. Vekatesan \cite[Prop.~2.1]{DJP}.

\begin{thm}
The isogeny graph of $\E$ is a Ramanujan graph.
\end{thm}
 
A first application, also to be found in \cite{DJP} is the following:

\begin{corollary}
Assume the $\E$ and  $\E'$ are elliptic curves over $\Fq$ with $\End(\E)=\End(\E')=\O$. Then there exists a subexponential algorithm which relates  the DLP in $\E(\Fq)$ to the DLP in $\E'(\Fq)$.
\end{corollary}

Hence we cannot weaken the DLP  by applying isogenies between elliptic curves with the same order.

\subsubsection{The System of Couveignes-Stolbunov}
A second application of isogeny graphs is constructive.    We sketch in the following work of Stolbunov \cite{stol} and Couveignes \cite{Cou}.  We use the results of \cref{lift} for an ordinary elliptic curve $\E_0$ over $\Fq$ with ring of endomorphism $\End(\E_0)=O$, which is an order in a quadratic imaginary field.

In analogy to the notation Theorem~\ref{lift} define  $\cS_{E_0}$ as   set of isomorphism classes of elliptic curves over $\ov{\F}_q$ with ring of endomorphisms  $O$. Then $\cS_{E_0}$ is a $\Pic (O)$-set.  Hence, we can use it for \textit{Key Exchange protocols}  as in \cref{Diffie-Hellman}:   \textit{The partner $P$ chooses  $c\in \Pic (O)$ and publishes the $j$-invariant of $c \cdot E_0$.}

The exchange is not as fast as DL-systems since we cannot use a \textit{double-and add}-algorithm but it is feasible since one finds enough isogenies that are composites of isogenies of small degree (smoothness); for an example see \cite{stol}.   The  \textbf{security } depends on the hardness of the following problem:

\begin{problem}
Find an isogeny between two given isogenous elliptic curves.
\end{problem}

The following gives an idea of the running time for the solution to this problem. 

\begin{prop}[Kohel, Galbraith, Hess, Smart et al.]
The expected number of \textbf{bit}-operations for the computation of an isogeny between ordinary elliptic curves over $\Fq$   with endomorphism ring   $O_{K_E}$ is
\[ \O(q^{1/4 + o(1) } \log^2(q) \log\log (q)).\]
\end{prop}
    
But recall: We are in the situation where an abelian group is acting on a set, and so there is a subexponential algorithm to solve the hidden-shift problem.    This means that we can only expect \emph{subexponential} security in the Q-bit world for the key exchange scheme; see results of    Childs, Jao, Soukharev in   \cite{MR3163097}.     Comparing this with the situation we have nowadays with respect to the widely tolerated  RSA-system this may be not so disastrous. \\

\subsubsection{The Key Exchange System of De Feo}\label{supersingular}
The suggestion is to use supersingular elliptic curves over $\F_{p^2}$ and their properties also stated in \cref{lift}.  Take 
\[  p = r^a\cdot s^b\cdot f - 1 \]
with  $p\equiv 1 \mod 4$.    Then 
\[ \E_0 :  \;  Y^2 Z=X^3+X Z^2 \]
is a supersingular elliptic curve over $\F_{p^2}$. We describe the key exchange scheme invented and implemented by  De Feo, Jao and Pl{\^ut} \cite{DJP} in the frame we have introduced in \cref{Diffie-Hellman}.  

The objects in the  categories $\CC_i$ \, $(i=1,2)$  have as  \textbf {objects} the isomorphism classes of supersingular curves $\E$ over $\F_{p^2}$ isogenous to $\E_0$ and hence with  
\[ | \E(\mathbb{F}_{p^2})|=(r^a\cdot s^b\cdot f)^2.\]
The \textbf{morphisms in $\CC_1$} are isogenies $\varphi$ with $|\ker(\varphi)|$  dividing $r^a$.  The \textbf{morphisms in $\CC_2$} are isogenies $\psi$ with  $|\ker(\psi)|$  dividing $s^b$. For these categories  {pushouts} exist.   For additional information  choose  $P_1, P_2$ of order $r^a$  and $Q_1,Q_2$ of order $s^b$  in $\E_0(\F_{p^2})$.   \\

\newpage
 
\noindent \textbf{Key Exchange:} 

\begin{itemize}
\item The Partner  $\mathcal{P}_1$
 chooses $n_1,n_2\in \Z/r^a$ and the isogeny 
\[ \eta : \E_0 \ra \E_0/  \< n_1 P_1 + n_2P_2 \> =: \E_1.\]
\item  
%
The Partner $\mathcal{P}_2$
 chooses $m_1, m_2\in \Z/s^b$ and computes the isogeny 
\[ \psi: \E_0 \ra \E_0/ \< m_1 Q_1 + m_2 Q_2 \> =: \E_2. \]
\item  
%
$\mathcal{P}_2$ sends  $(\E_2, \psi(P_1),\psi(P_2))$. 
 
\item  
$\mathcal{P}_1$ can compute the common secret, the pushout of $\eta$ and $\psi$ as 
\[ 
\E_3:= \E_2/ \< n_1 \psi (P_1) + n_2 \psi (P_2) \>. \\
\]
\end{itemize}

An analogous procedure enables $\mathcal{P}_2$ to compute the  isomorphy class of $\E_3$, which is the common secret of the partners.

\medskip

Again \textbf{security} depends on the hardness to compute an isogeny of two elliptic curves, but now the two elliptic curves are supersingular. \\
 
\noindent \textbf{State of the art}: The best known algorithms have exponential complexity {$p^{1/4}$ (bit-computer) resp. $p^{1/6}$ (quantum computer)}, and so one can hope that a prime $p$ with $768$ bit yields AES128 security level. So we have, compared with other post-quantum suggestions for key exchange schemes, a very small key size.

In contrast to the ordinary case the  groups around like the class groups of left ideals in maximal orders \textbf{are not abelian}, and so the hidden shift problem is not solved till now in subexponential time.

\bibliographystyle{amsalpha} 
\bibliography{ref}{}

\end{document}